\documentclass[preprint,10pt,Shah2021]{elsarticle}

\newcommand{\michele}[2][cyan]{\emph{\textcolor{#1}{#2}}}
\newcommand{\nirav}[2][red]{\emph{\textcolor{#1}{#2}}}

\newcommand{\rev}[2][black]{\textcolor{#1}{#2}} % Rev 1 red 
\newcommand{\revv}[2][black]{\textcolor{#1}{#2}} % Rev 2 blue
\newcommand{\revvv}[2][black]{\textcolor{#1}{#2}} % Rev 3 green
\usepackage{amssymb}
\usepackage[spanish,english]{babel}
\usepackage{type1cm}        % activate if the above 3 fonts are
\usepackage{makeidx}         % allows index generation
\usepackage{graphicx}        % standard LaTeX graphics tool
\usepackage{wrapfig}         % when including figure files
\usepackage{multicol}        % used for the two-column index
\usepackage[bottom]{footmisc}% places footnotes at page bottom
\usepackage{gensymb}
\usepackage{float}
\usepackage{amsmath}
\usepackage{newtxtext}       %
\usepackage{newtxmath}       % selects Times Roman as basic font
\usepackage{url}
\usepackage{bm}
\usepackage{tikz}
\usepackage{subcaption}
\usepackage{neuralnetwork}
\usepackage{comment}
\usepackage{tabularx}
\usepackage[margin=2.5cm]{geometry}
\usepackage{multirow}
\usepackage{mwe}
\usepackage{pgfplots}
\usetikzlibrary{babel}
\usetikzlibrary{patterns}
\usetikzlibrary{arrows, arrows.meta}
\usetikzlibrary{shapes.geometric}
\DeclareMathOperator{\Tr}{Tr}
\DeclareMathOperator{\spn}{span}

\DeclareMathOperator*{\argmin}{arg\,min}

\makeindex             % used for the subject index
\journal{Finite Elements in Analysis and Design}

\begin{document}

\begin{frontmatter}

%% Title, authors and addresses

%% use the tnoteref command within \title for footnotes;
%% use the tnotetext command for theassociated footnote;
%% use the fnref command within \author or \affiliation for footnotes;
%% use the fntext command for theassociated footnote;
%% use the corref command within \author for corresponding author footnotes;
%% use the cortext command for theassociated footnote;
%% use the ead command for the email address,
%% and the form \ead[url] for the home page:
%% \title{Title\tnoteref{label1}}
%% \tnotetext[label1]{}
%% \author{Name\corref{cor1}\fnref{label2}}
%% \ead{email address}
%% \ead[url]{home page}
%% \fntext[label2]{}
%% \cortext[cor1]{}
%% \affiliation{organization={},
%%            addressline={}, 
%%            city={},
%%            postcode={}, 
%%            state={},
%%            country={}}
%% \fntext[label3]{}

\title{Finite element based model order reduction for parametrized one-way coupled steady state linear thermo-mechanical problems}

%% use optional labels to link authors explicitly to addresses:
%% \author[label1,label2]{}
%% \affiliation[label1]{organization={},
%%             addressline={},
%%             city={},
%%             postcode={},
%%             state={},
%%             country={}}
%%
%% \affiliation[label2]{organization={},
%%             addressline={},
%%             city={},
%%             postcode={},
%%             state={},
%%             country={}}

\author[label1]{Nirav Vasant Shah}
\ead{shah.nirav@sissa.it}

\author[label1]{Michele Girfoglio}
\ead{michele.girfoglio@sissa.it}

\author[label2,label3]{Peregrina Quintela}
\ead{peregrina.quintela@usc.es}

\author[label1]{Gianluigi Rozza \corref{cor1}}
\ead{gianluigi.rozza@sissa.it}
\cortext[cor1]{Corresponding author}

\author[label4]{Alejandro Lengomin}
\ead{alejandro.lengomin@arcelormittal.com}

\author[label5]{Francesco Ballarin}
\ead{francesco.ballarin@unicatt.it}

\author[label2,label3]{Patricia Barral}
\ead{patricia.barral@usc.es}

\affiliation[label1]{organization={Scuola Internazionale Superiore di Studi Avanzati (SISSA)},%Department and Organization
            addressline={via Bonomea 265}, 
            city={Trieste},
            postcode={34136}, 
            %state={},
            country={Italy}}

\affiliation[label2]{organization={Instituto Tecnol\'oxico de Matem\'atica Industrial (ITMATI), currently integrated in CITMAGA},%Department and Organization
            addressline={s/n, Campus Vida, R\'ua de Constantino Candeira}, 
            city={Santiago de Compostela},
            postcode={15705}, 
            %state={},
            country={Spain}}            

\affiliation[label3]{organization={Departamento de Matem\'atica Aplicada, Universidade de Santiago de Compostela},
            %addressline={}, 
            city={Santiago de Compostela},
            postcode={15782}, 
            %state={},
            country={Spain}}            
            
\affiliation[label4]{organization={Primary \& By-Products Department, ArcelorMittal, Global R\&D Asturias},%Department and Organization
            addressline={P.O. Box 90}, 
            city={Avil\'es},
            postcode={33400}, 
            %state={Asturias},
            country={Spain}}

\affiliation[label5]{organization={Department of Mathematics and Physics, Catholic University of the Sacred Heart},%Department and Organization
            addressline={via Musei 41}, 
            city={Brescia},
            postcode={25121}, 
            %state={},
            country={Italy}}
            
\begin{abstract}
This contribution focuses on the development of Model Order Reduction (MOR) for one-way coupled steady state linear thermo-mechanical problems in a finite element setting. %Based on the nature of the mathematical framework, we use a segregated approach at both full and reduced order level: the temperature field, computed by solving the energy conservation equation, is inserted into the momentum conservation equation to compute the thermal stresses.%We are interested in investigating the steady state condition of the system.
%The approach is applied for the modeling of phenomena arising in blast furnace hearth walls. %The temperature profile, computed by solving the energy conservation equation, is inserted into the momentum conservation equation to compute the thermal stresses.
We apply Proper Orthogonal Decomposition (POD) for the computation of reduced basis space. On the other hand, for the evaluation of the modal coefficients, we use two different methodologies: the one based on the Galerkin projection (G) and the other one based on Artificial Neural Network (ANN). We aim to compare POD-G and POD-ANN in terms of relevant features including errors and computational efficiency. In this context, both physical and geometrical parametrization are considered. We also carry out a validation of the Full Order Model (FOM) based on customized benchmarks in order to provide a complete computational pipeline. The framework proposed is applied to a relevant industrial problem related to the investigation of thermo-mechanical phenomena arising in blast furnace hearth walls. \\%It is shown that POD-ANN framework allows to decouple errors estimates related to energy and momentum e Numerical experiments show that the accuracy of the POD-ANN method and show the substantial speed-up enabled at the online stage as compared to a traditional RB strategy.%The experimental results include error analysis, speedup analysis, eigenvalue decay and comparison of solutions at a given parameter.\\
\end{abstract}


\end{frontmatter}

%% \linenumbers

%% main text
\section{Introduction}\label{introduction_chapter}

Due to technological developments occurred in recent years, the high-fidelity numerical computations, based on the so-called Full Order Models (FOM) (e.g., finite element or finite volume methods), are required to be performed for many configurations. This puts the computational resources under considerable stress. In this context, Model Order Reduction (MOR) %also known as Reduced Basis (RB) method,
has been introduced as an efficient tool to accelerate the computations with ``affordable'' and ``controllable'' loss of accuracy. The faster computations obtained by MOR helped in many query contexts, e.g. quick transfer of computational results to industrial problems. %In particular, MOR in combination with geometric parametrization has emerged as an alternative to the shape optimization and has been used in many engineering applications.

The basic idea on which MOR is based is related to the fact that often the parametric dependence of the problem at hand has an intrinsic dimension much lower than the number of degrees of freedom associated to the governing FOM. The development of a MOR consists in two main steps. The first one is the so-called \emph{offline} stage when a database of several high-fidelity solutions is collected by solving the FOM for different values of physical and/or geometrical parameters. Then all the solutions are combined and compressed to extract a set of basis functions that approximate the low-dimensional manifold on which the solution lies. %build the space onto which we can project the solution manifold and efficiently compute the solutions for new parameters istances (\emph{online} stage).
The second step is the \emph{online} stage when the information obtained in the \emph{offline} stage is used to efficiently compute the solutions for new parameters instances. For a comprehensive review on MOR, the reader is referred to, e.g., \cite{mor_bader,mor_survey_benner,mor_survey_benner2,mor_book_benner3, haasdonk_chapter,RBniCS,mor_book_quarteroni,schilders}. %[29, 47, 8, 7, 3, 9].

In this work, we address the development of a MOR framework in a Finite Element (FE) environment for one-way coupled steady state linear thermo-mechanical problems. %The approach is applied for the modeling of phenomena arising in blast furnace hearth walls.
%\michele{(We should be sure of citing all the works available in the literature. VERY IMPORTANT. Please, add also the paper published on Computational Mechanics)} \nirav{Do you mean \cite{variational_multiscale}?If yes, already added in this paragraph.}
In the literature different MOR techniques have been proposed in the context of thermo-mechanical problems. Gu\'erin et al. \cite{rational_craig_halle} developed Rational Craig-Hale methodology for the investigation of thermo-mechanical coupling effects in turbomachinery. Benner et al. \cite{thermomechanical_pod_bt_mm_comparison} compared the performance of Proper Orthogonal Decomposition (POD), Balanced Truncation, Pad\'e approach and iterative rational Krylov algorithm for the approximation of the transient thermal field concerning an optimal sensor placement problem for a thermo-elastic solid body model. Zhang et al. \cite{variational_multiscale} introduced reduced order variational multiscale enrichment method and tested their approach on proper benchmark tests related to the thermo-mechanical loading applied to a 2D composite beam and a functionally graded composite beam. More recently, Hern\'andez-Beccero et al. \cite{kms_paper} used Krylov Modal Subspace method %\michele{(sorry Nirav, I am not expert...this method is different from iterative rational Krylov algorithm already mentioned when talking about of Benner et al.?)} \nirav{(Benner's paper considered Hermite interpolant of Transfer function whereas Beccero's paper considered Krylov subspace for steady state and eigenvectors for dynamic behavior)} 
for thermo-mechanical models as applicable to machine tools. \revvv{We highlight that all these works are focused on MOR for the efficient reconstruction of the time evolution of the thermo-mechanical field. Regarding steady state problems, such as those dealt with in this work, Hoang et al. \cite{greedy_paper} used a two-field reduced basis algorithm based on the greedy algorithm in a physical parametrization setting. However, here we consider not only physical parameters but also geometrical ones. }Concerning the methodology adopted, we use POD for the construction of reduced basis space. On the other hand, for the computation of the reduced coefficients, we consider two different approaches: a standard Galerkin projection (G) and Artificial Neural Network (ANN). POD-G aims to generate a MOR by a projection of the governing equations onto the POD space. So, the reduced coefficients associated to the POD bases are obtained by solving a set of algebraic equations. On the other hand, POD-ANN belongs to the category of machine learning methods based on systems that learn from data. In this framework, the reduced coefficients are computed by the employment of a properly trained neural network. There is a broad range of strategies for the development of new deep learning architectures %for numerical analysis \cite{gkn_paper,fno_paper,pinn_paper1,pinn_paper2} and specifically, for
for non-intrusive MOR approaches, i.e. without the need to access to FOM implementation: see, e.g., \cite{pod_nn_paper,ann_mor_transient_flow,non_intrusive_1,prnn_paper, federico_paper,nicola_demo_paper,wang_ann_paper,amsallem_paper,demo_strazzullo_paper,meneghetti_paper}. \revvv{We compare POD-G and POD-ANN in terms of relevant computational features including errors and speed up. From this viewpoint, the present contribution draws inspiration by \cite{pod_nn_paper} in which it has been shown that POD-ANN performs better than POD-G both in terms of efficiency and accuracy for steady state problems in heat transfer and fluid dynamics. We underline that whilst in \cite{pod_nn_paper} nonlinear problems are addressed, here we deal with a linear modeling framework by obtaining significantly different results. Our approach is applied within an industrial framework related to the investigation of thermo-mechanical phenomena arising in blast furnace hearth walls. Neural networks can in theory represent any functional relationship between inputs and output. However, many applications are still unexplored, and we retain that the use of POD-ANN method in the context of thermo-mechanical problems of industrial interest could open the door towards the application of deep learning techniques to new multiphysics scenarios.}% and we consider this work as an intermediate step towards the development of MOR for nonlinear thermo-mechanical problems. %We highlight that, to the best of our knowledge, geometrical parametrization has not yet been investigated within a thermo-mechanical modeling context %\nirav{Recently, application of deep learning methods to partial differential equations has seen significant development (\cite{gkn_paper,fno_paper,pinn_paper1,pinn_paper2,prnn_paper}).} TODO questi riferimenti li mettiamo nel paragrafo di ANN

The workflow of this paper is organized as follows. In Sec. \ref{chapter_physicalproblem} the physical problem related to the blast furnace hearth walls is introduced. The full order model in strong formulation is described in Sec. \ref{chapter_thermo-mechanical_model}.  %the relevant function spaces are introduced and weak formulation is derived.
The corresponding weak formulation is then derived in Sec. \ref{weak_form_chapter} and the finite element analysis is introduced in Sec. \ref{fem_chapter}. %The validation of the full order model is presented in Sec. \ref{chapter_benchmark_test}. 
Subsequently, in Sec. \ref{parametric_rb_chapter}, the MOR approach is described, and results obtained are shown and discussed. Conclusions and perspectives are drawn in Sec. \ref{remark_chapter}. \revv{Finally, we dedicate a wide Appendix to the FOM validation in order to show the functionality of the developed finite element code}. %provide a complete computational pipeline to be used for efficient and accurate numerical simulation of thermo-mechanical problems %We use Proper Orthogonal Decomposition (POD) for the construction of reduced basis space. Concerning the computation of the parameter dependent coefficients of the reduced basis, we apply both Galerkin projection-based intrusive and Artificial Neural Network (ANN)-based non intrusive techniques. The accuracy of such strategies is compared in terms of relevant computational features: errors, speed-up and eigenvalue decay.

\section{Physical problem}\label{chapter_physicalproblem}
\begin{figure}
\begin{center}
\includegraphics[width=0.3\linewidth]{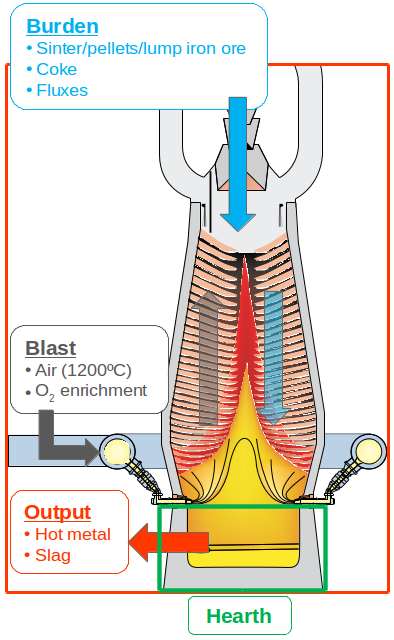}
\end{center}
\caption{Blast furnace [Courtesy: ArcelorMittal, Spain].}
\label{blast_furnace_schematic}
\end{figure}

%\michele{DA QUI!!!}
Steelmaking is a very old process that has contributed to the development of technological societies since ancient times. The previous stage is the ironmaking process, which is performed inside a blast furnace, whose general layout is shown in Figure \ref{blast_furnace_schematic}. It is a metallurgical reactor used to produce hot metal from iron ore. %(95 \% Fe, 4.5\% C \michele{Nirav, check the percentages, the last 0.5\% seems missed}) from iron ore.
For further details the reader is referred, e.g., to \cite{blast_furnace_book1,blast_furnace_book2}. %\michele{(Nirav, here I would cite one or two works at which the reader could refer in order to acquire more information about the process)}. \nirav{Done}%The burden contains iron ore, fluxes and coke. The process involves exothermic and endothermic reactions.

The blast furnace operates at a high temperature (up to 1500 \degree C). The associated thermal stresses significantly limit the overall blast furnace campaign period.
%\begin{wrapfigure}{r}{5cm}
%\end{wrapfigure}
%The numerical computation of thermal stresses requires coupling of thermal and mechanical models.
In this context, thermo-mechanical modeling has been used extensively either to support experimental campaign or to design various components. V\'azquez-Fern\'andez et al. \cite{sara_paper} simulated the stationary heat transfer in a trough of a blast furnace to ensure durability based on the location of critical isotherm. In Barral et al \cite{patricia_1,patricia_2}, the transient behaviour of the temperature during a tapping and a full campaign cycle was analyzed. Numerical modeling of heat flows in the blast furnace hearth lining was used by Swartling M. et al. \cite{Swartling} for improving experimental assessment. Thermo-mechanical modeling of blast furnace hearth was also developed by Brulin et al. \cite{Brulin2}: they used  micro-macro approach with homogenization method for replacing bricks and mortars by an equivalent material. %Notice that this review is far from complete. %\michele{Can we add some other references?}%\michele{Also here, we should be sure of citing all the important references}% and solve the associated coupled system.

%Here, we adopt a finite element method to address the problem. %that has successfully been applied in several applications involving thermo-mechanical modeling, such as flowform technology \cite{flowform_tech}, stress analysis for VLSI devices \cite{VLSI}, and additive manufacturing \cite{additive_manufacturing}. %For what concerning previous works related to the blast furnace hearth,
The blast furnace operates under different conditions, each of which is governed by a different mathematical model. Considering the objectives of the present work, the following simplifications are considered:

\begin{itemize}

\item Taphole operation is not part of this study. The perforation action of the taphole and the important pressures in the draining of the hot metal and slag produce important mechanical stresses located in the area that requires a deeper analysis and that is out of the scope of this work. %So, in this document, the detail of the tap on the furnace walls is omitted.

\item Since the objective is to be able to calculate in real time the effects of wall design on blast furnace operation, we focus on the steady state operations.

\item We assume that the hearth is made up of a single homogeneous, elastic, and isotropic material with temperature-independent material properties.

\item Heat transfer only by conduction within hearth walls will be considered. The temperature of the molten metal inside the hearth is assumed constant and known. Therefore, the fluid region will not be part of the problem.

\end{itemize}

\section{Full order model}\label{chapter_thermo-mechanical_model}

In this section, we discuss the mathematical formulation corresponding to the physical problem described in Sec. \ref{chapter_physicalproblem}. We present the thermo-mechanical model in cylindrical coordinates endowed with suitable boundary conditions in Sec. \ref{chapter_thermo-mechanical_model_2}. Then, we introduce the axisymmetric hypothesis and derive the axisymmetric thermo-mechanical model in Sec. \ref{chapter_thermo-mechanical_model_3}.

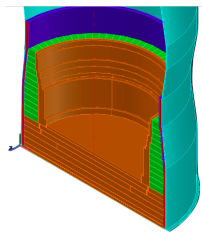
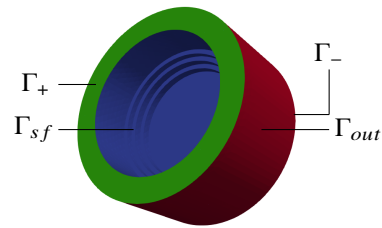
\begin{figure}[H]
    \centering
    \begin{subfigure}[t]{0.5\textwidth}
        \centering
        \includegraphics[height=1.2in]{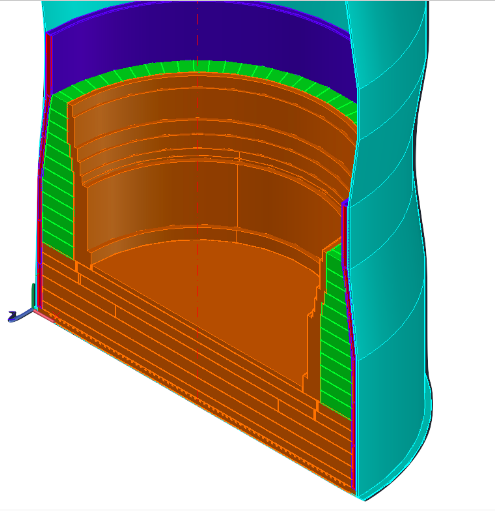}
        \caption{Section of hearth geometry [Courtesy : ArcelorMittal].}
    \end{subfigure}%
%     \begin{subfigure}[t]{0.45\textwidth}
%      \centering
%      \includegraphics[height=1.2in]{images/hearth_3d_computational_domain.png}
%      \caption{Simplified domain $\Omega$}
%      \end{subfigure}
\begin{subfigure}[t]{0.45\textwidth}
\centering
\begin{tikzpicture}
\node at (-0.5,-0.5) {\includegraphics[height=1.2in]{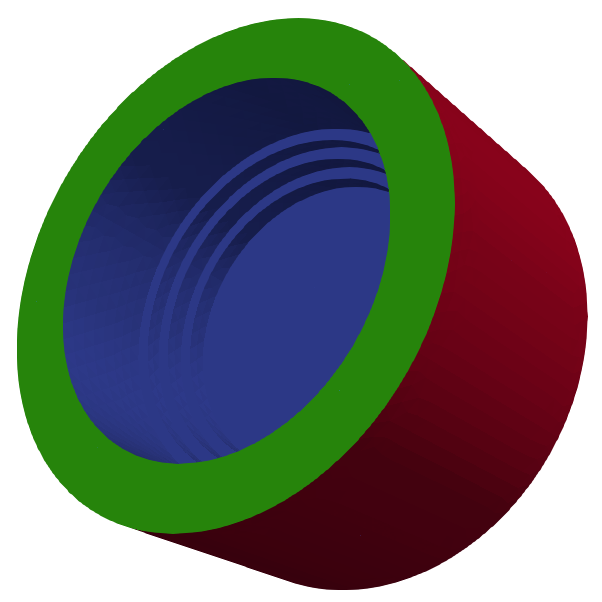}};
\node at (-2.5,-0.7) {$\Gamma_{sf}$};
\draw[-] (-1.2,-0.7)--(-2.2,-0.7);
\node at (-2.5,-0.1) {$\Gamma_+$};
\draw[-] (-1.7,-0.1)--(-2.2,-0.1);
\node at (1.8,-0.7) {$\Gamma_{out}$};
\draw[-] (0.5,-0.7)--(1.4,-0.7);
\node at (1.4,0.3) {$\Gamma_-$};
\draw[-] (0.93,-0.5)--(1.4,-0.5);
\draw[-] (1.4,-0.5)--(1.4,0.1);
\end{tikzpicture}
\caption{\revv{Simplified domain $\Omega$ and its boundaries.}}
\end{subfigure}
% \begin{minipage}[b]{0.4\linewidth}
%     \centering
%     \includegraphics[width=.5\linewidth]{images/gamma_plus.png}
%     \subcaption{Top boundary $\Gamma_+$}
%     \label{top_boundary_Gamma_plus}
%     \vspace{4ex}
%   \end{minipage}%%
%   \begin{minipage}[b]{0.4\linewidth}
%     \centering
%     \includegraphics[width=.5\linewidth]{images/gamma_sf.png}
%     \subcaption{Inner boundary $\Gamma_{sf}$}
%     \label{inner_boundary_Gamma_sf}
%     \vspace{4ex}
%   \end{minipage}
%   \begin{minipage}[b]{0.4\linewidth}
%     \centering
%     \includegraphics[width=.5\linewidth]{images/gamma_out.png}
%     \subcaption{Outer boundary $\Gamma_{out}$}
%     \label{out_boundary_Gamma_out}
%     \vspace{4ex}
%   \end{minipage}%%
%   \begin{minipage}[b]{0.4\linewidth}
%     \centering
%     \includegraphics[width=.5\linewidth]{images/gamma_minus.png}
%     \subcaption{Bottom boundary $\Gamma_-$}
%     \label{bottom_boundary_Gamma_minus}
%     \vspace{4ex}
%   \end{minipage}
%     \caption{Hearth geometry: three dimensional domains, the real one (a) and its simplification (b), as well as its boundaries (c-f).}
%         \label{hearth_3d_domain}
% \end{figure}
    \caption{Hearth geometry: three dimensional domains, the real one (a) and its simplification as well as its boundaries (b).}
        \label{hearth_3d_domain}
\end{figure}
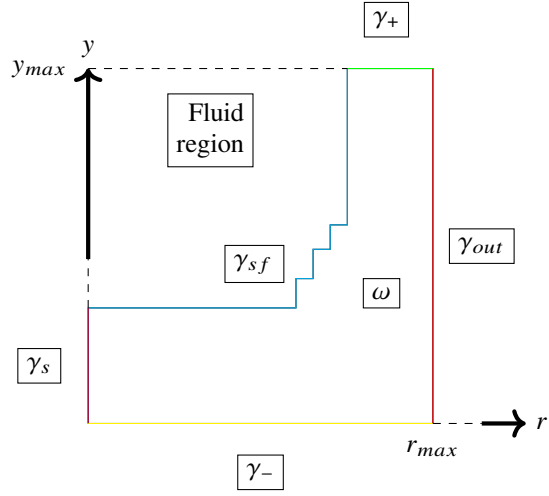
\begin{figure}[H]
\centering
    %\begin{subfigure}[t]{0.35\textwidth}
    %    \centering
    %    \includegraphics[height=2in]{images/hearth_2d_computational_domain.png}
        %\caption{Computational domain $\omega$}
    %    \caption{}
    %    \label{computational_domain_omega}
    %\end{subfigure}
    %\hspace{2cm}
    %\begin{subfigure}[t]{0.35\textwidth}
\begin{tikzpicture}[scale=0.65][
dot/.style={
  fill,
  circle,
  inner sep=2pt
  }
]
\draw (0,0) -- (7.05,0) -- (7.05,7.265) -- (5.3,7.265) -- (5.3,4.065) -- (4.95,4.065) -- (4.95,3.565) -- (4.6,3.565) -- (4.6,2.965) -- (4.25,2.965) -- (4.25,2.365) -- (0,2.365) -- (0,0);
\draw[red] (7.05,0)--(7.05,7.265);
\draw[green] (7.05,7.265)--(5.3,7.265);
\draw[cyan] (5.3,7.265) -- (5.3,4.065) -- (4.95,4.065) -- (4.95,3.565) -- (4.6,3.565) -- (4.6,2.965) -- (4.25,2.965) -- (4.25,2.365) -- (0,2.365);
\draw[purple] (0,2.365) -- (0,0);
\draw[yellow] (0,0) -- (7.05,0);
\draw[dashed] (7.05,0) -- (8.05,0);
\draw[->,ultra thick] (8.05,0)--(8.95,0) node[right]{$r$};
\draw[dashed] (0,2.365) -- (0,3.365);
\draw[->,ultra thick] (0,3.365)--(0,7.265) node[above]{$y$};
\node[draw,align=right] at (3.4,3.265) {$\gamma_{sf}$};
\node[draw,align=right] at (8.05,3.63) {$\gamma_{out}$};
\node[draw,align=right] at (-1,1.18) {$\gamma_s$};
\node[draw,align=right] at (3.525,-1) {$\gamma_-$};
\node[draw,align=right] at (6.1,8.265) {$\gamma_+$};
\node[draw,align=right] at (6,2.665) {$\omega$};
\node at (7.05,-0.5) {$r_{max}$};
\node at (-1,7.265) {$y_{max}$};
\draw[dashed] (0,7.265) -- (5.3,7.265);
\node[draw,align=right] at (2.5,6) {Fluid \\ region};
\end{tikzpicture}
%\caption{Boundaries $\partial \omega$ of computational domain $\omega$}
%\caption{}
%\label{hearth_2d_domain_boundaries}
%\end{subfigure}
\caption{\revv{Vertical section of the hearth geometry $\omega$ and its boundaries.}}
\label{hearth_2d_domain}
\end{figure}

\subsection{Thermo-mechanical model in cylindrical coordinates}\label{chapter_thermo-mechanical_model_2}

We consider the three dimensional domain $\Omega$ as in Figure \ref{hearth_3d_domain} (b), corresponding to the simplified hearth geometry. We represent the displacement vector field as $\overrightarrow{u}$ % = \left[u_r \ u_y \ u_{\theta} \right]$ 
and the temperature scalar field as $T$. Based on the simplifications listed in Sec. \ref{chapter_physicalproblem}, the energy and momentum conservation equations \cite{Bermudez,gurtin,gurtin_fried_anand_2010} can be written in a cylindrical coordinates system, $(r, y, \theta)$, with $(r,y) \in \omega$, the vertical cross section of $\Omega$ in $r-y$ plane, (see Figure \ref{hearth_2d_domain}), and $\theta \in  [0, 2 \pi)$, as,  

\begin{equation}\label{energy_equation_3d_cylindrical}
- \frac{1}{r} \frac{\partial}{\partial r} \left( r k \frac{\partial T}{\partial r} \right) - \frac{\partial}{\partial y} \left( k \frac{\partial T}{\partial y} \right) - \frac{1}{r} \frac{\partial}{\partial \theta} \left( \frac{k}{r} \frac{\partial T}{\partial \theta} \right) = Q \ , \ \text{in} \ \Omega \ ,
\end{equation}

\begin{flalign}\label{momentum_equation_3d_cylindrical}
\begin{split}
\frac{\partial \sigma_{rr}}{\partial r} + \frac{\partial \sigma_{ry}}{\partial y} + \frac{1}{r} \frac{\partial \sigma_{r\theta}}{\partial \theta} + \frac{\sigma_{rr} - \sigma_{\theta \theta}}{r} + f_{0,r} & = 0 \ , \ \text{in $\Omega$} \ ,\\
\frac{\partial \sigma_{r\theta}}{\partial r} + \frac{\partial \sigma_{\theta y}}{\partial y} + \frac{1}{r} \frac{\partial \sigma_{\theta \theta}}{\partial \theta} + 2 \frac{\sigma_{r \theta}}{r} + f_{0,\theta} & = 0\ , \ \text{in $\Omega$} \ ,\\
\frac{\partial \sigma_{ry}}{\partial r} + \frac{1}{r} \frac{\partial \sigma_{\theta y}}{\partial \theta} + \frac{\partial \sigma_{yy}}{\partial y} + \frac{\sigma_{ry}}{r} + f_{0,y} & = 0 \ , \ \text{in $\Omega$} \ .
\end{split}
\end{flalign}

Here $\overrightarrow{f}_0$ is the body force density term, $Q$ is the heat source term and $k>0$ is the thermal conductivity tensor. %In the case of isotropic material, the thermal conductivity $\bm{K}$ is expressed as,
%\begin{equation}\label{isotropic_case}
%\bm{K} = k \bm{I} \ \quad \text{with} \ k > 0 \ ,
%\end{equation}
%where $\bm{I}$ refers to the identity tensor. 
The thermo-mechanical stress tensor $\bm{\sigma}$ is related to the strain tensor $\bm{\varepsilon}$ through the Hooke's law:
\begin{equation}\label{stress_tensor}
\bm{\sigma} (\overrightarrow{u}) [T] = \lambda \Tr (\bm{\varepsilon}(\overrightarrow{u})) \bm{I} + 2 \mu \bm{\varepsilon}(\overrightarrow{u}) - (2 \mu + 3 \lambda) \alpha (T - T_0) \bm{I} \ ,
\end{equation}
where $T_0$ is the reference temperature, $\bm{\varepsilon}$ is defined as,
\begin{equation}\label{strain_tensor}
\begin{split}
\bm{\varepsilon} (\overrightarrow{u}) = \frac{1}{2} (\nabla \overrightarrow{{u}} + \nabla \overrightarrow{{u}}^T) \ ,
\end{split}
\end{equation}
$\alpha$ is the thermal expansion coefficient, and $\lambda$ and $\mu$ are the Lam\'e parameters of the material. These latter can be expressed in terms of Young modulus, $E$, and Poisson ratio, $\nu$, as:
\begin{equation}\label{lame_parameters}
\begin{split}
\mu = \frac{E}{2(1 + \nu)} \ , \ \lambda = \frac{E \nu}{(1 - 2 \nu)(1 + \nu)} \ .
\end{split}
\end{equation}

If $\bm{A}$ denotes the matrix,%fourth order tensor defined as:
\begin{equation}\label{elasticity_matrix_full}
\bm{A} = \frac{E}{(1 - 2 \nu)(1 + \nu)}
\begin{bmatrix}
1-\nu & \nu & \nu & 0 & 0 & 0 \\
\nu & 1-\nu & \nu & 0 & 0 & 0 \\
\nu & \nu & 1-\nu & 0 & 0 & 0 \\
0 & 0 & 0 & \frac{1 - 2 \nu}{2} & 0 & 0 \\
0 & 0 & 0 & 0 & \frac{1 - 2 \nu}{2} & 0 \\
0 & 0 & 0 & 0 & 0 & \frac{1 - 2 \nu}{2} \\
\end{bmatrix} \ ,
\end{equation}
the stress-strain relationship \eqref{stress_tensor} can be expressed in vector formulation as,
\begin{equation}\label{stress_cylindrical_coordinate}
\lbrace \bm{\sigma}(\overrightarrow{u})[T] \rbrace = \bm{A} \lbrace \bm{\varepsilon}(\overrightarrow{u}) \rbrace - (2 \mu + 3 \lambda) \alpha (T - T_0) \lbrace \bm{I} \rbrace \ ,
\end{equation}
where the following column vectors have been considered:
\begin{flalign}
\begin{split}
\lbrace \bm{\sigma} \rbrace & = 
\lbrace \begin{matrix}
\sigma_{rr} & \sigma_{yy} & \sigma_{\theta \theta} & \sigma_{y\theta} & \sigma_{r\theta} & \sigma_{ry} \\
\end{matrix} \rbrace^T \ , \\
\lbrace \bm{\varepsilon} \rbrace & = 
\lbrace \begin{matrix}
\varepsilon_{rr} & \varepsilon_{yy} & \varepsilon_{\theta \theta} & 2 \varepsilon_{y \theta} & 2 \varepsilon_{r \theta} & 2 \varepsilon_{ry} \\
\end{matrix} \rbrace^T \ , \\
\lbrace \bm{I} \rbrace & = 
\lbrace \begin{matrix}
1 & 1 & 1 & 0 & 0 & 0 \\
\end{matrix} \rbrace^T \ .
\end{split}
\end{flalign}

%As will be discussed in Sec. \ref{chapter_thermo-mechanical_model_3}, the simplifications assumed allow us to consider that $\Omega$ has axial symmetry. Then, the equations \eqref{energy_equation_3d_cylindrical}-\eqref{momentum_equation_3d_cylindrical} will be solved on its vertical section, $\omega$, shown in Figure \ref{hearth_2d_domain}.

We introduce the notations for the boundaries of the domain $\Omega$ (Figure \ref{hearth_3d_domain}), and its vertical cross section in $r-y$ plane, $\omega$ (Figure \ref{hearth_2d_domain}):

\begin{flalign}\label{eq:definitions}
\begin{split}
\Gamma_{out} = \partial \Omega \cap (r \equiv r_{max}) \ &= \gamma_{out} \times [0,2 \pi) \ ,\\
\Gamma_{+} = \partial \Omega \cap (y \equiv y_{max}) \ &= \gamma_{+} \times [0,2 \pi) \ ,\\
\Gamma_{-} = \partial \Omega \cap (y \equiv 0) \ &= \gamma_{-} \times [0,2 \pi) \ ,\\
\Gamma_{sf} = \partial \Omega \backslash (\Gamma_{out} \cup \Gamma_{+} \cup \Gamma_{-}) \ &= \gamma_{sf} \times [0,2 \pi) \ , \\
\gamma_s &= \partial \omega \cap (r \equiv 0) \ ,
\end{split}
\end{flalign}
where $r_{max} \in \mathbb{R}^+$ and $y_{max} \in \mathbb{R}^+$.

Moreover, on the boundary, we define normal force $\sigma_n$ and tangential force $\overrightarrow{\sigma_t}$ as follows:
\begin{equation}\label{normal_shear_stress}
\begin{split}
\sigma_n = (\bm{\sigma} \overrightarrow{n})\cdot \overrightarrow{n} \ , \ \overrightarrow{\sigma_t} = \bm{\sigma} \overrightarrow{n} - \sigma_n \overrightarrow{n} \ ,
\end{split}
\end{equation}
where $\overrightarrow{n}$ is the unit normal vector that is directed outwards from $\Omega$.

\begin{itemize}

\item On the upper boundary, $\Gamma_{+}$, the applied force, $\overrightarrow{g}_+$, and the density of heat flux, $q_+$, are known. Therefore, the following boundary conditions are considered:

\begin{equation}\label{boundary_upper}
\begin{split}
%(- \bm{K} \nabla T) \cdot \overrightarrow{n} = q_+ \ , \\
(- k \nabla T) \cdot \overrightarrow{n} = q_+ \ , \
\bm{\sigma} \overrightarrow{n} = \overrightarrow{g}_+ \ .
\end{split}
\end{equation}

\item On the bottom boundary, $\Gamma_{-}$, the convection heat transfer occurs with the heat exchanger %cooling medium %\michele{Physically, who is the "heat exchanger"?} \nirav{Removed "heat exchanger" with "cooling medium"}
at temperature $T_-$ and heat transfer coefficient $h_{c,-}$. The normal displacement is null, and we denote with $\overrightarrow{g}_-$ the shear force. %referred to the friction at the bottom,
%denoted with $\overrightarrow{g}_-$, is assumed to have zero normal component. %i.e. $\overrightarrow{g}_- \cdot \overrightarrow{n}=0$..
Therefore, it is verified:

\begin{equation}\label{boundary_bottom}
\begin{split}
%(- \bm{K} \nabla T) \cdot \overrightarrow{n} = h_{c,-}(T - T_-) \ , \\
(- k \nabla T) \cdot \overrightarrow{n} = h_{c,-}(T - T_-) \ , \
\overrightarrow{u} \cdot \overrightarrow{n} = 0 \ , \ \overrightarrow{\sigma_t} = \overrightarrow{g}_- \ .
\end{split}
\end{equation}
%where the $\overrightarrow{g}_-$ is assumed to have zero normal component, i.e. $\overrightarrow{g}_- \cdot \overrightarrow{n}=0$. %The shear force $\overrightarrow{g}_-$ refers to friction at the bottom surface. % case of frictionless surface, the shear force $\overrightarrow{g}_-$ becomes zero.

\item On the inner boundary, $\Gamma_{sf}$, convection heat transfer with the fluid phase occurs and hydrostatic pressure  ${g}_{sf}$ is acting. So, the following boundary conditions are considered:

\begin{equation}\label{boundary_inner}
\begin{split}
%(- \bm{K} \nabla T) \cdot \overrightarrow{n} = h_{c,f}(T - T_f) \ , \\
(- k \nabla T) \cdot \overrightarrow{n} = h_{c,f}(T - T_f) \ , \
\bm{\sigma} \overrightarrow{n} = \overrightarrow{g}_{sf} \ ,
\end{split}
\end{equation}
where $T_f$ is the fluid temperature, assumed to be known and constant at the steady state, and $h_{c,f}$ is the convective heat transfer coefficient. In addition, $\overrightarrow{g}_{sf} = -g_{sf} \overrightarrow{n}$.%, and $\overrightarrow{g}_{sf} = -p_h \overrightarrow{n}$ is the hydrostatic pressure.

\item On the outer boundary, $\Gamma_{out}$, a convective heat flux and known applied force $\overrightarrow{g}_{out}$ are assumed:

\begin{equation}\label{boundary_outer}
\begin{split}
%(-\bm{K} \nabla T) \cdot \overrightarrow{n} = h_{c,out} (T - T_{out}) \ , \\
(-k \nabla T) \cdot \overrightarrow{n} = h_{c,out} (T - T_{out}) \ , \
\bm{\sigma} \overrightarrow{n} = \overrightarrow{g}_{out} \ ,
\end{split}
\end{equation}
$h_{c,out}$ being the convective heat transfer coefficient and $T_{out}$ the ambient temperature.
\end{itemize}

\subsection{Axisymmetric thermo-mechanical model}\label{chapter_thermo-mechanical_model_3}

%In the context of blast furnace application,
%\michele{(here, can we add a reference to justify the following assumptions?)} \nirav{No I do not have any reference as this is an assertion and not a proposal/assumption.} \michele{Maybe it's better to talk about this.}
In the context of blast furnace application, the body force density term $\overrightarrow{f_0}$, %can be expressed as,
%\begin{gather}
%\overrightarrow{f_0} = \overrightarrow{f_0} (r,y) = f_{0,r} (r,y) \overrightarrow{e_r} + f_{0,y} (r,y) \overrightarrow{e_y} \ .
%\end{gather}
%and it depends only on $(r,y)$ coordinates.
as well as surface forces, $\overrightarrow{g}_+$, $\overrightarrow{g}_-$, $\overrightarrow{g}_{sf}$, $\overrightarrow{g}_{out}$, have zero component in $\overrightarrow{e_{\theta}}$ direction and they do not depend on $\theta$. Besides, the heat source term, $Q$, the heat flux density, $q_+$, the heat transfer coefficients, $h_{c,-} , h_{c,f}, h_{c,out}$, and temperatures $T_-, T_f, T_{out}$ are assumed to be only dependent on $(r,y)$ coordinates. Therefore, a symmetry hypothesis is applicable that leads to significant computational savings. %the three dimensional model \eqref{energy_equation_3d_cylindrical}-\eqref{momentum_equation_3d_cylindrical_bc} does not depend on $\theta$, and hence a symmetry hypothesis is applicable.
%The axisymmetry leads to significant computational savings as the three dimensional model \eqref{energy_equation_3d_cylindrical}-\eqref{momentum_equation_3d_cylindrical_bc} defined in $\Omega$ (Figure \ref{hearth_3d_domain})
%is replaced by the corresponding two dimensional model defined in its vertical section, $\omega$ ( Figure \ref{hearth_2d_domain}). %The unit outer normal vector to the boundary of this section will now be represented as $\overrightarrow{n} = n_r \overrightarrow{e_r} + n_y \overrightarrow{e_y}$. In the axisymmetric system, we represent the displacement $\overrightarrow{{u}}$ and temperature $T$, both independent of $\theta$, as
%\begin{equation}
%\begin{split}
%\overrightarrow{u} = u_r(r,y) \overrightarrow{e_r} + u_y (r,y) \overrightarrow{e_y} \ , \\
%T = T(r,y) \ .
%\end{split}
%\end{equation}

The associated axisymmetric model is reduced to consider conservation equations \eqref{energy_equation_3d_cylindrical}, \eqref{momentum_equation_3d_cylindrical} defined in $\omega$ and the boundary conditions \eqref{boundary_upper} - \eqref{boundary_outer}, where the $\Gamma$ boundaries are replaced by $\gamma$ such that $\Gamma = (\gamma \backslash \gamma_s) \times [0,2\pi)$ (see Figures \ref{hearth_3d_domain} and \ref{hearth_2d_domain}, and definitions \eqref{eq:definitions}). Moreover, the usual symmetry conditions on $\gamma_s$ are added:
\begin{equation}
\begin{split}
(-k \nabla T) \cdot \overrightarrow{n} = 0 \ , \
\overrightarrow{u} \cdot \overrightarrow{n} = 0 \ , \ \overrightarrow{\sigma_t}  = \overrightarrow{0} \ .
\end{split}
\end{equation}

Therefore, the axisymmetric thermo-mechanical model can be summarized as:

\begin{itemize}

\item \textbf{Thermal model}:
\begin{equation}\label{summarised_strong_energy_equations}
- \frac{1}{r} \frac{\partial}{\partial r} \left( r k \frac{\partial T}{\partial r} \right) - \frac{\partial}{\partial y} \left( k \frac{\partial T}{\partial y} \right) = Q \ , \ \text{in} \ \omega \ ,
\end{equation}

with boundary conditions:
\begin{flalign}\label{summarised_strong_energy_equations_boundary}
\begin{split}
\text{on $\gamma_+$ } & : -k \frac{\partial T}{\partial y} = q_+ \ , \\
\text{on $\gamma_-$ } & : k \frac{\partial T}{\partial y} = h_{c,-} (T-T_-) \ , \\
\text{on $\gamma_{sf}$ } & : -k \frac{\partial T}{\partial r} n_r - k \frac{\partial T}{\partial y} n_y  = h_{c,f} (T - T_f) \ , \\
\text{on $\gamma_{out}$ } & : -k \frac{\partial T}{\partial r} = h_{c,out} (T - T_{out}) \ , \\
\text{on $\gamma_s$ } & : \frac{\partial T}{\partial r} = 0 \ .
\end{split}
\end{flalign}

\item \textbf{Mechanical model}:

\begin{flalign}\label{summarised_strong_momentum_equations}
\begin{split}
\frac{\partial \sigma_{rr}}{\partial r} + \frac{\partial \sigma_{ry}}{\partial y}  + \frac{\sigma_{rr} - \sigma_{\theta \theta}}{r} + f_{0,r} & = 0 \ , \ \text{in} \ \omega \ , \\
\frac{\partial \sigma_{ry}}{\partial r} + \frac{\partial \sigma_{yy}}{\partial y} + \frac{\sigma_{ry}}{r} + f_{0,y} & = 0 \ , \ \text{in} \ \omega \ ,
\end{split}
\end{flalign}

%\michele{Check the following formulas with respect to eqs. (1.8)-(1.11) please. I seem that in the last ones you already have introduced cylindric system and axisymmeetry hipothesis that is of course not coherent with the context.}
In vector notation, axisymmetric stress-strain relationship can be expressed as,
\begin{flalign}
\begin{split}\label{stress_strain_relation_axisymmetric}
\lbrace \bm{\sigma}(\overrightarrow{u})[T] \rbrace & = \bm{A} \lbrace \bm{\varepsilon}(\overrightarrow{u}) \rbrace - (2 \mu + 3 \lambda) \alpha (T - T_0) \lbrace \bm{I} \rbrace \ , \\
\bm{A} & = \frac{E}{(1 - 2 \nu)(1 + \nu)}
\begin{bmatrix}
1-\nu & \nu & \nu & 0 \\
\nu & 1-\nu & \nu & 0 \\
\nu & \nu & 1-\nu & 0 \\
0 & 0 & 0 & \frac{1 - 2 \nu}{2} \\
\end{bmatrix} \ , \\
\lbrace \bm{\sigma} \rbrace & =
\lbrace \begin{matrix}
\sigma_{rr} & \sigma_{yy} & \sigma_{\theta \theta} & \sigma_{ry} \\
\end{matrix}\rbrace^T \ , \\
\lbrace \bm{\varepsilon} \rbrace & =
\lbrace \begin{matrix}
\varepsilon_{rr} & \varepsilon_{yy} & \varepsilon_{\theta \theta} & 2 \varepsilon_{ry} \\
\end{matrix}\rbrace^T \ , \\
\lbrace \bm{I} \rbrace & =
\lbrace \begin{matrix}
1 & 1 & 1 & 0 \\
\end{matrix}\rbrace^T \ ,
\end{split}
\end{flalign}

with the boundary conditions :

\begin{flalign}\label{summarised_strong_momentum_equations_boundary}
\begin{split}
\textrm{on} \ \gamma_+ \ & : \sigma_{ry} = g_{+,r} \ , \ \sigma_{yy} = g_{+,y} \ , \\
\textrm{on} \ \gamma_- \ & : \ u_y = 0, \ \sigma_{ry} = -g_{-,r} \ , \\
\textrm{on} \ \gamma_{sf} \ & : \ \sigma_{rr}n_r + \sigma_{ry}n_y = g_{sf,r} \ , \ \\ & \ \ \ \sigma_{ry}n_r + \sigma_{yy}n_y = g_{sf,y} \ , \\
\textrm{on} \ \gamma_{out} \ & : \ \sigma_{rr} = g_{out,r} \ , \ \sigma_{ry} = g_{out,y} \ , \\
\textrm{on} \ \gamma_s \ & : \ u_r = 0 \ , \ \sigma_{ry} = 0 \ .
\end{split}
\end{flalign}

\end{itemize}

%\nirav{It shoud be noted that $\varepsilon_{r \theta} = \varepsilon_{ry} = 0$ due to axisymmetry hypothesis.}

\section[Weak formulation]{Weak formulation of axisymmetric thermo-mechanical model}\label{weak_form_chapter}
In this section, we derive the weak formulation related to the axisymmetric thermo-mechanical model \eqref{summarised_strong_energy_equations}-\eqref{summarised_strong_momentum_equations_boundary}. % introduced in the section \ref{chapter_thermo-mechanical_model_3}.
First, in Sec. \ref{sec:sobolev}, we introduce the relevant function spaces for temperature and displacement fields as well as for the model data, including boundary conditions, physical properties, and source terms, such that the problem is well defined. Next, the weak formulations for the thermal and mechanical problems are reported in Sec. \ref{ref:weakT} and \ref{ref:weakU}, respectively. %we investigate the well-posedness of the problem, i.e. the existence and uniqueness of the weak solution.

\subsection{Functional spaces}\label{sec:sobolev}
%\michele{The definition of $L^\infty$ misses.}\nirav{added}

%\michele{I moved here the definition of function spaces for the mechanical problem for sake of order and coherence.} \nirav{Noted}

%\michele{Nirav, i would insert relevant reference(s) for function spaces introduced.} \nirav{Added}

For data of both thermal and mechanical problem, we introduce the weighted Sobolev spaces, $L^2_r(\omega)$, with norm $||\cdot||_{L^2_r(\omega)}$ as,
\begin{equation}\label{norms_thermal_model_1}
\begin{split}
L^2_r(\omega) = \bigg\lbrace f: \omega \mapsto \mathbb{R} \ , \ \int_{\omega} f^2 r dr dy < \infty \bigg\rbrace \ , \\
||f||_{L^2_r(\omega)}^2 = \int_{\omega} f^2 r dr dy \ .
\end{split}
\end{equation}
Analogously, given $\gamma$ a subset of $\partial \omega$, the boundary of $\omega$,  
\begin{equation}\label{norms_thermal_model_2}
L^2_r(\gamma) = \bigg\lbrace g:\gamma \mapsto \mathbb{R}, \int_{\gamma} g^2 r d\gamma < \infty \bigg\rbrace \ .
\end{equation}

Let $L^{\infty}(\omega)$ be the space
\begin{equation}\label{norms_thermal_model_infty}
\begin{split}
L^{\infty}(\omega) = \lbrace f : \omega \mapsto \mathbb{R} \ , \sup_{\omega}|f| \leq C \ , \ C \geq 0 \ \rbrace \ , \\
||f||_{L^{\infty}(\omega)} = \sup_{\omega}|f| \ .
\end{split}
\end{equation}
Analogously, $L^{\infty}(\gamma)$ is defined.

For the temperature, we introduce the weighted Sobolev space, $H^1_r(\omega)$, with norm $||\cdot||_{H^1_r(\omega)}$ as,
\begin{equation}\label{norms_thermal_model_2}
\begin{split}
H^1_r(\omega) = \bigg\lbrace \psi : \omega \mapsto \mathbb{R} \ , \ \int_{\omega} \left(\psi^2 + \left( \frac{\partial \psi}{\partial r} \right)^2 +  \left( \frac{\partial \psi}{\partial y} \right)^2 \right) r dr dy < \infty \bigg\rbrace \ , \\
||\psi||_{H^1_r(\omega)}^2 = \int_{\omega} \left( \psi^2 + \left( \frac{\partial \psi}{\partial r} \right)^2 +  \left( \frac{\partial \psi}{\partial y} \right)^2 \right) r dr dy \ .
\end{split}
\end{equation}

%Notice that $H^1_r(\omega)$ is closely related to the usual Sobolev space $H^1(\Omega)$. Indeed, if $\psi$ belongs to $L^2_r(\omega)$, and we extend it to $\Omega = \omega \times [0,2\pi)$ as an axisymmetric function $\bar{\psi}$ defined as: $\bar{\psi}(r,y,\theta)=\psi(r,y)$, then:
%\begin{equation}\label{isometry}
%\begin{split}
%||\bar{\psi}||_{L^2(\Omega)} = \sqrt{2 \pi} ||\psi||_{L^2_r(\omega)} \ , \\
%||\bar{\psi}||_{H^1(\Omega)} = \sqrt{2 \pi} ||\psi||_{H^1_r(\omega)} \ .
%\end{split}
%\end{equation}
%The reciprocal is also true for all axisymmetric function $\bar{\psi}$ with respect to the cylindrical coordinates $(r,y,\theta)$. Let $\bar{H}^1(\Omega) \subset H^1(\Omega)$ be the subspace of all axisymmetric functions in $H^1(\Omega)$ with respect to $y-$axis. Then the following properties are verified (see \cite{axisymmetric_polygonal}):
%\begin{theorem}
%\begin{itemize}
%\item $\bar{H}^1(\Omega)$ is isometric to $H^1_r(\omega)$.
%\item The space of smooth functions $C^{\infty}(\omega)$ is dense in $L^2_r(\omega)$ and in $H^1_r(\omega)$.
%\item By isometry, it can be concluded that the embedding of $H^1_r(\omega)$ in $L^2_r(\omega)$ is compact.
%\end{itemize}
%\end{theorem}

On the other hand, the following space $\mathbb{V}$ for the displacement is considered:
\begin{equation}
\mathbb{V} = (H^1_r(\omega) \cap L^2_{1/r} (\omega)) \times H^1_r(\omega) \ .
\end{equation}
It will be equipped with the inner product, %$<\overrightarrow{u},\overrightarrow{\phi}>_{\mathbb{V}} = $
\begin{align}
<\overrightarrow{u},\overrightarrow{\phi}>_{\mathbb{V}} = \int_{\omega} \left( \phi_r u_r + \phi_y u_y + \frac{\partial u_r}{\partial r} \frac{\partial \phi_r}{\partial r} + \frac{\partial u_r}{\partial y} \frac{\partial \phi_r}{\partial y} + \frac{u_r}{r} \frac{\phi_r}{r} + \frac{\partial u_y}{\partial r} \frac{\partial \phi_y}{\partial r} + \frac{\partial u_y}{\partial y} \frac{\partial \phi_y}{\partial y} + \right.
\left. \frac{\partial u_r}{\partial y} \frac{\partial \phi_y}{\partial r} + \frac{\partial u_y}{\partial r} \frac{\partial \phi_r}{\partial y} \right) r dr dy \ ,
\end{align}
and with the norm,
\begin{gather}\label{norms_mechanical_model}
||\overrightarrow{\phi}||_{\mathbb{V}}^2 = <\overrightarrow{\phi},\overrightarrow{\phi}>_{\mathbb{V}} \ .
\end{gather}
Its closed and convex subspace $\mathbb{U}$,
\begin{equation}
\begin{split}
\mathbb{U} = \lbrace \overrightarrow{\phi} = \left( \phi_r \ \phi_y \right) \in  \mathbb{V} \ , \ \phi_y = 0 \ \text{on} \ \gamma_- \ , \ \phi_r = 0 \ \text{on} \ \gamma_s \ \rbrace \ ,
\end{split}
\end{equation}
will be the set of admissible displacements. The subspace $\mathbb{U}$ is equipped with the same norm as spave $\mathbb{V}$ i.e. $||\overrightarrow{\phi}||_{\mathbb{U}} = ||\overrightarrow{\phi}||_{\mathbb{V}}, \forall \overrightarrow{\phi} \in \mathbb{U}$.

Finally, the function space for stress tensor is defined as, %\michele{$\sigma$\_{$\alpha$ 3}...who is $\alpha$?!} \nirav{Corrected},
\begin{equation}
\mathbb{S} = \lbrace \bm{\sigma} = [\sigma_{ij}] \in [L_r^2(\omega)]^{3 \times 3}, \sigma_{ij} = \sigma_{ji}, \sigma_{\alpha 3} = 0, \alpha = 1,2  \rbrace \ .
\end{equation}

For more details the reader is referred, e.g., to \cite{axisymmetric_ivan_korn,axisymmetric_ivan,axisymmetric_polygonal}.

\subsection{Weak formulation for thermal model}\label{ref:weakT}

We assume the following hypotheses on the data:

%\begin{itemize}
\begin{description}

\item[(TH1)] The heat source term, $Q$, is such that
\begin{equation*}
Q \in L^2_r(\omega) \ .
\end{equation*}

\item[(TH2)] The convection temperatures, $T_{sf}$, $T_-$ and $T_{out}$, as well as the boundary heat flux $q_{+}$ are such that
\begin{equation*}
T_{sf} \in L^2_r(\gamma_{sf}) \ , \ T_- \in L^2_r(\gamma_-) \ , \ T_{out} \in L^2_r(\gamma_{out}) \ , \ q_+ \in L_r^2(\gamma_+) \ .
\end{equation*}

%\item[(TH3)] The boundary heat flux $q_{+}$ is such that
%\begin{equation*}
%q_+ \in L_r^2(\gamma_+) \ .
%\end{equation*}

%There exists constants $k_0 > 0$, $h_{c,f,0} > 0$, $h_{c,out,0} > 0$, $h_{c,-,0} > 0$ such that %There exists a constant $k_0 > 0$, such that the thermal conductivity, $k(r,y)$ satisfies,
%\begin{equation*}
%k(r,y) \in L^\infty (\omega) \ , \ k(r,y) > k_0 > 0 \ , \ \text{in} \ \omega \ .
%\end{equation*}
%Also, there exist constants $h_{c,f,0} > 0, h_{c,out,0} > 0, h_{c,-,0} > 0$ such that,
\item[(TH3)] The thermal conductivity $k(r,y)$ and the convective heat transfer coefficients, $h_{c,f}(r,y)$, $h_{c,out}(r,y)$ and $h_{c,-}(r,y)$ are such that %  There exists constants $k_0 > 0$, $h_{c,f,0} > 0$, $h_{c,out,0} > 0$, $h_{c,-,0} > 0$ such that %There exists a constant $k_0 > 0$, such that the thermal conductivity, $k(r,y)$ satisfies,
\begin{gather*}
k(r,y) \in L^\infty (\omega) \ , \ k(r,y) > k_0 > 0 \ , \\ %\ \text{in} \ \omega \ , \\
h_{c,f}(r,y) \in L^\infty(\gamma_{sf}) \ , \ h_{c,f}(r,y) > h_{c,f,0} > 0 \ , \\% \ \text{on} \ \gamma_{sf} \ , \\
h_{c,out}(r,y) \in L^\infty(\gamma_{out}) \ , \ h_{c,out}(r,y) > h_{c,out,0} > 0 \ , \\ %\text{on} \ \gamma_{out} \ , \\
h_{c,-}(r,y) \in L^\infty(\gamma_-) \ , \ h_{c,-}(r,y) > h_{c,-,0} > 0 \ , %\ \text{on} \ ,% \gamma_- \ ,
\end{gather*}
where $k_0$, $h_{c,f,0}$, $h_{c,out,0}$ and $h_{c,-,0}$ are suitable constants.
\end{description}
%\end{itemize}

In order to propose a weak formulation for the thermal model \eqref{summarised_strong_energy_equations} - \eqref{summarised_strong_energy_equations_boundary}, %\eqref{summarised_strong_energy_equations} and \eqref{summarised_strong_energy_equations_boundary},
in the following we assume sufficient regularity to perform the calculations. We multiply the energy equation \eqref{summarised_strong_energy_equations} by $r\psi(r,y)$, integrate over the domain $\omega$ with respect to $(r,y)$ variables, apply Gauss divergence theorem and use boundary conditions \eqref{summarised_strong_energy_equations_boundary} to obtain:
\begin{equation}\label{axisymteric_energy_equation_weak_3}
\begin{split}
\int_{\omega} rk \left( \frac{\partial T}{\partial y} \frac{\partial \psi}{\partial y} + \frac{\partial T}{\partial r} \frac{\partial \psi}{\partial r} \right) dr dy + \int_{\gamma_{sf}} \psi h_{c,f} T r d \gamma + \int_{\gamma_{out}} \psi h_{c,out} T r d \gamma + \\ \int_{\gamma_-} \psi h_{c,-} T r d \gamma  = \int_{\omega} \psi Q r dr dy + \int_{\gamma_{sf}} \psi h_{c,f} T_f r d \gamma + \\ \int_{\gamma_{out}} \psi h_{c,out} T_{out} r d \gamma + \int_{\gamma_{-}} \psi h_{c,-} T_- \ r d \gamma - \int_{\gamma_+} \psi q^+ r d \gamma \ .
\end{split}
\end{equation}

%Therefore, we propose the following weak formulation for thermal problem $(T1)$: \begin{problem} Weak thermal formulation $(WT1)$ :
%Under the assumptions (TH1)-(TH4), find $T \in H^1_r(\omega)$, such that equality \eqref{axisymteric_energy_equation_weak_3} is verified for all $\psi \in  H^1_r(\omega)$. \end{problem}

It is to be noted that under assumptions (TH1)-(TH3) all integrals of the proposed weak formulation are well defined  for all $T,\psi \in H^1_r(\omega)$. The left hand side of equation \eqref{axisymteric_energy_equation_weak_3} is bilinear and symmetric. So, we define in $H^1_r(\omega) \times H^1_r(\omega)$ the operator:
\begin{equation}\label{bilinear_form_thermal}
\begin{split}
a_T(T,\psi) = \int_{\omega} rk \left( \frac{\partial T}{\partial y} \frac{\partial \psi}{\partial y} + \frac{\partial T}{\partial r} \frac{\partial \psi}{\partial r} \right) dr dy + \int_{\gamma_{sf}} \psi h_{c,f} T r d \gamma + \int_{\gamma_{out}} \psi h_{c,out} T r d\gamma + \\ \int_{\gamma_-} \psi h_{c,-} T r d\gamma \ .
\end{split}
\end{equation}
The right hand side of equation \eqref{axisymteric_energy_equation_weak_3} is linear and the following operator defined on $H^1_r(\omega)$ is introduced:
\begin{equation}\label{linear_form_thermal}
\begin{split}
l_T(\psi) = \int_{\omega} \psi Q r dr dy + \int_{\gamma_{sf}} \psi h_{c,f} T_f r d \gamma + \int_{\gamma_{out}} \psi h_{c,out} T_{out} r d\gamma + \\ \int_{\gamma_{-}} \psi h_{c,-} T_- \ r d\gamma -  \int_{\gamma_+} \psi q^+ r d\gamma \ .
\end{split}
\end{equation}

Then, we can define the following problem: %\begin{problem}
\begin{itemize}
\item \textbf{Weak thermal model $\textbf{(WT)}$} : Under the assumptions (TH1)-(TH3), find $T \in H^1_r(\omega)$ such that,
\begin{equation}\label{bilinear_linear_form_thermal}
a_T(T,\psi) = l_T(\psi) \ , \ \forall \ \psi \in H^1_r(\omega) \ .
\end{equation} \end{itemize}

By using Cauchy-Schwarz inequality, the trace operator properties, and Friedrich's inequality  \cite{brezis, necas}, it can be shown that, under the assumptions (TH1)-(TH3), $a_T(T,\psi)$ and $l_T(\psi)$ are continuous on $H^1_r(\omega) \times H^1_r(\omega)$ and $H^1_r(\omega)$, respectively, and $a_T(\psi,\psi)$ is coercive on  $H^1_r(\omega) \times H^1_r(\omega)$. %In other words, there exists a constant, $C_T > 0$, referred as continuity constant of $a_T(\cdot,\cdot)$, such that,
%\begin{equation}
%a_T(T,\psi) \leq C_T ||T||_{H^1_r(\omega)} ||\psi||_{H^1_r(\omega)} \ , \ \forall T,\psi \in H^1_r(\omega) \ .
%\end{equation}
%Besides, since \michele{(What means "mes"?!)} \michele{"mes" is tipically used?} \nirav{the length of $\gamma_- \cup \gamma_{out} \cup \gamma_{sf}$ is non-zero} i.e. $mes(\gamma_- \cup \gamma_{out} \cup \gamma_{sf}) > 0$, by Friedrich's inequality (see theorem 1.9 of \cite{necas}) as applied to $H^1(\Omega)$, using equation \eqref{isometry}, and hypotheses $(TH4)$ on the data, one can prove that the bilinear form $a_T(\psi,\psi)$ is coercive on $H^1_r(\omega) \times H^1_r(\omega)$. %there exists a constant, $c_T > 0$, referred as coercivity constant of $a_T(\cdot,\cdot)$, such that,
%\begin{equation}
%\begin{split}
%c_T ||\psi||^2_{H^1_r(\omega)} \leq a_T(\psi,\psi) \ , \ \forall \psi \in H^1_r(\omega) \ .
%\end{split}
%\end{equation}
%In other words, the bilinear form $a_T(\psi,\psi)$ is coercive on $H^1_r(\omega) \times H^1_r(\omega)$.
Hence the conditions of the Lax-Milgram theorem \cite{brezis} are satisfied and accordingly the weak thermal model $(WT)$ has a unique solution. %$T$. %Given the equivalence between $(WT1)$ and $(WT2)$ formulations, $T$ is the unique solution of weak thermal formulation $(WT1)$.

\subsection{Weak formulation of the mechanical model}\label{ref:weakU}

We assume the following hypotheses on the data:
\begin{description}
%\begin{itemize}

\item[(MH1)] The body force density, $\overrightarrow{f_0}$, is such that
\begin{equation*}
\overrightarrow{f_0} \in [L^2_r(\omega)]^2 \ .
\end{equation*}

\item[(MH2)] The boundary forces $\overrightarrow{g}_+$, $\overrightarrow{g}_{sf}$, $\overrightarrow{g}_{out}$ and $\overrightarrow{g}_-$ are such that
\begin{gather*}
\overrightarrow{g}_+ \in [L^2_r(\gamma_+)]^2 , \ \overrightarrow{g}_{sf} \in [L^2_r(\gamma_{sf})]^2 , \ \overrightarrow{g}_{out} \in [L^2_r(\gamma_{out})]^2 ,\ \overrightarrow{g}_- \in [L^2_r(\gamma_-)]^2. %\text{and} \ \overrightarrow{n}\cdot\overrightarrow{g}_- = 0 \  \text{on $\gamma_-$} .
\end{gather*}

\item[(MH3)] The Young's modulus $E(r,y)$, the coefficient of thermal expansion $\alpha(r,y)$, and the Poisson's ratio $\nu(r,y)$ are such that
%There exist constants $E_0 > 0$ and $\alpha_0 > 0$ such that, the Young's modulus, $E(r,y)$, and the coefficient of thermal expansion, $\alpha(r,y)$ satisfy:
\begin{gather*}
E(r,y) \in L^\infty (\omega) \ , \ E > E_0 > 0 \ , \\
\alpha(r,y) \in L^\infty (\omega) \ , \ \alpha > \alpha_0 > 0 \ , \\ %\ \text{in} \ \omega \ , \\
\nu(r,y) \in L^\infty(\omega) \ , \ \nu_0 < \nu < \nu_1 \ , \ \nu_0 > 0 \ , \ %\ \nu_1 < 0.5 \ , \ \text{in} \ \omega \ ,
\end{gather*}
where $E_0$, $\alpha_0$, $\nu_0$ and $\nu_1$ are suitable constants.

%\item[(MH4)] %There exist constants $\nu_0 > 0, \nu_1 < 0.5$ such that, the Poisson's ratio, $\nu(r,y)$, satisfies:
%\begin{equation*}
%\nu(r,y) \in L^\infty(\omega) \ , \ \nu_0 < \nu < \nu_1 \ , \ \nu_0 > 0 \ , \ \nu_1 < 0.5 \ , \ \text{in} \ \omega \ ,
%\end{equation*}
%where $\nu_0$ and $\nu_1$ are suitable constants.

%\end{itemize}
\end{description}

%It can be seen that in case $T=T_0$, we have $\bm{\sigma} \in \mathbb{S}$. %in equation \eqref{stress_tensor}, %\eqref{stress_strain_relation_axisymmetric},
%the stress tensor $\bm{\sigma}(\overrightarrow{\phi})$ belongs to the space $\mathbb{S}$.
%Besides, since $T$ and $T_0$ $\in H^1_r(\omega)$, the stresses generated belong to $\mathbb{S}$ \michele{$\mathbb{S}$?} \nirav{Done} too. In other words, the stresses generated due to mechanical effects lies also in the same space as that of stresses generated due to thermo-mechanical effects.

Analogously to what has been done for the thermal model, to propose a weak formulation of the mechanical model \eqref{summarised_strong_momentum_equations} - \eqref{summarised_strong_momentum_equations_boundary}, in the following we assume sufficient regularity to
perform the calculations. %we assume that all functions have sufficiently regularity as necessary for the following calculations.
Given a function $\overrightarrow{\phi} =(\phi_r, \phi_y)$, we multiply the first equation of \eqref{summarised_strong_momentum_equations} by $r\phi_r(r,y)$, the second one by $r\phi_y(r,y)$, we sum both, integrate over $\omega$, apply Green formula and use equation \eqref{stress_tensor} and \eqref{summarised_strong_momentum_equations_boundary} %and \eqref{summarised_strong_momentum_equations_boundary}
to obtain %and consider $\overrightarrow{\phi} \cdot \overrightarrow{n} = 0$ on $\gamma_- \cup \gamma_s$, to obtain
\begin{equation}\label{axisymteric_momentum_equation_weak_3}
\begin{split}
\int_{\omega} \bm{A} \lbrace \bm{\varepsilon}(\overrightarrow{u}) \rbrace \cdot  \lbrace \bm{\varepsilon}(\overrightarrow{\phi})\rbrace r dr dy = \int_{\omega} (2 \mu + 3 \lambda) \alpha (T - T_0) \lbrace \bm{I} \rbrace \cdot \lbrace \bm{\varepsilon}(\overrightarrow{\phi}) \rbrace r dr dy + \\ \int_{\omega} \left( \phi_r f_{0,r} + \phi_y f_{0,y}\right) r dr dy + \int_{\gamma_{sf}} \overrightarrow{\phi} \cdot \overrightarrow{g}_{sf} r d\gamma + \int_{\gamma_{out}} \overrightarrow{\phi} \cdot \overrightarrow{g}_{out} r d\gamma + \\ \int_{\gamma_-} \overrightarrow{\phi} \cdot \overrightarrow{g}_- r d\gamma + \int_{\gamma_+} \overrightarrow{\phi} \cdot \overrightarrow{g}_+ r d\gamma \ , \ \forall \overrightarrow{u}, \overrightarrow{\phi} \in \mathbb{U} \ ,
\end{split}
\end{equation}
where $T$ is assumed to be the solution of the weak thermal model $(WT)$.

Notice that under assumptions (MH1)-(MH3), and since $T \in H^1_r(\omega)$, all integrals in \eqref{axisymteric_momentum_equation_weak_3} are well defined for all $\overrightarrow{u}, \overrightarrow{\phi} \in \mathbb{U}$. %Therefore, we propose the following weak formulation for the mechanical model $(M1)$:
%\begin{problem}
%Weak mechanical formulation $(WM1)$ : Let $T \in H^1_r(\omega)$ be the solution of the weak thermal model $(WT2)$. Under assumptions (MH1)-(MH4), find $\overrightarrow{u} \in \mathbb{U}$, such that equality \eqref{axisymteric_momentum_equation_weak_3} is verified for all $\overrightarrow{\phi} \in \mathbb{U}$.
%\end{problem}

The left hand side of equation \eqref{axisymteric_momentum_equation_weak_3},
\begin{equation}
a_M(\overrightarrow{u},\overrightarrow{\phi}) = \int_{\omega} \bm{A} \lbrace \bm{\varepsilon}(\overrightarrow{u}) \rbrace \cdot  \lbrace \bm{\varepsilon}(\overrightarrow{\phi})\rbrace r dr dy \ ,
\end{equation}
is bilinear in $\mathbb{V} \times \mathbb{V}$, while the right hand side,
\begin{equation}
\begin{split}
l_M[T](\overrightarrow{\phi}) = \int_{\omega} (2 \mu + 3 \lambda) \alpha (T - T_0) \lbrace \bm{I} \rbrace \cdot \lbrace \bm{\varepsilon}(\overrightarrow{\phi}) \rbrace r dr dy + \int_{\omega} \left( \phi_r f_{0,r} + \phi_y f_{0,y}\right) r dr dy + \\ \int_{\gamma_{sf}} \overrightarrow{\phi} \cdot \overrightarrow{g}_{sf} r d\gamma + \int_{\gamma_{out}} \overrightarrow{\phi} \cdot \overrightarrow{g}_{out} r d\gamma + \int_{\gamma_-} \overrightarrow{\phi} \cdot \overrightarrow{g}_- r d\gamma + \int_{\gamma_+} \overrightarrow{\phi} \cdot \overrightarrow{g}_+ r d\gamma \ ,
\end{split}
\end{equation}
is linear in $\mathbb{V}$. By using these two operators, the weak formulation of the mechanical problem can be expressed as follows: %\begin{problem}

\begin{itemize}
\item \textbf{Weak mechanical problem $\textbf{(WM)}$} : Let $T \in H^1_r(\omega)$ be the solution of the weak thermal model $(WT)$. Under the assumptions (MH1)-(MH3), find $\overrightarrow{u} \in \mathbb{U}$ such that,
\begin{equation}\label{bilinear_linear_form_momentum}
a_M(\overrightarrow{u},\overrightarrow{\phi}) = l_M[T](\overrightarrow{\phi}) \ , \ \forall \ \overrightarrow{\phi} \in \mathbb{U} \ .
\end{equation} \end{itemize}
%\michele{\nirav{Done} Nirav, use everywhere (W)T1 and (W)M1. I seem that for istance in the following they are missed. Thanks}
%The mechanical problem includes the temperature, which is the solution of the thermal model. We consider a segregated approach to address such a coupling: the temperature field is computed in advance, then it is inserted into the mechanical problem. % to compute the thermal
%stresses. %a one way coupling between the thermal and mechanical problem. Therefore, the mechanical problem is The mechanical problem includes the temperature, which is the solution of the thermal model.
%The last one is solved in advance independently.
%Therefore, a one way coupling between both models is being considered.

By using Lemmas 2.2 and 2.3 of \cite{axisymmetric_ivan}, under the assumption $(MH3)$, one can show that $a_M(\overrightarrow{u},\overrightarrow{\phi})$ and $l_M[T](\overrightarrow{\phi})$ are continuous in $\mathbb{V} \times \mathbb{V}$ and $\mathbb{V}$, respectively, and $a_M$ is $\mathbb{U}$-coercive.  %Lemma 2.2 of \cite{axisymmetric_ivan} shows that $a_M(\overrightarrow{u},\overrightarrow{\phi})$ is continuous in $\mathbb{V} \times \mathbb{V}$. Similarly, under assumptions $(MH1)$ and $(MH2)$, and assuming that $T \in H^1_r(\omega)$, Lemma 2.3 of \cite{axisymmetric_ivan} shows that $l_M[T](\overrightarrow{\phi})$ is continuous in $\mathbb{V}$ and it is bounded.
Hence the conditions of the Lax-Milgram theorem are satisfied and accordingly the weak mechanical formulation $(WM)$ has a unique solution. %$%\overrightarrow{u} \in \mathbb{U}$.}

Notice that here we could use the principle of superposition: the net displacement at any point in the domain $\overrightarrow{u}$ is the sum of the displacement due to purely mechanical effects $\overrightarrow{u}_M \in \mathbb{U}$ and the displacement due to purely thermal effects $\overrightarrow{u}_T \in \mathbb{U}$:

\begin{equation}
\overrightarrow{u} = \overrightarrow{u}_M + \overrightarrow{u}_T \ .
\end{equation}

Therefore, the problem (\emph{WM}) could be split in two sub-problems:

\begin{itemize}
\item \textbf{Weak mechanical problem $\textbf{(WM1)}$} : Under the assumptions (MH1)-(MH3), find $\overrightarrow{u}_M \in \mathbb{U}$ such that,
\begin{equation}\label{bilinear_linear_form_momentum}
a_M(\overrightarrow{u}_M,\overrightarrow{\phi}) = l_M[T_0](\overrightarrow{\phi}) \ , \ \forall \overrightarrow{\phi} \in \mathbb{U} \ , \\
\end{equation}
\item \textbf{Weak mechanical problem $\textbf{(WM2)}$} : Let $T \in H^1_r(\omega)$ be the solution of the weak thermal model $(WT)$. Under the assumptions (MH1)-(MH3), find $\overrightarrow{u}_T \in \mathbb{U}$ such that,
\begin{equation}\label{bilinear_linear_form_momentum}
a_M(\overrightarrow{u}_T,\overrightarrow{\phi}) = l_M[T](\overrightarrow{\phi}) - l_M[T_0](\overrightarrow{\phi}) \ , \ \forall \overrightarrow{\phi} \in \mathbb{U} \ , \\
\end{equation} \end{itemize}
%\begin{equation}\label{superposition_principle_weak_form}
%\begin{split}
%a_M(\overrightarrow{u}_M,\overrightarrow{\phi}) = l_M[T_0](\overrightarrow{\phi}) \ , \ \forall \overrightarrow{\phi} \in \mathbb{U} \ , \\
%a_M(\overrightarrow{u}_T,\overrightarrow{\phi}) = l_M[T](\overrightarrow{\phi}) - l_M[T_0](\overrightarrow{\phi}) \ , \ \forall \overrightarrow{\phi} \in \mathbb{U} \ . \\
%\overrightarrow{u} = \overrightarrow{u}_M + \overrightarrow{u}_T \ .
%\end{split}
%\end{equation}

\section[Finite element discretization]{Finite element discretization of the thermo-mechanical model}\label{fem_chapter}
%\michele{Let'us abbreviations: FE for finite element and MOR for model order reduction}
%We now look at the weak forms $(WT1)$ and $(WM1)$ introduced in Section \ref{weak_form_chapter}
In this section, we are going to briefly describe the Lagrange finite element method used to approximate the problems $(WT)$, and $(WM)$, $(WM1)$ and $(WM2)$. % at the aim to make suitable formulations for their numerical analysis and simulation.

We introduce the $n_h-$ dimensional space $H_{r,h}^1(\omega)\subset H_r^1(\omega)$ and $m_h-$ dimensional space ${\mathbb{U}}_h\subset\mathbb{U}$, %and $\mathbb{U}_h$, %of continuous spaces $H_r^1(\omega),\mathbb{U}$ respectively, i.e.,
%\begin{gather*}
%H_{r,h}^1(\omega) \subset H_r^1(\omega) \ , \\
%\mathbb{U}_h \subset \mathbb{U} \ .
%\end{gather*}
%The $n_h-$dimensional space $H_{r,h}^1(\omega)$ and $m_h-$dimensional space $\mathbb{U}_h$ are specified by,
\begin{gather}
H_{r,h}^1(\omega) = \spn \lbrace \psi_{1,h}, \psi_{2,h}, \dotsc, \psi_{n_h,h} \rbrace \ ,
\end{gather}
%and the $2n_h-$dimensional space $\mathbb{U}_h$ is specified by,
\begin{gather}
\mathbb{U}_h = \spn \lbrace \overrightarrow{\phi}_{1,h}, \overrightarrow{\phi}_{2,h}, \ldots, \overrightarrow{\phi}_{m_h,h} \rbrace \ .
\end{gather}
%$\lbrace \psi_{i,h} \rbrace_{i=1}^{n_h}$ and $\lbrace \overrightarrow{\phi}_{i,h} \rbrace_{i=1}^{2n_h}$ are the basis of the spaces $H_{r,h}^1(\omega)$ and $\mathbb{U}_h$, respectively.

%We seek the solutions $T_h \in H_{r,h}^1(\omega)$ and $\overrightarrow{u}_h \in \mathbb{U}_h$ of the discretized models corresponding to \michele{abbvs.} \nirav{$(WT1)$ and $(WM1)$}, respectively. %It is to be noted that $T_h \in H_{r,h}^1(\omega)$ and $\overrightarrow{u}_h \in \mathbb{U}_h$ are the approximation of $T \in H^1_r(\omega)$ and $\overrightarrow{u} \in \mathbb{U}$ respectively.

Based on the Galerkin method of weighted residuals \cite{finite_element_solids_book}, we can express the approximated solutions $T_h$ and $\overrightarrow{u}_h$ as follows:
\begin{gather}
T_h = \sum_{i=1}^{n_h} T_h^i \psi_{i,h} \ , \\
\overrightarrow{u}_h = \sum_{i=1}^{m_h} u_h^i \overrightarrow{\phi}_{i,h} \ ,
\end{gather}
where $T_h^i$ and $u_h^i$ are the nodal temperature and nodal displacement, respectively, and %, are obtained by solving linear systems of equations.
the basis functions $\psi_{i,h}$ and $\overrightarrow{\phi}_{i,h}$ are piecewise polynomial of degree $p \geq 1$ in $(r,y)$ space. %The nodal values, $\lbrace T_h^i \rbrace_{i=1}^{n_h}$ and $\lbrace u_h^i \rbrace_{i=1}^{2n_h}$, are common only to elements which are sharing the same node. Hence, the computational stencil of the system matrices does not extend beyond immediate neighboring elements and accordingly the system matrices are sparse matrices.

Therefore, the approximation of the problem $(WT)$ in the finite dimensional space $H^1_{r,h}$ can be stated as,

\begin{itemize}
\item \textbf{Problem $(WT)_h$} : Under the assumptions $(TH1)-(TH3)$, find $T_h \in H^1_{r,h}(\omega)$ such that,
\begin{equation}\label{fem_bilinear_thermal}
a_{T}(T_h,\psi_h) = l_{T}(\psi_h) \ , \ \forall \psi_h \in H^1_{r,h}(\omega) \ .
\end{equation} \end{itemize}

%The discrete model $(WT1)_h$ can be equivalently reformulated as a linear system:
%\begin{equation}\label{linear_form_discrete_thermal}
%\bm{A}_{T,h} \bm{T}_h = \bm{F}_{T,h}
%\end{equation}

Similarly, the approximation of the problems $(WM)$, $(WM1)$ and $(WM2)$ in the finite dimensional space $\mathbb{U}_h$ can be stated as:

\begin{itemize}
\item \textbf{Problem $(WM)_h$}: Let $T_h \in H^1_{r,h}(\omega)$, be the solution of the discretized thermal model $(WT)_h$. Under the assumptions $(MH1)-(MH3)$, find $\overrightarrow{u}_h \in \mathbb{U}_h$ such that,
%\end{problem}
\begin{equation}\label{fem_mechanical_equation}
a_{M}(\overrightarrow{u}_h,\overrightarrow{\phi}_h) = l_{M}[T_h](\overrightarrow{\phi}_h) \ , \ \forall \overrightarrow{\phi}_h \in \mathbb{U}_h \ .
\end{equation} \end{itemize}

\begin{itemize}
\item \textbf{Problem $(WM1)_h$}: Under the assumptions $(MH1)-(MH3)$, find $\overrightarrow{u}_{M,h} \in \mathbb{U}_h$ such that,
%\end{problem}
\begin{equation}\label{fem_mechanical_equation1}
a_{M}(\overrightarrow{u}_{M,h},\overrightarrow{\phi}_h) = l_{M}[T_0] (\overrightarrow{\phi}_h) \ , \ \forall \overrightarrow{\phi}_h \in \mathbb{U}_h \ .
\end{equation} \end{itemize}

\begin{itemize}
\item \textbf{Problem $(WM2)_h$}: Let $T_h \in H^1_{r,h}(\omega)$, be the solution of the discretized thermal model $(WT)_h$. Under the assumptions $(MH1)-(MH3)$, find $\overrightarrow{u}_{T,h} \in \mathbb{U}_h$ such that,
%\end{problem}
\begin{equation}\label{fem_mechanical_equation2}
a_{M}(\overrightarrow{u}_{T,h},\overrightarrow{\phi}_h) = l_{M}[T_h] (\overrightarrow{\phi}_h) - l_{M}[T_0] (\overrightarrow{\phi}_h) \ , \ \forall \overrightarrow{\phi}_h \in \mathbb{U}_h \ .
\end{equation} \end{itemize}

\section{Model order reduction}\label{parametric_rb_chapter}

In this section, we present our MOR framework. Firstly, we introduce the parameter space related to the problem under investigation (Sec. \ref{parspace}). Then, in Sec. \ref{rb_section}, we describe the POD algorithm that is used for the construction of reduced basis space as well as the two methods adopted for the computation of the reduced degrees of freedom, Galerkin projection (G) and Artificial Neural Network (ANN). Finally, in Sec. \ref{ROM_results}, we show some numerical tests with the aim to validate our approach. The MOR computations have been carried out using RBniCS \cite{RBniCS,RBniCS_website}, an in-house open source python library employing several reduced order techniques based on FEniCS \cite{FEniCS,FEniCS_website}, and PyTorch \cite{pytorch,PyTorch_website}, a python machine learning library. %based on the Torch library. %Finally, in \cite{}, we present some numerical tests in order to validate our approach.  %. We next discuss the construction of reduced basis space using Proper orthogonal decomposition. The reduced degrees of freedom are computed by two methods : Galerkin projection and Artificial Neural Network.

\subsection{Parameter space}\label{parspace}

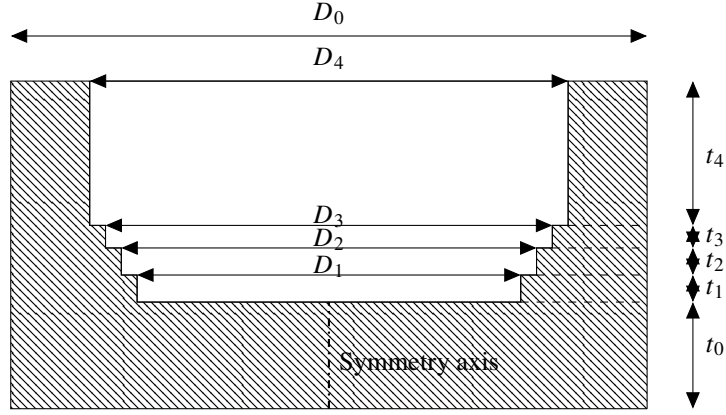
\begin{figure}
\centering
\begin{tikzpicture}[scale=0.6][
dot/.style={
  fill,
  circle,
  inner sep=2pt
  }
]
\draw[pattern=north west lines, pattern color=black] (7.05,0) -- (7.05,7.265) -- (5.3,7.265) -- (5.3,4.065) -- (4.95,4.065) -- (4.95,3.565) -- (4.6,3.565) -- (4.6,2.965) -- (4.25,2.965) -- (4.25,2.365) -- (-4.25,2.365) -- (-4.25,2.965) -- (-4.6,2.965) -- (-4.6,3.565) -- (-4.95,3.565) -- (-4.95,4.065) -- (-5.3,4.065) -- (-5.3,7.265) -- (-7.05,7.265) -- (-7.05,0) -- (7.05,0);
\draw[>=triangle 45, <->] (-7.05,8.265) -- (7.05,8.265);
\node at (0,8.765) {$D_0$};
\draw[>=triangle 45, <->] (8.065,2.365) -- (8.065,0);
\node at (8.565,1.38) {$t_0$};
\draw[dashed] (4.25,2.365) -- (7.05,2.365);
\draw[>=triangle 45, <->] (-5.3,7.265) -- (5.3,7.265);
\node at (0,7.765) {$D_4$};
\draw[>=triangle 45, <->] (8.065,2.965) -- (8.065,2.365);
\node at (8.565,2.665) {$t_1$};
\draw[dashed] (4.6,2.965) -- (7.05,2.965);
\draw[>=triangle 45, <->] (-4.25,2.965) -- (4.25,2.965);
\node at (0,3.165) {$D_1$};
\draw[>=triangle 45, <->] (8.065,2.965) -- (8.065,3.565);
\node at (8.565,3.265) {$t_2$};
\draw[dashed] (4.95,3.565) -- (7.05,3.565);
\draw[>=triangle 45, <->] (-4.6,3.565) -- (4.6,3.565);
\node at (0,3.765) {$D_2$};
\draw[>=triangle 45, <->] (8.065,3.565) -- (8.065,4.065);
\node at (8.565,3.815) {$t_3$};
\draw[dashed] (5.3,4.065) -- (7.05,4.065);
\draw[>=triangle 45, <->] (-4.95,4.065) -- (4.95,4.065);
\node at (0,4.265) {$D_3$};
\draw[>=triangle 45, <->] (8.065,4.065) -- (8.065,7.265);
\node at (8.565,5.565) {$t_4$};
\draw (-5.3,7.265) -- (-5.3,4.065) -- (-4.95,4.065) -- (-4.95,3.565) -- (-4.6,3.565) -- (-4.6,2.965) -- (-4.25,2.965) -- (-4.25,2.365) -- (4.25,2.365) -- (4.25,2.965) -- (4.6,2.965) -- (4.6,3.565) -- (4.95,3.565) -- (4.95,4.065) -- (5.3,4.065) -- (5.3,7.265) -- (-5.3,7.265);
\draw[thick, dash dot] (0,0) -- (0,2.365);
%\node at (0,6) {Volume of hearth};
\node at (2,1) {Symmetry axis};
\end{tikzpicture}
\caption{Hearth geometric parameters \revv{of $\omega$}\revv{: diameters $D_0,D_1,D_2,D_3,D_4$ and thicknesses $t_0,t_1,t_2,t_3,t_4$}.}
\label{hearth_parameters}
\end{figure}

Let $\mathbb{P} \subset \mathbb{R}^{d}$ be the parameter space having dimensionality $d$ with $\Xi \in \mathbb{P}$ a tuple of parameters. For the problem of the hearth blast furnace, the relevant parameters are related both to the physical properties and the geometry of the domain $\omega$. %As discussed in Sec. \ref{benchmark_test_section}, the domain is divided into $n_{su}$ triangular subdomains $\omega_i$ with $i = {1, 2, \dots, n_{su}}$ that represent a partition of $\omega$. In the context of geometric parametrization, domain $\omega$ is characterized by the geometric parameters.
The physical parameters are the thermal conductivity of the material, $k$, the thermal expansion coefficient, $\alpha$, the Young modulus, $E$, and the Poisson's ratio, $\nu$. On the other hand, the geometric parameters are the diameter of each section of the hearth $D_0,D_1,D_2,D_3, D_4$, and the thickness of each section of the hearth $t_0,t_1,t_2,t_3, t_4$ (see Figure \ref{hearth_parameters}). So for the problem under consideration, in the most general case (i.e., when all the parameters are considered), we have $\Xi = \{\Xi_p, \Xi_g\}$ where $\Xi_p = \{k, \alpha, E, \nu\} \in \mathbb{P}_p \subset \mathbb{R}^{4}$ is the physical parameters tuple and $\Xi_g = \{D_0, D_1, D_2, D_3, D_4, t_0, t_1, t_2, t_3, t_4\}  \in \mathbb{P}_g \subset \mathbb{R}^{10}$ is the geometric parameters tuple, and $d = 14$. %We note that the geometric parameters in the parameter tuple $\Xi$ uniquely identified the domain $\omega(\Xi)$.

Let us consider a geometrical parameters tuple $\bar{\Xi}_g$ % \in \mathbb{P}$ 
and the corresponding domain $\hat{\omega} = \omega(\bar{\Xi}_g)$. We refer to $\hat{\omega}$ 
as the reference domain. %As discussed in Sec. \ref{benchmark_test_section}, 
The domain $\hat{\omega}$ is divided into $n_{su}$ non-overlapping triangular subdomains (see, e.g., \cite{geometric_para_1, geometric_para_2}), i.e. 
\begin{equation}
\hat{\omega} = \bigcup\limits_{i=1}^{n_{su}} \hat{\omega}_i \ , \ \hat{\omega}_i \cap \hat{\omega}_j = \emptyset \ \text{for} \ i \neq j \ , \ 1 \leq i,j \leq n_{su} \ .
\end{equation}
For each of the subdomains $\hat{\omega}_i$, one can consider an invertible mapping $\bm{F}_i$,
\begin{equation}
\begin{split}
\bm{F}_i : \hat{\omega}_i \times \mathbb{P}_g \rightarrow \omega_i \ , %\ i = \lbrace 1,\ldots,n_{su} \rbrace,
\end{split}
\end{equation}
of the form
\begin{equation}\label{affine_F}
\overrightarrow{x} = \bm{F}_i(\overrightarrow{\hat{x}},\Xi_g) = \bm{G}_{F,i}(\Xi_g)\overrightarrow{\hat{x}} + \overrightarrow{c}_{F,i}(\Xi_g) ; \ , \ \forall \overrightarrow{\hat{x}} \in \hat{\omega}_i \ , \forall \overrightarrow{x} \in \omega_i (\Xi_g) \ , \\
\end{equation}
where
\begin{equation}
\begin{split}
\bm{G}_{F,i}=
  \begin{bmatrix}
    G_{F,i,11} & G_{F,i,12} \\
    G_{F,i,21} & G_{F,i,22}
  \end{bmatrix} \ ,
  \
  \overrightarrow{x} = \lbrace r \ y \rbrace^T \ , \ \overrightarrow{\hat{x}} = \lbrace \hat{r} \ \hat{y} \rbrace^T \ , \overrightarrow{c}_{F,i} = \lbrace c_{F,i,1} \ c_{F,i,2} \rbrace^T \ .
  \end{split}
\end{equation}

Equation \eqref{affine_F} highlights that the Jacobian matrix $\bm{G}_{F,i}$ and translation vector $\overrightarrow{c}_{F,i}$ are dependent only on the geometric parameters tuple $\Xi_g$ and do not vary over a given subdomain. In the following, the domains $\omega$ will be the image by eq. \eqref{affine_F} of the reference domain for the tuples of geometric parameters considered. % in $\mathbb{P}$.
% be the parameter space having dimensionality $d_p$. %We also introduce discrete and finite dimensional subset of parameter space, $\mathbb{P}_h \subset \mathbb{P}$. %Let $d_p$ be the dimensionality of the parameter space. Different sampling procedures can be used to sample parameter tuples which form the parameter space $\mathbb{P}_h$.
%\begin{equation*}
%\Xi \in \mathbb{P}_h \subset \mathbb{P} \subset \mathbb{R}^{d_p} \ .
%\end{equation*}
%\subsubsection{Geometric parametrization}\label{geometric_parameter_subsection}

%\nirav{The spaces $H^1_{r,h}(\omega)$ and $\mathbb{U}_h$, respectively, are the spaces of solution manifolds $T_h(\cdot;\Xi)$ and $\overrightarrow{u}(\cdot;\Xi)$ under variation of parameter tuple $\Xi \in \mathbb{P}$. The spaces $H^1_{r,h}(\hat{\omega})$ and $\hat{\mathbb{U}}_h$ correspond to the solution fields $T(\cdot,\bar{\Xi})$ and $\overrightarrow{u}(\cdot;\bar{\Xi})$, defined on reference domain $\hat{\omega}$, respectively.}

\subsection{Main ingredients of MOR}\label{rb_section}
The basic idea of MOR is the assumption that solutions live in a low dimensional manifold. Thus, any solution can be approximated based on a reduced number of global basis functions.

We seek the reduced basis approximations $T_h^{rb} \in H^{1,rb}_{r,h} (\omega)$ and $\overrightarrow{u}_h^{rb} \in \mathbb{U}_h^{rb}$ of $T_h \in H^1_{r,h} (\omega)$ and $\overrightarrow{u}_h \in \mathbb{U}_h$, respectively. The reduced basis spaces $H^{1,rb}_{r,h}(\omega) \subset H^1_{r,h}(\omega)$ and $\mathbb{U}_h^{rb} \subset \mathbb{U}_h$ are given by,
\begin{equation*}
\begin{split}
H^{1,rb}_{r,h} (\omega) = \spn \lbrace \psi_h^1, \ldots , \psi_h^{N_T} \rbrace \ , \\
\mathbb{U}_h^{rb} = \spn \lbrace \overrightarrow{\phi}_h^1, \ldots , \overrightarrow{\phi}_h^{N_M} \rbrace \ ,
\end{split}
\end{equation*}
where $N_M << m_h$ and $N_T << n_h$ are the number of basis functions forming the reduced basis spaces $H^{1,rb}_{r,h}(\omega)$ and $\mathbb{U}_h^{rb}$, respectively. Then we can represent $\overrightarrow{u}_h^{rb}$ and $T_h^{rb}$ by
\begin{equation}
T_h^{rb} = \sum_{i=1}^{N_T} \zeta_{T}^i \psi_h^i \ , \\
%\overrightarrow{u}_h = \sum_{i=1}^{N_T}  \overrightarrow{\phi}_{i,h} \ .
\end{equation}
\begin{equation}
\overrightarrow{u}_h^{rb} = \sum_{i=1}^{N_M} \zeta_{M}^i \overrightarrow{\phi}_h^i \ , \\
%\overrightarrow{u}_h = \sum_{i=1}^{N_M}  \overrightarrow{\phi}_{i,h} \ .
\end{equation}
where $\lbrace \zeta_{T}^i \rbrace_{i=1}^{N_T}$ and $\lbrace \zeta_{M}^i \rbrace_{i=1}^{N_M}$ are the temperature and displacement degrees of freedom, respectively.

We also construct the reduced basis spaces for displacement fields $\overrightarrow{u}_{T,h}$ and $\overrightarrow{u}_{M,h}$, introduced in \eqref{fem_mechanical_equation1} and \eqref{fem_mechanical_equation2}, as,
\begin{equation*}
\begin{split}
\mathbb{U}_{T,h}^{rb} = \spn \lbrace \overrightarrow{\phi}_{T,h}^1, \ldots , \overrightarrow{\phi}_{T,h}^{N_{M,T}} \rbrace \ , \\
\mathbb{U}_{M,h}^{rb} = \spn \lbrace \overrightarrow{\phi}_{M,h}^1, \ldots , \overrightarrow{\phi}_{M,h}^{N_{M,M}} \rbrace \ .
\end{split}
\end{equation*}

So, the reduced basis approximations $\overrightarrow{u}_{M,h}^{rb} \in \mathbb{U}_{M,h}^{rb}$ of $\overrightarrow{u}_{M,h}\in \mathbb{U}_h$ and $\overrightarrow{u}_{T,h}^{rb} \in \mathbb{U}_{T,h}^{rb}$ of $\overrightarrow{u}_{T,h}\in \mathbb{U}_h$ can be represented as,
\begin{equation}
\overrightarrow{u}_{M,h}^{rb} = \sum_{i=1}^{N_{M,M}} \zeta_{M,M}^i \overrightarrow{\phi}_{M,h}^i \ , \\
\overrightarrow{u}_{T,h}^{rb} = \sum_{i=1}^{N_{M,T}} \zeta_{M,T}^i \overrightarrow{\phi}_{T,h}^i \ .
\end{equation}

\subsubsection{POD algorithm}\label{pod_section}

In the literature, one can find several techniques to generate the reduced basis spaces, e.g. Proper Orthogonal Decomposition (POD), the Proper Generalized Decomposition (PGD) and the Reduced Basis (RB) with a greedy sampling strategy. See, e.g., \cite{mor_bader,mor_book_benner3,haasdonk_chapter,RBniCS,mor_book_quarteroni}. %\michele{Nirav, in the following some citations shoud be added here, together to ones already present (such [19]), thanks!}.
%REFERENCES%
% G. Rozza, D. B. P. Huynh, A. T. Patera, Reduced Basis Approximation and a Posteriori Error Estimation for Affinely Parametrized Elliptic Coercive Partial Differential Equations, Archives of Computational Methods in Engineering 15 (3) (2008) 229. doi:10.1007/s11831-008-9019-9.
% F. Chinesta, A. Huerta, G. Rozza, K. Willcox, Model Order Reduction, Encyclopedia of Computational Mechanics, Elsevier Editor.
% I. Kalashnikova, M. F. Barone, On the stability and convergence of a Galerkin reduced order model (ROM) of compressible flow with solid wall and far-field boundary treatment, International Journal for Numerical Methods in Engineering 83 (10) (2010) 1345–1375.
% F. Chinesta, P. Ladeveze, E. Cueto, A Short Review on Model Order Reduction Based on Proper Generalized Decomposition, Archives of Computational Methods in Engineering 18 (4) (2011) 395.
%A. Dumon, C. Allery, A. Ammar, Proper General Decomposition (PGD) for the resolution of Navier-Stokes equations, Journal of Computational Physics 230 (4) (2011) 1387–1407.
%V. Tsiolakis, M. Giacomini, R. Sevilla, C. Othmer, A. Huerta, Parametric solutions of turbulent incompressible flows in openfoam via the proper generalised decomposition.
%A. Quarteroni, A. Manzoni, F. Negri, Reduced Basis Methods for Partial Differential Equations, Springer International Publishing, 2016.
%P. Benner, W. Schilders, S. Grivet-Talocia, A. Quarteroni, G. Rozza, L. M. Silveira, Model Order Reduction, De Gruyter, Berlin, Boston, 2020.
In this work, the reduced basis spaces are constructed by POD that is able to capture the ``dominant'' modes by exploiting the information contained in the full order snapshots. %. POD exploits the information contained in the snapshots to construct the low dimensional reduced basis space which can approximate the solution within desirable accuracy.

We are going to describe the procedure for the computation of the reduced basis space $H^{1,rb}_{r,h} (\omega)$. The reduced basis spaces $\mathbb{U}_h^{rb}$, $\mathbb{U}_{T,h}^{rb}$ and $\mathbb{U}_{M,h}^{rb}$ are constructed in an analogous way.  First, $n_s$ parameter tuples, $\lbrace \Xi_k \rbrace_{k=1}^{n_s}$, are considered that form the training set. %and training set $\lbrace \Xi_k \rbrace_{k=1}^{n_s}$ is formed. 
We compute the snapshots $T_h(\Xi_k)$ related to each parameter tuple in the training set. Then a matrix $\bm{C}_T \in \mathbb{R}^{n_s \times n_s}$ is constructed,
\begin{equation}
\left( \bm{C}_T \right)_{kl} = <T_h(\Xi_k),T_h(\Xi_l)>_{H^1_{r,h}(\hat{\omega})} \ , \ 1 \leq k,l \leq n_s \ ,
\end{equation}
where transformation \eqref{affine_F} is considered. Next, $N_T$ largest eigenvalues $\lbrace \theta_{T}^{i} \rbrace_{i=1}^{N_T}$ of the  matrix $\bm{C}_T$, sorted in descending order, $\theta_{T}^{1} \geq \theta_{T}^{2} \geq \ldots \geq \theta_{T}^{N_T}$, and corresponding eigenvectors $\lbrace \bm{V}_{T}^{i} \rbrace_{i=1}^{N_T} , \bm{V}_{T}^{i} \in \mathbb{R}^{n_s}$, are computed,
\begin{equation}
\bm{C}_T \bm{V}_{T}^{i} = \theta_{T}^{i} \bm{V}_{T}^{i} \ .
\end{equation}
%It is assumed that the eigenvalues are sorted in descending order, $\theta_{T}^{1} \geq \theta_{T}^{2} \geq \ldots \geq \theta_{T}^{N_T}$.
The reduced basis are then given by
\begin{gather}\label{pod_basis_1}
\psi_h^{i} = \frac{\sum\limits_{k=1}^{n_s} (\bm{V}_{T}^{i})_{k} T_h(\Xi_k)}{||\sum\limits_{k=1}^{n_s} (\bm{V}_{T}^{i})_{k} T_h(\Xi_k) ||_{H^1_{r,h}(\hat{\omega})}}  \ , \ 1 \leq i\leq N_T \ .
%H^{1,rb}_{r,h}(\omega) = \lbrace \psi_h^{i} \rbrace_{i=1}^{N_T} \ .
\end{gather}

In order to determine the admissibility of a given eigenvector into the POD space, we refer to the following criterion % \cite{} \michele{let's insert a relevant citation} \nirav{Sorry, Michele but I do not have references anymore and I seriously do not think we can afford to delay our work for adding more references. There are many references now in the bibliography and it is time that we move onto this.} \michele{Nirav, it is not necessary added NEW references, it's sufficient report here some relevant references already cited (ones where this common criterion is proposed)}: %The eigenvalue related to a given eigenvector is used to determine admissibility of the eigenvector itself into the POD space. The relative eigenvalue $\theta_{T,rel}^i$ is defined as,
\begin{equation}\label{eigenvalue_pod_criteria}
%\theta_{T,rel}^i =
\frac{\theta_T^i}{\theta_T^1}  \geq 1e-4 \ .
\end{equation}

%\michele{Maybe the following part is not necessary...i see after}

%\michele{Nirav, if we do not need more of $\bm{B}$ i would comment the following part}

%We observe that each reduced basis \eqref{pod_basis_1} could be expressed by means of a linear combination of basis belonging to the discretized finite element space:
%\begin{equation}
%\psi_h^{i} = \bm{B}_{T,ij} \psi_{j,h} \ ,% \ \bm{B}_T \in \mathbb{R}^{N_T \times n_h} \ , \ \psi_{j,h} \in H_{r,h}^1(\omega) \ , \ \psi_h^{i} \in H_{r,h}^{1,rb}(\omega) \ ,
%\end{equation}
%where $\bm{B}_{T,i}$ is a \michele{report the dimensions of such a matrix} matrix containing proper weights.

%\nirav{We repeat above algorithm using snapshots of $\overrightarrow{u}_h, \overrightarrow{u}_{M,h}$ and $\overrightarrow{u}_{T,h}$ to construct the operators $\bm{B}_u, \bm{B}_{u,M}$ and $\bm{B}_{u,T}$ and corressponding reduced basis spaces $\mathbb{U}_h^{rb},\mathbb{U}_{M,h}^{rb}$ and $\mathbb{U}_{T,h}^{rb}$ respectively.}

%The reduced basis space $\overrightarrow{u}_h^{rb}$ is similarly constructed.

In the following subsections, we describe two different approaches for the computation of the degrees of freedom: Galerkin projection (G) and Artificial Neural Network (ANN).

\subsubsection{Galerkin projection}

Here we choose to consider $(WM1)_h$ and $(WM2)_h$ in the place of $(WM)_h$ because \revv{the displacements due to purely mechanical effects and those due to purely thermal effects} may have different \revv{characteristic} scales. \revv{So treating these systems separately ensures that the system with bigger scale effects does not dominate the global system, avoiding a significant impact on the accuracy of MOR}: see, e.g.,  \cite{greedy_paper,mor_book,variational_multiscale}. %So the models $(WT)_h$, $(WM1)_h$ and $(WM2)_h$ are projected on the reduced basis spaces in order to obtain an algebraic equations system whose unknowns are the degrees of freedom.
% \nirav{because treating these systems separately ensures that the system with bigger scale effects does not dominate the system and the reduced basis space also involves representation of the system with smaller scale effect.}

We consider an affine parametric dependence, i.e. the bilinear forms $a_T(\cdot,\cdot;\Xi)$ and $a_M(\cdot,\cdot;\Xi)$ are expressed as weighted sum of $n_{a_T}$ and $n_{a_M}$ parameter independent bilinear forms. Similarly, the linear forms $l_T(\cdot;\Xi)$ and $l_M[T](\cdot;\Xi)$ are expressed as weighted sum of $n_{l_T}$ and $n_{l_M}$ parameter independent linear forms. We have:

\begin{equation}\label{eq:affine1}
\begin{split}
a_T(T,\psi;\Xi) = \sum_{i=1}^{n_{a_T}} \theta_{a_{T,i}}(\Xi) a_{T,i}(T,\psi;\bar{\Xi}) \ , \\
l_T(\psi;\Xi) = \sum_{i=1}^{n_{l_T}} \theta_{l_{T,i}} (\Xi) l_{T,i}(\psi;\bar{\Xi}) \ ,
\end{split}
\end{equation}
and
\begin{equation}\label{eq:affine2}
\begin{split}
a_M(\overrightarrow{u},\overrightarrow{\phi};\Xi) = \sum_{i=1}^{n_{a_M}} \theta_{a_{M,i}}(\Xi) a_{M,i}(\overrightarrow{u},\overrightarrow{\phi};\bar{\Xi}) \ , \\
l_M[T](\overrightarrow{\phi};\Xi) = \sum_{i=1}^{n_{l_M}} \theta_{l_{M,i}} (\Xi) l_{M,i}[T](\overrightarrow{\phi};\bar{\Xi}) \ .
\end{split}
\end{equation}

%In other words, the bilinear forms $a_T(\cdot,\cdot;\Xi)$ and $a_M(\cdot,\cdot;\Xi)$ are expressed as weighted sum of $n_{a_T}$ and $n_{a_M}$ parameter independent bilinear forms. Similarly, the linear forms $l_T(\cdot;\Xi)$ and $l_M[T](\cdot;\Xi)$ are expressed as weighted sum of $n_{l_T}$ and $n_{l_M}$ parameter independent linear forms.
The affine expansion of operators is essentially a change of variables and has been widely addressed in the literature: see, e.g., \cite{mor_book_benner3,haasdonk_chapter,RBniCS}. %\nirav{By introducing the geometric parametrization (section \ref{parspace}) in the operators, the problem defined on the parametrized domain $\omega$ can be recast on the reference domain $\hat{\omega}$. The effect of variation in geometry is incorporated through the parametrized transformations.}
The affinity assumption is particularly important as it leads to considerable efficiency. This is mainly due to the fact that the evaluation of bilinear forms, $a_{M,i}(\overrightarrow{u},\overrightarrow{\phi})$ and $a_{T,i}(T,\psi)$, and linear forms, $l_{T,i}(\psi)$ and $l_{M,i}[T](\overrightarrow{\phi})$ are not required for each new tuple of parameters.

So for what concerns the model $(WT)_h$, the bilinear form $a_{T,h}:H^1_{r,h}(\omega) \times H^1_{r,h}(\omega) \rightarrow \mathbb{R}$ is restricted to the reduced basis space as $a_{T,h}^{rb} : H^{1,rb}_{r,h}(\omega) \times H^{1,rb}_{r,h}(\omega) \rightarrow \mathbb{R}$. In the same way, the linear form $l_{T,h}:H^1_{r,h}(\omega) \rightarrow \mathbb{R}$ is restricted to the reduced basis space as $l_{T,h}^{rb} : H^{1,rb}_{r,h}(\omega) \rightarrow \mathbb{R}$. So the reduced basis approximation $T_h^{rb}$ at a given parameter tuple $\Xi^*$ is obtained by solving
\begin{equation}\label{reduced_weak_form_temperature}
a_{T,h}^{rb}(T_h^{rb},\psi_h^{rb};\Xi^*) = l_{T,h}^{rb}(\psi_h^{rb};\Xi^*) \ , \ \forall \psi_h^{rb} \in H^{1,rb}_{r,h}(\omega) \ .
\end{equation}

%We introduce the coercivity constant $c_{T,h}^{rb} > 0$ and continuity constant $C_{T,h}^{rb} > 0, \forall \psi_h^{rb}, T_h^{rb} \in H^{1,rb}_{r,h} (\omega)$.
%\begin{gather}
%c_{T,h}^{rb} ||\psi_h^{rb}||^2_{H^1_{r,h}(\omega)} \leq a_{T,h}^{rb}(\psi_h^{rb},\psi_h^{rb}) \ , \\
%a_{T,h}^{rb}(T_h^{rb},\psi_h^{rb}) \leq C_{T,h}^{rb} ||T_h^{rb}||_{H^1_{r,h}(\omega)} ||\psi_h^{rb}||_{H^1_{r,h}(\omega)}  \ .
%\end{gather}

%The reduced basis approximation $T_h^{rb}$ is orthogonal projection of finite element approximation $T_h$ with respect to $a_{T,h}(\cdot,\cdot)$,
%\begin{equation}
%a_{T,h}(T_h - T^{rb}_h,\psi_h^{rb}) = 0 \ , \ \forall \psi_h^{rb} \in H^{1,rb}_{r,h}(\omega) \ .
%\end{equation}

%By following procedure for deriving error estimate \eqref{error_estimate_h1r_temperature}, following error estimate can be derived,
%\begin{equation}\label{error_estimate_h1r_temperature_rb}
%||T_h-T_h^{rb}||_{H^1_{r,h}(\omega)} \leq \sqrt{\frac{C_{T,h}+c^{rb}_{T,h}-c_{T,h}}{c^{rb}_{T,h}}} ||T_h-\psi^{rb}_h||_{H^1_{r,h}(\omega)} \ , \ \forall \psi^{rb}_h \in H^{1,rb}_{r,h}(\omega) \ .
%\end{equation}

%\michele{Nirav, in the POD-G method you consider the same displacement basis for the following sub problems? If yes, they are referred to the POD applied to mechanical snapshots (the ones related to the first subproblem?)} \nirav{No. I performed POD spearately and hence, different spaces were constructed. I added 2 paragraphs below. I would request to check them and remove this paragraph.}

On the other hand, for what concerns the model $(WM1)_h$, the bilinear form $a_{M,h}:\mathbb{U}_h \times \mathbb{U}_h \rightarrow \mathbb{R}$ is restricted to the reduced basis space as $a_{M,h}^{rb} : \mathbb{U}_{M,h}^{rb} \times \mathbb{U}_{M,h}^{rb} \rightarrow \mathbb{R}$. The linear form $l_{M,h}:\mathbb{U}_h \rightarrow \mathbb{R}$ is restricted to the reduced basis space as $l_{M,h}^{rb} : \mathbb{U}_{M,h}^{rb} \rightarrow \mathbb{R}$. So the reduced basis approximation $\overrightarrow{u}_{M,h}^{rb}$ at a given parameter tuple $\Xi^*$ is obtained by solving,
\begin{equation}\label{reduced_weak_form_mechanical1}
a_{M,h}^{rb}(\overrightarrow{u}_{M,h}^{rb},\overrightarrow{\phi}_{M,h}^{rb};\Xi^*) = l_{M,h}^{rb}[T_0] (\overrightarrow{\phi}_{M,h}^{rb};\Xi^*) \ , \ \forall \overrightarrow{\phi}_{M,h}^{rb} \in \mathbb{U}_{M,h}^{rb} \ .
\end{equation}

Similarly, for the model $(WM2)_h$ the bilinear form $a_{M,h}:\mathbb{U}_h \times \mathbb{U}_h \rightarrow \mathbb{R}$ is restricted to the reduced basis space as $a_{M,h}^{rb} : \mathbb{U}_{T,h}^{rb} \times \mathbb{U}_{T,h}^{rb} \rightarrow \mathbb{R}$. The linear form $l_{M,h}:\mathbb{U}_h \rightarrow \mathbb{R}$ is restricted to the reduced basis space as $l_{M,h}^{rb} : \mathbb{U}_{T,h}^{rb} \rightarrow \mathbb{R}$. So, the reduced basis approximation $\overrightarrow{u}_{T,h}^{rb}$ at a given parameter $\Xi^*$ is obtained by solving,

\begin{equation}\label{reduced_weak_form_mechanical2}
a_{M,h}^{rb}(\overrightarrow{u}_{T,h}^{rb},\overrightarrow{\phi}_{T,h}^{rb};\Xi^*) = l_{M,h}^{rb}[T_h^{rb}] (\overrightarrow{\phi}_{T,h}^{rb};\Xi^*) - l_{M,h}^{rb}[T_0] (\overrightarrow{\phi}_{T,h}^{rb};\Xi^*) \ , \ \forall \overrightarrow{\phi}_{T,h}^{rb} \in \mathbb{U}_{T,h}^{rb} \ .
\end{equation}

\subsubsection{Artificial Neural Network}\label{ann_intro}

%The second method that we use for the computation of $\lbrace \zeta_{T}^i \rbrace_{i=1}^{N_T}$ and $\lbrace \zeta_{M}^i \rbrace_{i=1}^{N_M}$ is based on

Artificial Neural Network (ANN) is a computational model that takes inspiration from the human brain consisting of an interconnected network of simple processing units that can learn from experience by modifying their connections (see, e.g., \cite{haykin_ann_book,neural_network_book_raul}). %Once a proper training phase is performed, ANN is able to accurately predict the output from a given set of input. %, i.e. to mimic the functional relationship between input and output by ``learning'' from data. %without explicitly using the analytical relationship between input and output. In other words, ANN is used to mimic the functional relationship between input and output by ``learning'' from data. Here, we consider a feed-forward neural network whose scheme is displayed in Figure \ref{neural_network_schematic}. In this kind of network, the information flows from left to right.

Recently, the application of deep learning methods to partial differential equations has shown promising capabilities: see, e.g.,  \cite{gkn_paper,fno_paper,pinn_paper,prnn_paper,raissi_velocity_pressure_flows,sinha_paper,balla_paper}. Concerning the application of the ANN approach in a MOR context, the reader is referred, e.g., to \cite{pod_nn_paper,ann_mor_transient_flow,non_intrusive_1,prnn_paper, federico_paper,wang_ann_paper,demo_strazzullo_paper,meneghetti_paper}. %\michele{Nirav,  here insert also the other ROM NN references that we have cited. Moreover, inert also the following reference (also in the introduction), thanks!}\nirav{Done}
%REFERENCES
%https://www.sciencedirect.com/science/article/abs/pii/S1007570419301364
%Here, we test its capabilities in a linear framework.
We highlight that, unlike the Galerkin projection, ANN is a data-driven approach, i.e. based only on data and does not
require the knowledge of the original equations describing the system. It is also non-intrusive, in the sense that no modification of the simulation software is required. %it is a non-intrusive approach, i.e. eliminates the need to assemble and solve any system for reduced basis approximation.
\begin{figure}
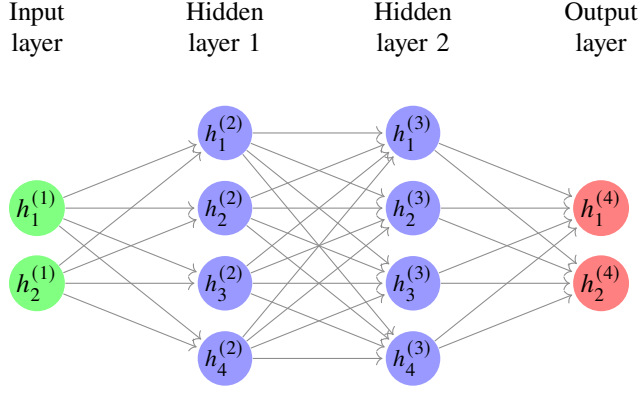

\begin{center}
\begin{neuralnetwork}[height=5]
    \newcommand{\x}[2]{$h^{(1)}_#2$}
    \newcommand{\y}[2]{$h^{(4)}_#2$}
    \newcommand{\hfirst}[2]{\small $h^{(2)}_#2$}
    \newcommand{\hsecond}[2]{\small $h^{(3)}_#2$}
    \inputlayer[count=2, bias=false, title=Input\\layer, text=\x]
    \hiddenlayer[count=4, bias=false, title=Hidden\\layer 1, text=\hfirst] \linklayers
    \hiddenlayer[count=4, bias=false, title=Hidden\\layer 2, text=\hsecond] \linklayers
    \outputlayer[count=2, title=Output\\layer, text=\y] \linklayers
\end{neuralnetwork}
%\caption{Sketch of a feed-forward ANN with $n_l = 4$.}
\caption{Sketch of a feed-forward ANN with \revv{$2$ hidden layers}.}
\label{neural_network_schematic}
\end{center}
\end{figure}

In this work, we use a feed-forward ANN consisting of input layer, two hidden layers \revvv{whose depth $H$ (i.e., the number of neurons constituting the hidden layer)} is determined by trial and error \cite{pod_nn_paper}, and output layer. %with $n_l$ the total number of layers and $d_l$ the number of unit cells (the so-called neurons) of the $l-$th layer. The neurons of each layer are connected to the neurons of the next layer by synapses. 
See Figure \ref{neural_network_schematic} for an illustrative representation of a feed-forward ANN. The weights as well as the biasing parameters of the network are iteratively adjusted by the backpropagation process using an optimization algorithm \cite{adam_reference}. Concerning the activation function, we consider for hidden layers the Sigmoid function \cite{ann_forecasting} whilst for the initial and final layer we use the identity function. 

Unlike what done for the Galerkin projection approach, concerning the mechanical problem, we consider the model $(WM)_h$. % and assess whether the scale effects are taken care by the ANN.
This choice is due to the fact that ANN suffers from high offline cost because of the training phase, so it is beneficial to train one only model instead of training two models. 

We consider $N_t^T$ parameter tuples $\lbrace \Xi_k \rbrace_{k=1}^{N_t^T}$ and compute the temperature field $T_h(\Xi_k)$ by solving problem $(WT1)_h$ at each parameter tuple $\Xi_k$. Next, the temperature field $T_h(\Xi_k)$ is projected on the reduced basis space so to obtain the projected solution $T_h^{\pi} (\Xi_k)$ and corresponding degrees of freedom $\bm{\zeta}_{T,\pi}(\Xi_k)$:
\begin{equation}\label{temperature_solution_projection}
T_h^{\pi} (\Xi_k) = \argmin_{\psi_h^{rb} \in H^{1,rb}_{r,h}(\omega)} ||T_h(\Xi_k) - \psi_h^{rb}||_{H^1_{r,h}(\omega)} = \sum\limits_{i=1}^{N_T} \zeta_{T,\pi}^i(\Xi_k) \psi_h^i \ . % \ \bm{\zeta}_{T,\pi}(\Xi_k) = \lbrace \zeta_{T,\pi}^i(\Xi_k) \rbrace_{i=1}^{N_T} \ .
\end{equation}

Similarly, we consider $N_t^M$ parameter tuples $\lbrace \Xi_k \rbrace_{k=1}^{N_t^M}$ and compute the displacement fields $\overrightarrow{u}_h(\Xi_k)$ by solving problem $(WM)_h$ at each parameter tuple $\Xi_k$. Next, the displacement field $\overrightarrow{u}_h(\Xi_k) \in \mathbb{U}_h$ is projected on the reduced basis space so to obtain the projected solution $\overrightarrow{u}_h^{\pi} (\Xi_k)$ and corresponding degrees of freedom $\bm{\zeta}_{M,\pi}(\Xi_k)$:
\begin{equation}\label{displacement_solution_projection}
\overrightarrow{u}_h^{\pi} (\Xi_k) = \argmin_{\overrightarrow{\phi}_h^{rb} \in \mathbb{U}_h^{rb}} ||\overrightarrow{u}_h(\Xi_k) - \overrightarrow{\phi}_h^{rb}||_{\mathbb{U}_h} = \sum\limits_{i=1}^{N_M} \zeta_{M,\pi}^i (\Xi_k) \overrightarrow{\phi}_h^i \ . % \ \bm{\zeta}_{M,\pi}(\Xi_k) = \lbrace \zeta_{M,\pi}^i(\Xi_k) \rbrace_{i=1}^{N_M} \ .
\end{equation}

So we consider two collections of (known) training input-desidered output pairs, $\lbrace \Xi_k, \zeta_{T,\pi}(\Xi_k) \rbrace_{k=1}^{N_t^T}$ and $\lbrace \Xi_k , \zeta_{M,\pi}(\Xi_k) \rbrace_{k=1}^{N_t^M}$. The goal is to approximate the functions $f_T$ and $f_M$ that map these training input-desidered output pairs. %and use the functions $f_T$ and $f_M$, for the reduced approximations $T_h^{rb}$ and $\overrightarrow{u}_h^{rb}$, respectively.
After training the two ANNs, we consider them as black boxes that can then be used to compute the POD coefficients related to a new parameter instance $\Xi^*$.

%\begin{equation}
%h_i^{(1)} = \Xi^* \ , \ h_i^{(n_l)} = \zeta_T \ , \ f_T : \mathbb{R}^{d_p} \rightarrow \mathbb{R}^{N_T} \ , \ \zeta_T(\Xi^*) = f_T ({\Xi^*}) \ , \ \zeta_T = \lbrace \zeta_T^i \rbrace_{i=1}^{N_T} \ .
%\end{equation}

%\begin{equation}
%h_i^{(1)} = \Xi^* \ , \ h_i^{(n_l)} = \zeta_M \ , \ f_M : \mathbb{R}^{d_p} \rightarrow \mathbb{R}^{N_M} \ , \ \zeta_M(\Xi^*) = f_M ({\Xi^*}) \ , \ \zeta_M = \lbrace \zeta_M^i \rbrace_{i=1}^{N_M} \ .
%\end{equation}

%In general, the training phase of ANN is an iterative process and requires balancing training error and generalisation error \cite{haykin_ann_book}.

We split the full order data into two parts: one to be used for training and one to be used for validation. While the training data are used to adjust weights and biasing parameters of the ANN, the validation data are used to measure its accuracy. A common issue is that ANN may perform better on training data but may not perform well on data other than training data. To avoid this overfitting phenomenon, we use the early stopping criteria \cite{haykin_ann_book}: the training is stopped when the the mean squared error 
\begin{equation}\label{mse_error}
\begin{split}
\epsilon_T = \frac{\sum\limits_{i=1}^{N_T} \left( \zeta^i_{T,\pi} (\Xi_k) - \zeta^i_T (\Xi_k) \right)^2}{N_T} \ , \\
\epsilon_M = \frac{\sum\limits_{i=1}^{N_M} \left( \zeta^i_{M,\pi} (\Xi_k) - \zeta^i_M(\Xi_k) \right)^2}{N_M} \ ,
\end{split}
\end{equation}

% $\epsilon_T$ and $\epsilon_M$ at given parameter $\Xi_k$.
% \begin{equation}\label{mse_error}
% \begin{split}
% \epsilon_T = \frac{\sum\limits_{i=1}^{N_T} \left( \zeta^i_{T,\pi} (\Xi_k) - \zeta^i_T (\Xi_k) \right)^2}{N_T} \ , \\
% \epsilon_M = \frac{\sum\limits_{i=1}^{N_M} \left( \zeta^i_{M,\pi} (\Xi_k) - \zeta^i_M(\Xi_k) \right)^2}{N_M} \ ,
% \end{split}
% \end{equation}
as measured on validation data starts to increase.

\subsection{\revv{MOR results}}\label{ROM_results}
%In this section we test the performances of the two MOR frameworks presented in Sec. \ref{parametric_rb_chapter}, POD-G and POD-ANN.

%\subsubsection{Design of experiment}
%We explain now the mesh considered and the range of parameter values.
The coordinates  of the 12 vertices constituting the reference domain %illustrated in Figure %\ref{hearth_2d_domain} 
are reported in Table \ref{domain_geometry_table}. %The domain considered in this section % corresponds to a real furnace, whose vertical section is a polygon in the $r-y$ plane (see Figure \ref{hearth_2d_domain}). Its geometrical features are reported in Table \ref{domain_geometry_table}.
\begin{comment}
\begin{table}[ht]
\centering
\begin{tabular}{|c|}
\hline
{($r$ ,$y$)} [m] \\%\textbf{Variables}
\hline
 (0   , 0    )\\
 (7.05, 0    ) \\
 (7.05, 7.265)\\
 (5.30, 7.265) \\
 (5.30, 4.065)\\
 (4.95, 4.065) \\
 (4.95, 3.565)\\
 (4.6 , 3.565) \\
 (4.6 , 2.965)\\
 (4.25, 2.965) \\
 (4.25, 2.365)\\
 (0   , 2.365) \\
\hline
%Maximum radius & $r_{max}$ & $7.05$ \\
%Maximum height & $y_{max}$ & $7.265$ \\
\end{tabular}
\caption{Vertices coordinates of the domain $\omega$ (see Figure \ref{hearth_2d_domain}).}
\label{domain_geometry_table}
\end{table}
\end{comment}
\begin{comment}
\begin{table}[ht]
\centering
\begin{tabular}{|c|c|}
\hline
$r$ & $y$ \\
\hline
0 & 0 \\
\hline
0 & 0 \\
\hline
7.05 & 0 \\
\hline
7.05 & 7.265\\
\hline
5.30 & 7.265 \\
\hline
5.30 & 4.065\\
\hline
4.95 & 4.065\\
\hline
4.95 & 3.565\\
\hline
4.6 & 3.565\\
\hline
4.6 & 2.965\\
\hline
4.25 & 2.965\\
\hline
4.25 & 2.365\\
\hline
0 & 2.365 \\
\hline
\end{tabular}
\caption{coordinates (in m) of the vertices of domain $\omega$ (see Figure \ref{hearth_2d_domain}).}
\label{domain_geometry_table}
\end{table}
\end{comment}

\begin{table}[ht]
\centering
\begin{tabular}{|c|c|c|c|c|c|c|c|c|c|c|c|c|}
\hline
$r$ & 0 & 7.05 & 7.05 & 5.30 & 5.30 & 4.95 & 4.95 & 4.6 & 4.6 & 4.25 & 4.25 & 0 \\
\hline
$y$ & 0 & 0 & 7.265 & 7.265 & 4.065 & 4.065 & 3.565 & 3.565 & 2.965 & 2.965 & 2.365 & 2.365 \\
\hline
\end{tabular}
\caption{coordinates (in m) of the vertices of domain \revv{$\hat{\omega}$} (see Figure \revv{\ref{mesh_details}}).}
\label{domain_geometry_table}
\end{table}

%\begin{center}

\begin{figure}
\centering
\begin{subfigure}{0.45\textwidth}
\centering
\includegraphics[height=2in]{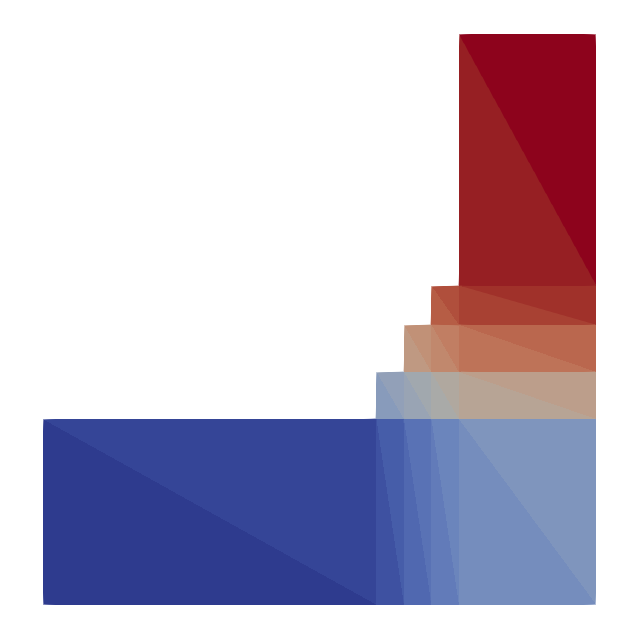}
\caption{Subdomains decomposition.}
\label{subdomains_ch4}
\end{subfigure}
\begin{subfigure}{0.45\textwidth}
\centering
\includegraphics[height=2in]{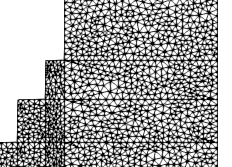}
\caption{Close up of the mesh.}
\label{mesh_elements}
\end{subfigure}
%\begin{subfigure}{0.33\textwidth}
\centering
%\includegraphics[height=1.5in]{images/mesh_quality.png}
%\caption{Condition number.}%\michele{I would change the color bar style (try "rainbow") and make the label greater}\nirav{Done except "Rainbow" style. I think this is good. Do you think so?} \michele{No because the scale [1,4] does not allows to appreciate anything (it seems a unique colour). Let's try to do better (for istance you can get down the maximum value of the reference scale, try 2)}\nirav{Revised. Bringing maximum value down to 2 will hide elements with worst condition number.}}
%\label{mesh_quality}
%\end{subfigure}
%\begin{subfigure}{0.45\textwidth}
%\centering
%\includegraphics[height=1.5in]{images/mesh.png}
%\caption{Mesh \michele{Nirav, here we should find a way to show well the mesh (now the elements are not visible). We can try to show only a part of the mesh for istance} \nirav{See figure  \ref{mesh_elements}}}%$\mathcal{T}$ of domain $\omega$}
%\label{mesh_omega}
%\end{subfigure}
\caption{Discretization of the domain \revv{$\hat{\omega}$}.}
\label{mesh_details}
\end{figure}

\begin{figure}
\begin{subfigure}[t]{0.5\textwidth}
    \centering
    \includegraphics[height=2.1in]{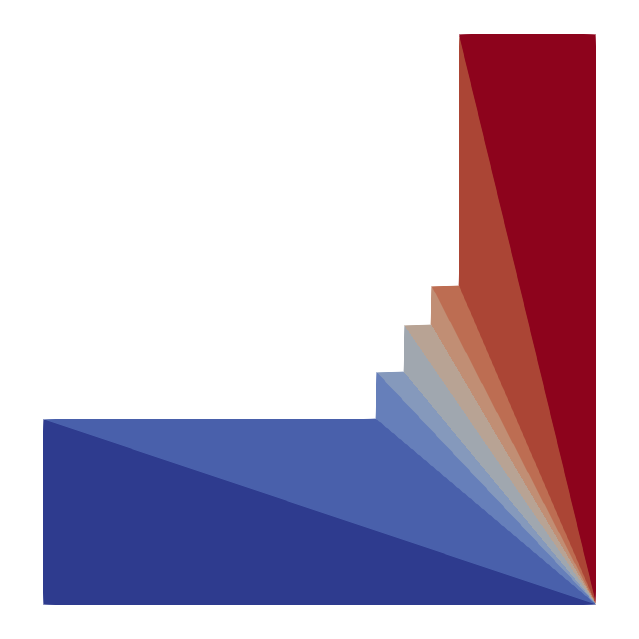}
    \caption{Subdomains decomposition.}
    \label{mesh_2_subdomains}
\end{subfigure}
%\hfill
%\begin{subfigure}[t]{0.32\textwidth}
%    \centering
%    \includegraphics[height=2in]{images/mesh_degenerate.png}
%    \caption{Mesh.}
%    \label{mesh_degenerate}
%\end{subfigure}
%\hfill
\begin{subfigure}[t]{0.4\textwidth}
    \centering
    \includegraphics[height=2in]{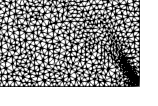}
    \caption{Close up of the mesh at bottom right.}
    \label{mesh_degenerate_bottom_right}
\end{subfigure}
\caption{\rev{Improper domain decomposition: poor mesh quality at bottom right.}}
\label{mesh_condition_number_degenerate}
\end{figure}

\begin{figure}
\begin{subfigure}[t]{0.5\textwidth}
    \centering
    \includegraphics[height=2.1in]{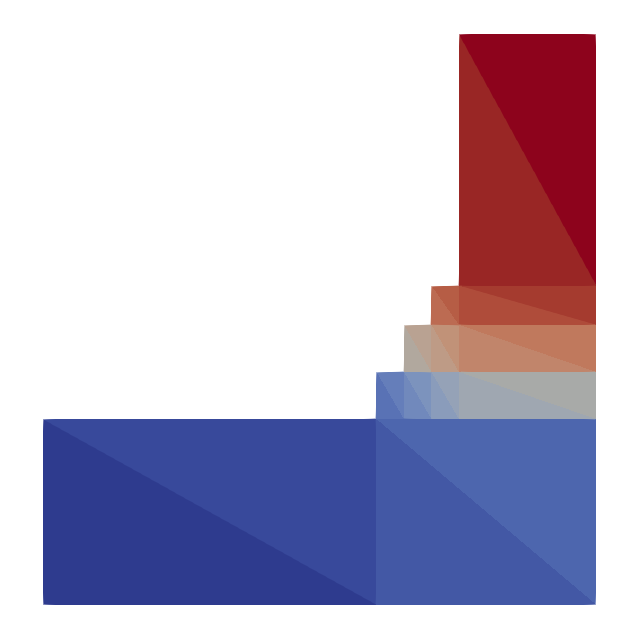}
    \caption{Subdomains decomposition.}
    \label{mesh_1_subdomains}
\end{subfigure}
%\hfill
%\begin{subfigure}[t]{0.32\textwidth}
%    \centering
%    \includegraphics[height=2in]{images/mesh_deformed_2.png}
%    \caption{Mesh on deformed domain.}
%    \label{mesh_deformed_2}
%\end{subfigure}
%\hfill
\begin{subfigure}[t]{0.4\textwidth}
    \centering
    \includegraphics[height=2in]{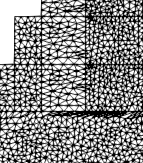}
    \caption{Close up of the region affected by poor mesh quality.}
    \label{mesh_degenerate_poor_quality}
\end{subfigure}
\caption{\rev{Improper domain decomposition: poor mesh quality under variation of the diameter $D_4$.}}
\label{mesh_deformed}
\end{figure}

%As we can observe, the geometry of $\omega$ can be described by $12$ vertices.
%Since one of the main objectives of this work is to develop MOR able to work in a geometrical parametrization setting, before proceeding to the meshing, we consider a decomposition of the domain into $n_{su}$ triangular subdomains as (see, e.g., \cite{geometric_para_1, geometric_para_2}),
%\begin{equation}
%\omega = \bigcup\limits_{i=1}^{n_{su}} \omega_i \ , \ \omega_i \cap %\omega_j = \emptyset \ \text{for} \ i \neq j \ , \ 1 \leq i,j \leq %n_{su} \ .
%\end{equation}
\rev{We set  $n_{su} = 30$ (see Figure \ref{subdomains_ch4}).} The considered mesh is compliant with the triangular subdomains and contains $8887$ triangular elements and $4608$ vertices (see Figure \ref{mesh_elements}). \rev{%The number of subdomains is determined in such a way that under the variation of the parameters hanging nodes are not created and regions with bad mesh quality are not developed. %The domain decomposition shown in figure \ref{subdomains_ch4} is not unique. 
A proper domain decomposition ensures that a suitable mesh can be used for the discretization on each of the subdomains and also throughout one subdomain and the other. As shown in Figure \ref{mesh_condition_number_degenerate}, an improper domain decomposition could generate regions with poor mesh quality. Another issue that might arise from an improper domain decomposition is that, under the variation of the geometrical parameters, poor mesh quality might occur in some regions: an example is showed in Figure \ref{mesh_deformed}.}%, when the diameter $D_4$ is varied, some regions with poor mesh quality might emerge.}

% \begin{figure}
% \begin{subfigure}{0.45\textwidth}
%     \centering
%     \includegraphics[height=2in]{images/mesh_2_subdomains.png}
%     \caption{Subdomains decomposition.}
%     \label{mesh_2_subdomains}
% \end{subfigure}
% \begin{subfigure}{0.45\textwidth}
%     \centering
%     \includegraphics[height=2in]{images/mesh_degenerate.png}
%     \caption{Mesh.}
%     \label{mesh_degenerate}
% \end{subfigure}
% \caption{\rev{Domain decomposition: Poor mesh quality at bottom right.}}
% \label{mesh_condition_number_degenerate}
% \end{figure}

% \begin{figure}
% \begin{subfigure}{0.45\textwidth}
%     \centering
%     \includegraphics[height=2in]{images/mesh_1_subdomains.png}
%     \caption{Subdomains decomposition.}
%     \label{mesh_1_subdomains}
% \end{subfigure}
% \begin{subfigure}{0.45\textwidth}
%     \centering
%     \includegraphics[height=2in]{images/mesh_deformed_2.png}
%     \caption{Mesh on deformed domain.}
%     \label{mesh_deformed_2}
% \end{subfigure}
% \caption{\rev{Domain decomposition: Poor mesh quality under variation of geometric parameter.}}
% \label{mesh_deformed}
% \end{figure} 

The minimum and maximum mesh size, that is measured as distance between vertices of an element of the mesh, are $0.047$m and $0.16$m respectively. %We show the mesh condition number in Figure \ref{mesh_quality} ranging from $1$ to $\infty$, with $1$ being a perfectly shaped element.
The quality of each mesh element, $q_e$, could be estimated by using the following formula \cite{sara_paper}:
\begin{equation}\label{mesh_quality_formula}
q_e = \frac{4 \sqrt{3} A}{l_1^2+l_2^2+l_3^2} \ ,
\end{equation}
where $A$ is the area of the element, and $l_1$, $l_2$ and $l_3$ are the lengths of its three edges. The minimum value of $q_e$ was $0.25$, that is sufficiently far from zero. 
%We use a mesh of $\hat{\omega}$ containing $8887$ triangular elements and $4608$ vertices. %\nirav{The mesh size is measured as distance between vertices of an element of the mesh}. 
%The minimum mesh size is $0.047$ m and the maximum one is $0.16$ m. Its minimum quality is $q_e = 0.25$ (eq. \ref{mesh_quality_formula}). 

Notice that we use a coarser mesh with respect to the one used for the FOM benchmark tests in \revv{\ref{appendix}}. Such a choice is justified by the fact that the FOM solution is required to be solved at many parameters values, so using a fine mesh can be very costly and make prohibitive the collection of the high-fidelity database. %We highlight that this does not affect the We also would like to clarify that the time spent for the collection of the FOM snapshots does
%not affect the efficiency of the ROM, because snapshot collection is included in the offline
%stage only.%Besides, as explained earlier, the aim of the reduced basis method is to approximate the full order is model solution rather than approximating the actual solution (see section \ref{rb_section}, equations \ref{rb_error_bound_energy} and \ref{rb_error_bound_momentum}).
%While, the mesh subdomains are same as the section \ref{benchmark_test_section}, the condition number and the mesh are shown in figure \ref{rb_mesh_details}. The mesh contains $8887$ triangular elements and $4608$ vertices.The minimum stance between any two vertices of a cell was $0.047$ and the maximum distance between any two vertices of a cell was $0.16$. The minimum quality among all elements was $0.25$ and the average quality over all elements was $0.94$ (equation \ref{mesh_quality_formula}).
%\begin{figure}
%\centering
%\begin{subfigure}{0.3\textwidth}
%\centering
%\includegraphics[height=1.5in]{images/rb_mesh.png}
%\caption{Mesh $\mathcal{T}$ of domain $\omega$}
%\label{mesh_omega}
%\end{subfigure}
%\begin{subfigure}{0.3\textwidth}
%\centering
%\includegraphics[height=1.5in]{images/condition_number_rb_mesh.png}
%\caption{Mesh quality}
%\label{mesh_quality}
%\end{subfigure}
%\caption{Discretization of domain $\omega$}
%\label{rb_mesh_details}
%\end{figure}

The ranges of physical and geometrical parameters for training and testing are reported in Table \ref{range_of_parameters}. \revv{The sampling is carried out by using a Latin Hypercube Sampling (LHS) approach \cite{LHS_reference} which is a statistical method for generating near-random samples of parameter values from a multidimensional distribution. LHS divides the parameter space into equal partitions and samples parameters from each partition. In this manner, it is ensured that the patterns from entire parameter space are represented. The process has been repeated multiple times in order to ensure that random nature of samplings do not affect the final result.} ANN has been trained by using the 70\% of the total data provided by the full order model whilst the remaining 30\% is used for the validation.

\begin{table}[H]
\centering
%\resizebox{0.95\textwidth}{!}{
\begin{tabular}{|c|c|c|}
\hline
Parameter & Minimum value & Maximum value \\\hline
$t_0$ & 2.3 & 2.4\\\hline
$t_1$ & 0.5 & 0.7\\\hline
$t_2$ & 0.5 & 0.7\\\hline
$t_3$ & 0.4 & 0.6\\\hline
$t_4$ & 3.05 & 3.35\\\hline
$D_0$ & 13.5 & 14.5\\\hline
$D_1$ & 8.3 & 8.7\\\hline
$D_2$ & 8.8 & 9.2\\\hline
$D_3$ & 9.8 & 10.2\\\hline
$D_4$ & 10.4 & 10.8\\\hline
$k$ & 9.8 & 10.2\\\hline
$\mu$ & 1.9e9 & 2.5e9 \\\hline
$\lambda$ & 1.2e9 & 1.8e9 \\\hline
$\alpha$ & 0.8e-6 & 1.2e-6\\\hline
% x.x & DD/MM/YYYY  &   &  \\\hline
\end{tabular}
%}
\caption{Parameters ranges used for MOR training and testing.}
\label{range_of_parameters}
\end{table}

%\michele{DA TENERE IN CONTO DIRETTAMENTE NEL CAPITOLO DEI RISULTATI}
The accuracy of our MOR approach is quantified by the relative error defined as follows
\begin{equation}\label{rel_error}
\epsilon_{rel,X_h} = \frac{||X_h-X_h^{rb}||}{||X_h||},
\end{equation}
where $X_h$ and $X_h^{rb}$ are the finite element solution and the corresponding reduced basis solution, respectively.
%\michele{We also should define the projection error that we use as benchmark for MOR accuracy.}
%\revv{If $X_h^{\pi}$ is the projection of finite element solution $X_h$ on the reduced basis space (equations \eqref{temperature_solution_projection},\eqref{displacement_solution_projection}),} 
We consider the projection error \revv{between the finite element solution $X_h$ and its projection on the reduced basis space $X_h^{\pi}$}, %If $X_h$ is the finite element solution and $X_h^{\pi}$ is its projection on the reduced basis space, the projection error is defined as follows}: %, in relevant norm $||\cdot||$, as defined above, between finite element solution $X_h$ and its projection onto reduced basis space as,
\begin{equation}\label{projection_error}
\epsilon_{proj,X_h} = \frac{||X_h- X_h^{\pi}||}{||X_h||} \ ,
\end{equation}
as benchmark for the relative error. $||\cdot||$ is the relevant norm ($||\cdot||_{H^1_{r,h}(\omega)}$ and $||\cdot||_{\mathbb{U}_h}$).
%$\bm{B}^T$ is the transpose of the projection operators $\bm{B}_T$ for the temperature field and $\bm{B}_u$ for the displacement field.

%On the other hand, the computational efficiency is referred to the speed-up $SU$, defined as the ratio between the time taken to solve the full order model and the time taken to solve the reduced algebraic system of equations, at a given parameter tuple.
%\begin{equation}\label{SU}
%SU = \dfrac{t_{FOM}}{t_{on}},
%\end{equation}
%where $t_{FOM}$ is the time taken to solve the full order model and $t_{on}$ is the time taken to solve the reduced algebraic system of equations, at a given parameter tuple.
%\end{itemize}

\subsubsection{Thermal model}

\begin{figure}[H]
\begin{center}
\includegraphics[width=0.6\textwidth]{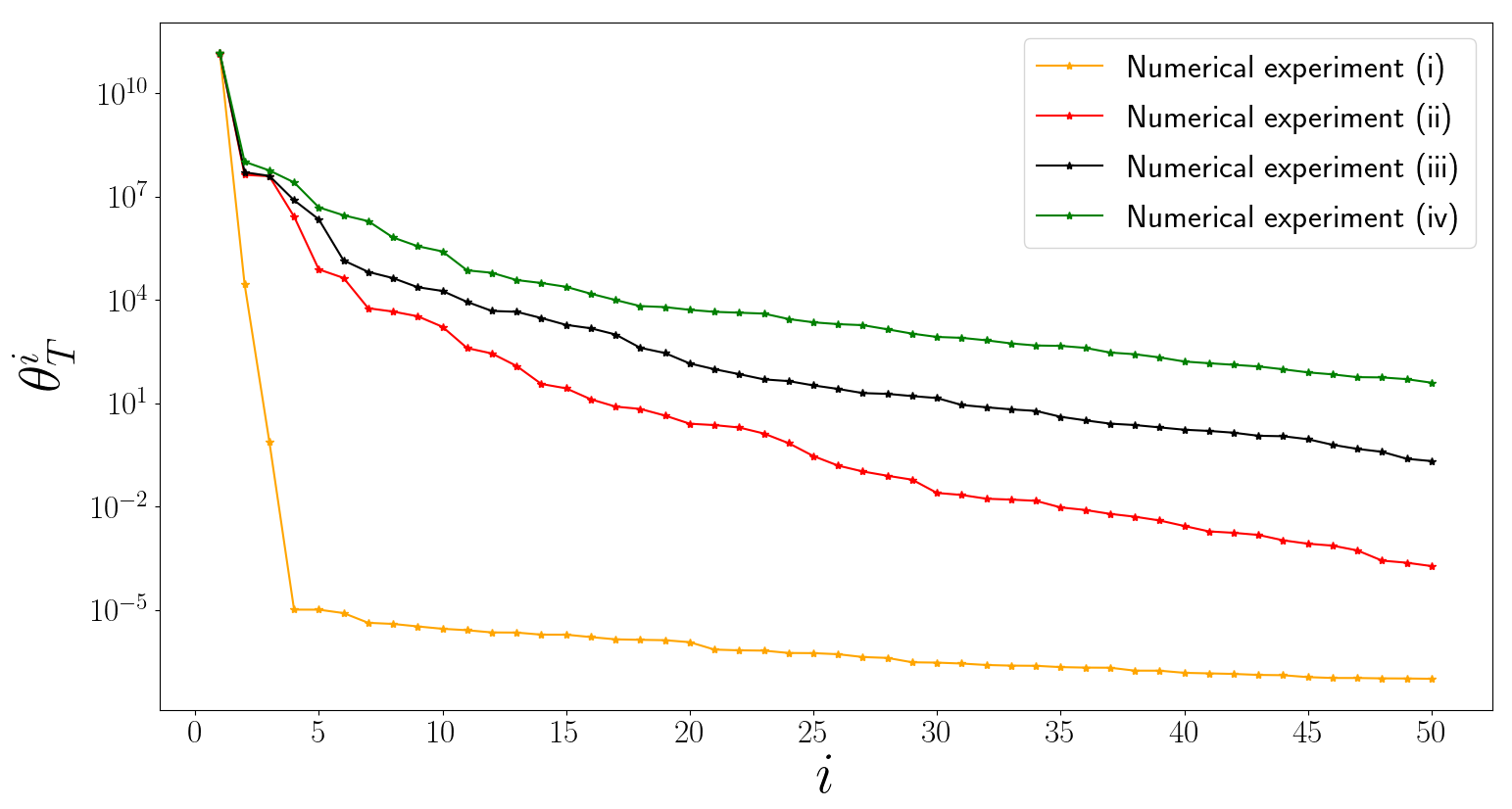}
\end{center}
\caption{Thermal model: \revv{plot of the} eigenvalues \revv{$\lbrace \theta_T^i \rbrace_{i=1}^{50}$ sorted in descending order} for all the \revvv{numerical} experiments considered.}
\label{eigenvalue_thermal_active}
\end{figure}

%In this section, we compare the performance of our MOR approaches for the thermal model. %The minimum admissible relative eigenvalue for reduced basis space was kept at $1e-4$.
We consider four \revvv{numerical} experiments that differ in terms of kind (physical and/or geometrical) and number of the parameters considered: % with respect to which the training phase is performed: %with parameter tuple containing different number of physical and geometric parameters.
\begin{itemize}
%\item \emph{experiment (i)}: 4 physical parameters $\Xi = \lbrace k, \mu, \lambda, \alpha \rbrace$.
\item \emph{\revvv{Numerical} experiment (i)}: 1 physical parameter: $\Xi = \lbrace k \rbrace$.
\item \emph{\revvv{Numerical} experiment (ii)}: 1 physical parameter and 3 geometric parameters: $\Xi = \lbrace k, t_0, D_2, D_4 \rbrace$.
\item \emph{\revvv{Numerical} experiment (iii)}: 1 physical parameter and 6 geometric parameters: $\Xi = \lbrace k, t_0, t_2, t_4, D_0, D_2, D_4 \rbrace$.
\item \emph{\revvv{Numerical} experiment (iv)}: 1 physical parameter and all (10) geometric parameters: $\Xi = \lbrace k, t_0, t_1, t_2, t_3, t_4, D_0, D_1, D_2, D_3, D_4 \rbrace$.
\end{itemize}
\revvv{Table \ref{thermal_training_validation_testing_table} shows the number of samples provided by the full order model, $n_{tr}$, as well as the number of samples used for training and testing of ANN}. Regarding the computation of POD space, for \emph{\revvv{numerical} experiment (i)}, $50$ FOM snapshots were considered while for the other ones $1000$. %POD space was constructed from these snapshots and applying criteria as per equation \eqref{eigenvalue_pod_criteria}.}
The eigenvalues decay is shown in Figure \ref{eigenvalue_thermal_active}. We see that the decay related to the \revvv{numerical} experiment (iv) is the slowest. %Therefore, a higher number of basis functions needs to be considered in order to obtain an accurate reconstruction of the field.
This is due to the fact that in the \revvv{numerical} experiment (iv) we consider a larger number of parameters, so the system exhibits a greater complexity, and the modal content is more wide. %is more complex.

\begin{figure}[H]
        \centering
        \begin{subfigure}[b]{0.475\textwidth}
            \centering
            \includegraphics[width=\textwidth]{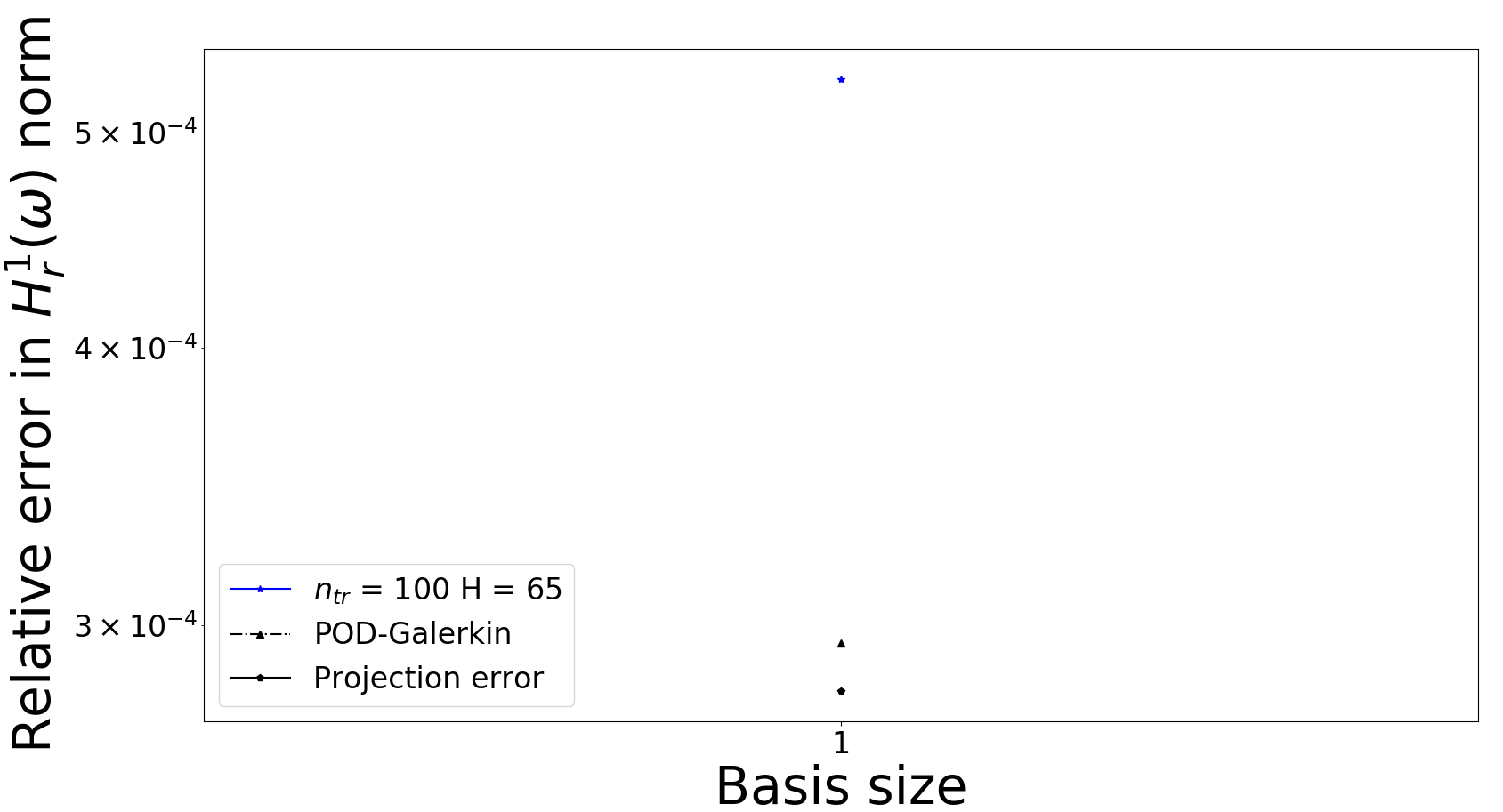}
            \caption[]%
            {{\small \revvv{Numerical} experiment (i).}}
            \label{thermal_model_1_material_error}
        \end{subfigure}
        \hfill
        \begin{subfigure}[b]{0.475\textwidth}
            \centering
            \includegraphics[width=\textwidth]{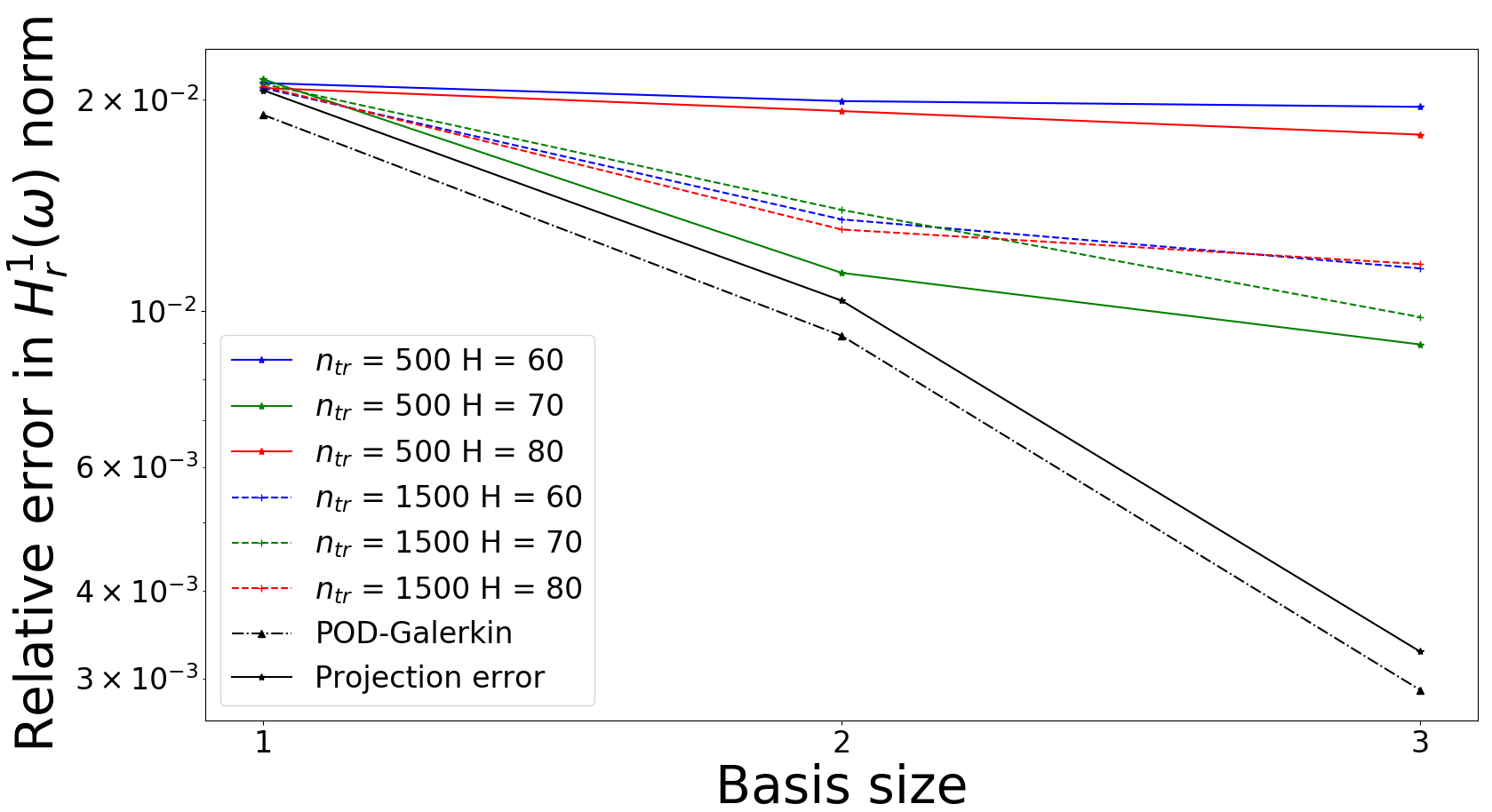}
            \caption[]%
            {{\small \revvv{Numerical} experiment (ii).}}
            \label{thermal_model_1_material_3_geometric_error}
        \end{subfigure}
        \vskip\baselineskip
        \begin{subfigure}[b]{0.475\textwidth}
            \centering
            \includegraphics[width=\textwidth]{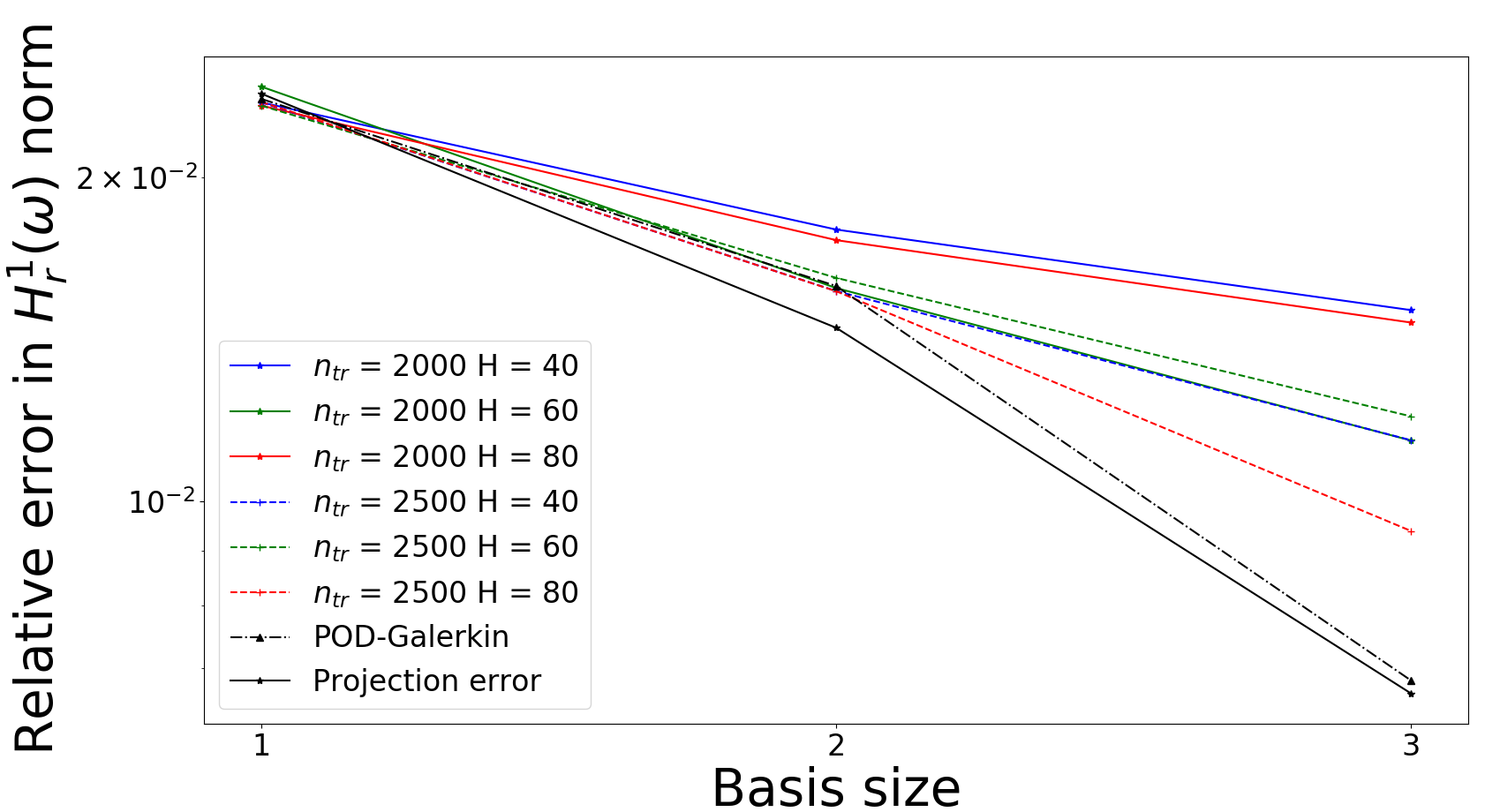}
            \caption[]%
            {{\small \revvv{Numerical} experiment (iii).}}
            \label{thermal_model_1_material_6_geometric_error}
        \end{subfigure}
        \hfill
        \begin{subfigure}[b]{0.475\textwidth}
            \centering
            \includegraphics[width=\textwidth]{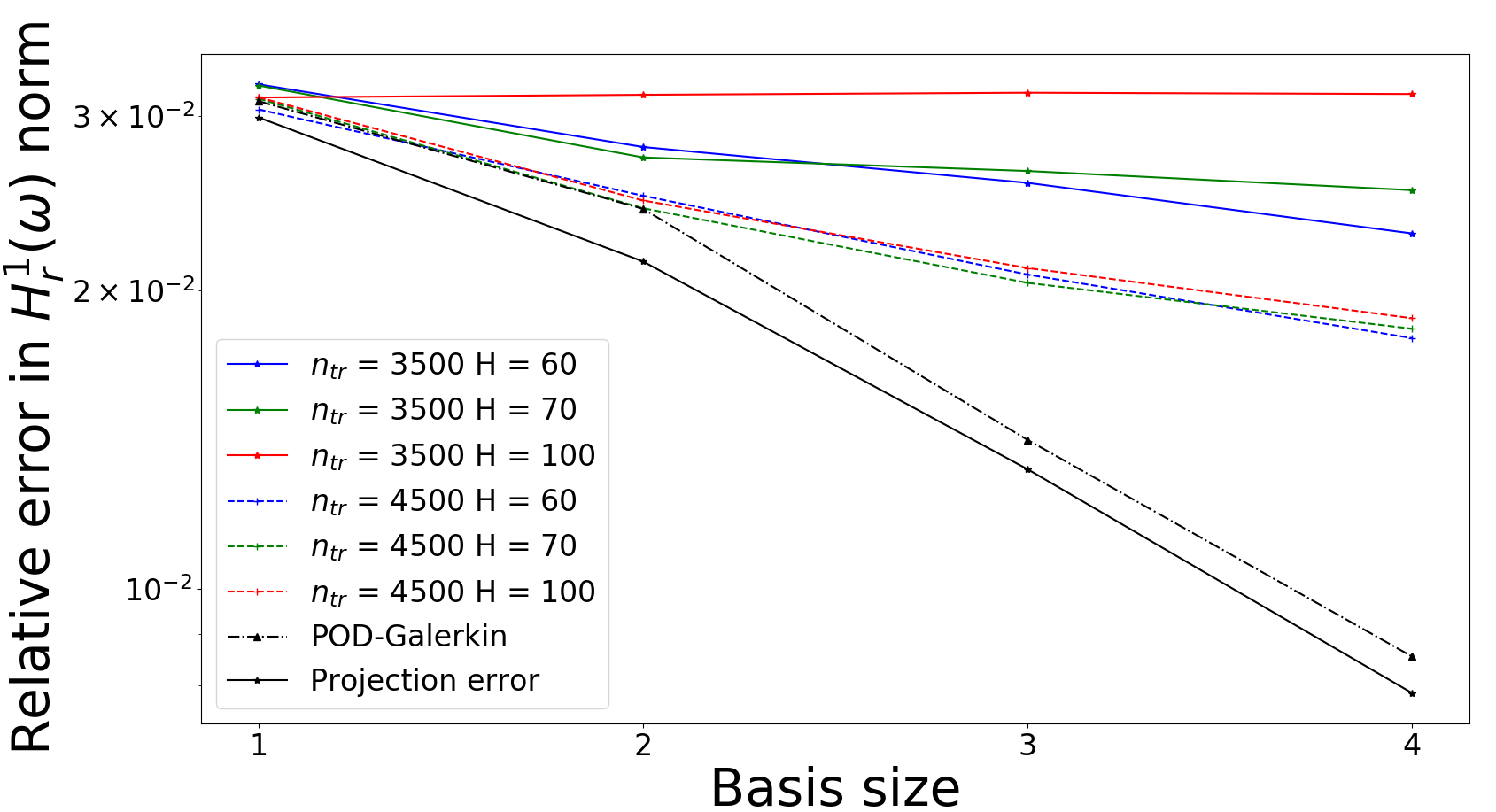}
            \caption[]%
            {{\small \revvv{Numerical} experiment (iv).}}
            \label{thermal_model_all_parameters_error}
        \end{subfigure}
        \caption[]%
        {Thermal model: error analysis for POD-G and POD-ANN for all the \revvv{numerical} experiments carried out. %\michele{I would increase a bit the basis size.} \nirav{The basis size was considered to achieve accuracy of the order of ~1e-2. Beyond this the training error dominates and not the projection error. The ~1e-2 error also matches with POD minimum admissible eigenvalue of $1e-4$. I am not sure why do we need to increase the number of basis function. }
}
        \label{thermal_error_analysis_plots}
    \end{figure}

Figure \ref{thermal_error_analysis_plots} shows the relative error (\ref{rel_error}) both for POD-ANN, related to different values $n_{tr}$ and depth of hidden layers $H$, %\michele{(This variable in the previous chapter has been called $d_l$)}\nirav{See new sentence added below Sigmoid activateion function.}
and POD-G. We also report the projection error \eqref{projection_error}. We observe that the performance of the POD-ANN method crucially depends on the values of $n_{tr}$ and $H$. As expected, if we expand the training set and increase the depth of hidden layers, we obtain more accurate predictions when the number of parameters considered starts to get significative (\revvv{numerical} experiments (iii) and (iv)). %This result is in agreement with that reported for in  \cite{pod_nn_paper} for a nonlinear modelling framework. 
\revvv{Unlike \cite{pod_nn_paper}, we observe that the POD-G method results to be in general more accurate than the POD-ANN method. This could be justified by considering that in the nonlinear framework, investigated in \cite{pod_nn_paper},  the affine expansion could not be enforced and an Empirical Interpolation Method (EIM) \cite{eim_reference} is used within the POD-G approach. Its implementation introduces interpolation error during the assembling of the reduced equations system by significantly affecting the accuracy of the POD-G method.} % TODO poi dobbiamo parlare dei tempi computationali Moreover, %S

\begin{figure}[H]
        \centering
        \begin{subfigure}[b]{0.3\textwidth}
            \centering
            \includegraphics[width=\textwidth]{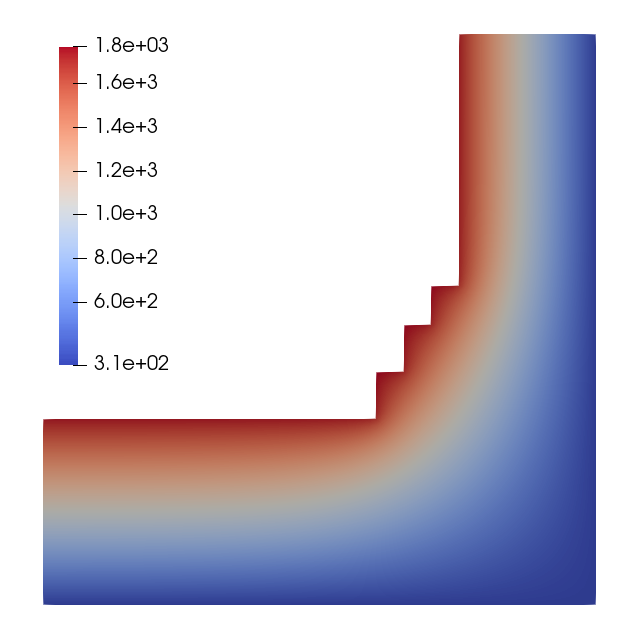}
            \caption[]%
            {{\small FOM solution}}
            \label{fem_solution_ref}
        \end{subfigure}
        \hfill
        %\begin{subfigure}[b]{0.35\textwidth}
         %   \centering
         %   \includegraphics[width=\textwidth]{images/projected_solution_ref.png}
         %   \caption[]%
         %   {{\small Projected solution}}
         %   \label{projected_solution_ref}
        %\end{subfigure}
        %\vskip\baselineskip
        \begin{subfigure}[b]{0.3\textwidth}
            \centering
            \includegraphics[width=\textwidth]{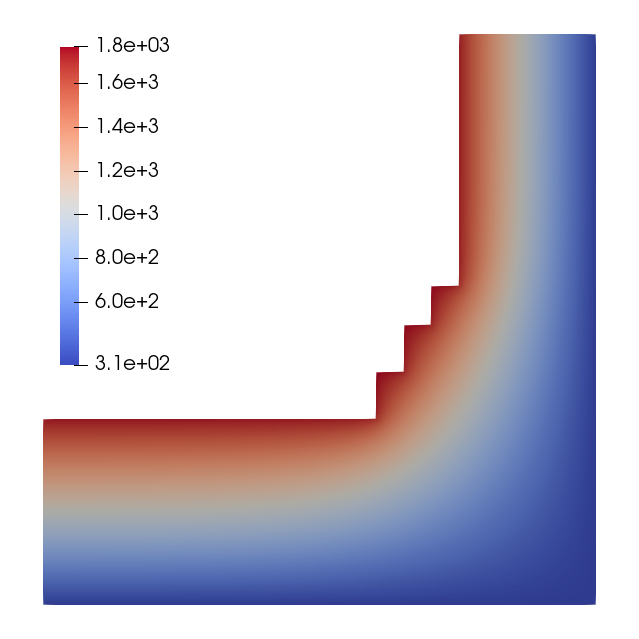}
            \caption[]%
            {{\small POD-G solution}}
            \label{podg_solution_ref}
        \end{subfigure}
        \hfill
        \begin{subfigure}[b]{0.3\textwidth}
            \centering
            \includegraphics[width=\textwidth]{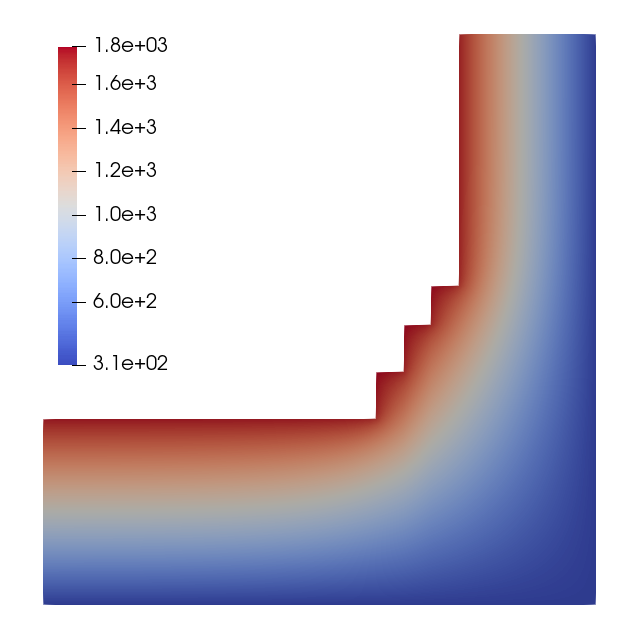}
            \caption[]%
            {{\small POD-ANN solution}}
            \label{pod_ann_solution_ref}
        \end{subfigure}
        \caption[]%
        {Thermal model: comparison between the temperature field (in K) computed by the FOM and by the POD-G and POD-ANN methods related to the \revvv{numerical} experiment (iv) for $\Xi = \lbrace 2.365, 0.6, 0.6, 0.5, 3.2, 14.10, 8.50, 9.2, 9.9, 10.6, 10 \rbrace$. We consider 4 POD modes. For POD-ANN, we set $n_{tr} = 4500$ and $H = 70$.}
        \label{solutions_ref}
    \end{figure}

Illustrative representations of the computed FOM and MOR are displayed in Figure \ref{solutions_ref} related to the \revvv{numerical} experiment (iv) for the parameters tuple %\nirav{based on \emph{experiment (iv)}}: %and \ref{solutions_para} %We now present the full order model solution, the full order solution projected on reduced basis space, POD-Galerkin solution and POD-ANN solution (figures \ref{solutions_ref},\ref{solutions_para}) for the problem presented in the section \ref{benchmark_actual} for the parameter values :
\begin{gather*}
\Xi = \lbrace 2.365, 0.6, 0.6, 0.5, 3.2, 14.10, 8.50, 9.2, 9.9, 10.6, 10 \rbrace \ . %\\
%\Xi_2 = \lbrace 2.165, 0.6, 0.6, 0.45, 3.2, 14.10, 8.30, 9.2, 9.9, 10.6, 10. \rbrace \ .
\end{gather*}
We use 4 POD basis. The POD-ANN solution was computed with $n_{tr}=4500$ and $H=70$. % \michele{inserire un commento in piu}. %The reduced basis space contained $4$ basis functions \michele{inserire un commento in piu}.
As we can see from Figure \ref{solutions_ref}, both MOR approaches are able to provide a good reconstruction of the temperature field.

We conclude by proving some information about the efficiency of our MOR approach. \rev{We report in Table \ref{thermal_offline_time} some estimations related to the offline time for all the numerical experiments carried out. We observe that the time taken by POD for numerical experiments (ii)-(iv) is much larger than the numerical experiment (i). This is fully justified by the fact that for the numerical experiments (ii)-(iv) we consider a larger number of snapshots (1000 instead of 50 as discussed above) for the computation of the reduced space. On the other hand, it should also be noted that the ANN training is faster for the numerical experiment (i) where only physical parameters are involved. This could be attributed to the fact that the introduction of geometric parameters increases the complexity of the input-output map that ANN is expected to learn. If on one hand the offline cost of the POD-G method is most composed of time taken by the computation of the snapshots from which the reduced space is extracted and the time taken by the computation of the POD modes, on the other hand the one related to the POD-ANN method is mainly associated to the computation of training data. So, when the parameter space is large (as for the numerical experiments (ii)-(iv)), the total offline cost of the POD-ANN method could be importantly greater than the one related to the POD-G method.} % whilst the one related to the POD-G method is almost. 
%For the most complex case represented by the numerical experiment (iv), we obtain about s for POD-G and about for POD-ANN. totally importantly greater than the one associated to the POD-G method because of the computation of training data by FOM.}
%as well as the one taken by the ANN training are much larger for numerical experiments (ii)-(iv) than the numerical experiment (i). This is expected because for the numerical experiments (ii)-(iv) we consider both a larger number of snapshots for the computation of the reduced space and a larger value of $n_{tr}$. %because for these numerical experiments we consider a larger number of snapshots (1000 instead of 50). %when geometrical parameters are considered.} %TThe time taken by POD is much less than the time for computation of the snapshots. Table \ref{thermal_offline_time} also reports time taken for the training of ANN after computation of the training data. While, the time taken for the training of ANN is higher than time taken for the eigenvalue decomposition, 
%
% We observe that the most dominant factor which  between the POD-G and POD-ANN methods in terms of offline cost is due to the computation of training data by FOM is the most dominant factor which makes   for the training cost of ANN.} %It should also be noted that the training of ANN is very fast when only material parameters are involved. This can be attributed to the fact that introduction of geometric parameters increases complexity of the analytical map that ANN is expected to learn.} 
We report in Table \ref{thermal_rb_time_plots} the online time related to the POD-G and POD-ANN methods for all the \revvv{numerical} experiments carried out. % I would like to avoid offline time and speedup since the comparison is architecture, code and number of degrees of freedom dependent. I would like to avoid it. Comparison of two MOR methods is less architecture and problem dependent. Kindly give it a serious thought.} \michele{I manily agree with you, but we should give some information about this topic. I mean, we provide the online time but this parameter, alone, does give no information about the efficiency of our ROM approach. I agree with you that we can avoid to provide a lot of details related to the computational cost of the offline phase (time taken by computation of POD modes, etc...). So, you are saying me that for solving the FOM for a parameter tupla you spend 0.08s? For which model? I would have an estimate for all the models under investigation, thanks! (0.08 s is a good value, the online time are significantly lower and we obtain SU of 10/100 as order of magnitude that are ok. Thanks}% as well as the time spent by a FOM simulation, $t_{FOM}$. %\ref{thermal_rb_time_plots}.
\revvv{As can be seen, the online time of POD-G method increases significantly in presence of geometric parameters by moving from $7e-4$ s (\revvv{numerical} experiment (i)) to $1.3$/$1.5e-2$ s (\revvv{numerical} experiments (ii)-(iv)). On the other hand, the time taken by POD-ANN online stage remains relatively constant for all the \revvv{numerical} experiments under investigation, around $5e-4$. So the computational efficiency of POD-ANN is much higher, of almost two order of magnitude, than POD-G when geometrical parametrization is considered.} %In terms of speed-up we obtain for POD-G and for POD-ANN in the worst cases.

\begin{table}[H]
\centering
\begin{tabular}{|c|c|c|c|}
\hline
& $n_{tr}$ & Training & Testing \\ %& Testing cases \\
\hline
\multirow{1}{0.23\textwidth}{Numerical experiment (i)} & 100 & 70 & 30 \\%& 50 \\
\hline
\multirow{2}{0.23\textwidth}{Numerical experiment (ii)} & 500 & 350 & 150 \\%& 100 \\
\cline{2-4}
& 1500 & 1050 & 450 \\% & 100 \\
\hline
\multirow{2}{0.23\textwidth}{Numerical experiment (iii)} & 2000 &  1400 & 600 \\% & 100 \\
\cline{2-4}
& 2500 & 1750 & 750 \\% & 100 \\
\hline
\multirow{2}{0.23\textwidth}{Numerical experiment (iv)} & 3500 & 2450 & 1050 \\ %& 100 \\
\cline{2-4}
& 4500 & 3150 & 1350 \\ % & 100 \\
\hline
\end{tabular}
\caption{\revvv{Thermal model: number of total samples $n_{tr}$ by FOM and number of samples used for training and testing of ANN.}}
\label{thermal_training_validation_testing_table}
\end{table}

\begin{table}[H]
\centering
\begin{tabular}{|c|c|c|c|c|c|c|}
\hline
& $t_{off}^{POD-G}$  & $t_{off}^{POD-ANN}$ & $t_{POD}$ & $t_{tr}$ & $t_{FOM}$ & $t_{proj}$\\
\hline
\multirow{1}{16em}{Numerical experiment (i)} & $\approx$ 5 & $\approx$ 1.2e1 & \multirow{1}{3em}{8.4e-1} & 2.6e-1 & \multirow{4}{2em}{8e-2} & 6.0e-4 \\
\cline{1-5} \cline{7-7}
\multirow{1}{16em}{Numerical experiment (ii)} & $\approx$ 9e1 & $\approx$ 1.3e2 & \multirow{1}{3em}{1.2e1} & 5.9e-1 &  & 7.5e-4 \\
\cline{1-5} \cline{7-7}
\multirow{1}{16em}{Numerical experiment (iii)} & $\approx$ 9e1 & $\approx$ 3.7e2 & \multirow{1}{3em}{1.2e1} & 7.4e1 & & 7.1e-4\\
\cline{1-5} \cline{7-7}
\multirow{1}{16em}{Numerical experiment (iv)} & $\approx$ 9e1 & $\approx$ 5.0e2 & \multirow{1}{3em}{1.2e1} & 4.5e1 & & 7.3e-4\\
\hline
\end{tabular}
\caption{\rev{Thermal model: time (in s) taken by (i) the entire offline stage ($t_{off}$), (ii) the computation of the POD modes ($t_{POD}$), (iii) the training of ANN ($t_{tr}$), (iv) the computation of a FOM solution ($t_{FOM}$) and (v) the projection of a FOM solution on the POD space ($t_{proj}$). Concerning POD-ANN, we use $n_{tr} = 100, H = 65$ for the \revvv{numerical} experiment (i), $n_{tr} = 500, H=70$ for the \revvv{numerical} experiment ii), $n_{tr} = 2500, H=80$ for the \revvv{numerical} experiment (iii) and  $n_{tr} = 4500, H = 70$ for the \revvv{numerical} experiment (iv).}}
\label{thermal_offline_time}
\end{table}

% \begin{table}
% \centering
% \begin{tabular}{|c|c|}
% \hline
% Problem & $(WT)_h$ \\
% \hline
% Offline time & 8e-2\\
% \hline
% \end{tabular}
% \caption{\revv{Thermal model: offline time (in s).}}
% \label{thermal_offline_time}
% \end{table}

\begin{table}[H]
\centering
\begin{tabular}{|c|c|c|c|}
\hline
& Basis size & POD-G & POD-ANN \\
\hline
\multirow{1}{16em}{\revvv{Numerical} experiment (i)} & 1 & 7.0e-4 & 4.9e-4\\
\hline
\multirow{1}{16em}{\revvv{Numerical} experiment (ii)} %& 1 & 1.3e-2 & 4.9e-4\\
%& 2 & 1.3e-2 & 4.8e-4\\
& 3 & 1.3e-2 & 4.8e-4\\
\hline
\multirow{1}{16em}{\revvv{Numerical} experiment (iii)} %& 1 & 1.5e-2 & 4.9e-4\\
%& 2 & 1.4e-2 & 5e-4\\
& 3 & 1.5e-2 & 4.9e-4\\
\hline
\multirow{1}{16em}{\revvv{Numerical} experiment (iv)} %& 1 & 1.3e-2 & 5.1e-4\\
%& 2 & 1.3e-2 & 5e-4\\
%& 3 & 1.3e-2 & 5e-4\\
& 4 & 1.3e-2 & 5.1e-4\\
\hline
\end{tabular}
\caption{Thermal model: online time (in s) for all the \revvv{numerical} experiments under investigation. Concerning POD-ANN, we use $n_{tr} = 100, H = 65$ for the \revvv{numerical} experiment (i), $n_{tr} = 500, H=70$ for the \revvv{numerical} experiment ii), $n_{tr} = 2500, H=80$ for the \revvv{numerical} experiment (iii) and  $n_{tr} = 4500, H = 70$ for the \revvv{numerical} experiment (iv).} %\michele{i would re-organize such results in terms of table.} \nirav{Done.}}
\label{thermal_rb_time_plots}
\end{table}

\begin{comment}
      \begin{figure}[H]
        \centering
        \begin{subfigure}[b]{0.35\textwidth}
            \centering
            \includegraphics[width=\textwidth]{images/fem_solution_para.png}
            \caption[]%
            {{\small FEM solution}}
            \label{fem_solution_para}
        \end{subfigure}
        \hfill
        \begin{subfigure}[b]{0.35\textwidth}
            \centering
            \includegraphics[width=\textwidth]{images/projected_solution_para.png}
            \caption[]%
            {{\small Projected solution}}
            \label{projected_solution_para}
        \end{subfigure}
        \vskip\baselineskip
        \begin{subfigure}[b]{0.35\textwidth}
            \centering
            \includegraphics[width=\textwidth]{images/pod_g_solution_para.png}
            \caption[]%
            {{\small POD-Galerkin solution}}
            \label{podg_solution_para}
        \end{subfigure}
        \hfill
        \begin{subfigure}[b]{0.35\textwidth}
            \centering
            \includegraphics[width=\textwidth]{images/nn_rb_solution_para.png}
            \caption[]%
            {{\small POD-ANN solution}}
            \label{pod_ann_solution_para}
        \end{subfigure}
        \caption[]%
        {Solution computed at parameter tuple $\Xi_2$}
        \label{solutions_para}
    \end{figure}
\end{comment}

\subsubsection{Mechanical model}

%Similar to thermal model, we compare the performance for POD-Galerkin and POD-ANN approach for mechanical and thermo-mechanical model for POD-Galerkin and mechanical model for POD-ANN.
We remark that for POD-ANN we refer to the $(WM)_h$ model, whilst we consider $(WM1)_h$ and $(WM2)_h$ models for POD-G. %\TODO%The minimum admissible relative eigenvalue for reduced basis space was kept at $1e-4$.
As done for the thermal model, we consider four different \revvv{numerical} experiments having different kinds and numbers of parameters: %that differ primarily in terms of kind (physical and/or geometrical) and number of the parameters with respect to which the training phase is performed: %with parameterWe consider four set of experiments with parameter tuple containing different number of geometric parameters.
\begin{itemize}
\item \emph{\revvv{Numerical} experiment (i)}: 4 physical parameters: $\Xi = \lbrace k, \mu, \lambda, \alpha \rbrace$.
\item \emph{\revvv{Numerical} experiment (ii)}: 4 physical parameters and all 3 geometric parameters: $\Xi = \lbrace k, \mu, \lambda, \alpha, t_0, D_2, D_4 \rbrace$.
\item \emph{\revvv{Numerical} experiment (iii)}: 4 physical parameters and 6 geometric parameters: $\Xi = \lbrace k, \mu, \lambda, \alpha, t_0, t_2, t_4, D_0, D_2, D_4 \rbrace$.
\item \emph{\revvv{Numerical} experiment (iv)}: 4 physical parameters and all (10) geometric parameters: \\ $\Xi = \lbrace k, \mu, \lambda, \alpha, t_0, t_1, t_2, t_3, t_4, D_0, D_1, D_2, D_3, D_4 \rbrace$.
\end{itemize}

\begin{figure}[H]
	\centering
	\begin{subfigure}[b]{0.475\textwidth}
		\centering
		\includegraphics[width=\textwidth]{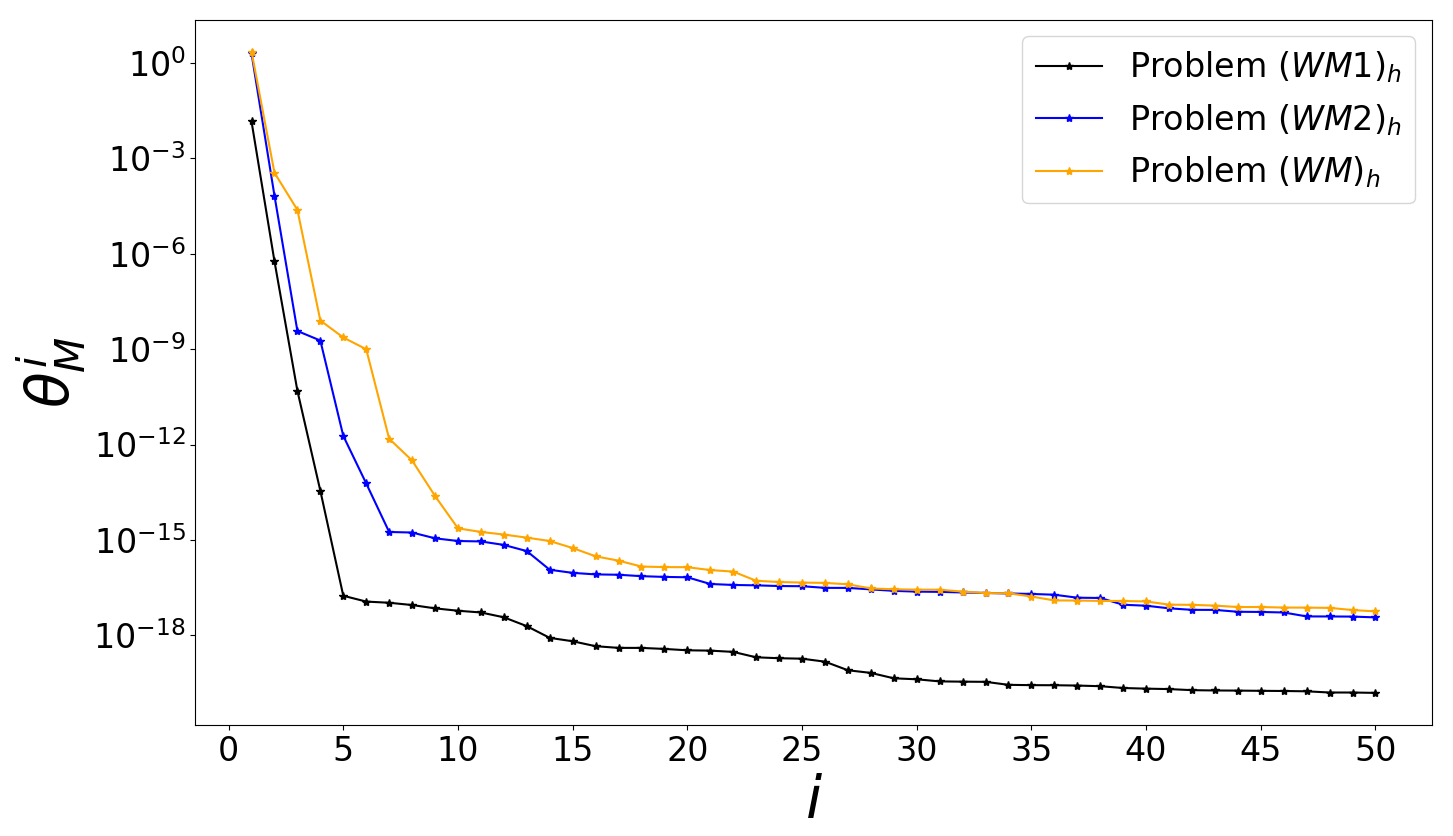}
		\caption[]%
		{{\small \revvv{Numerical} experiment (i).}}
		\label{mechanical_model_4_materials_eigenvalues}
	\end{subfigure}
	\hfill
	\begin{subfigure}[b]{0.475\textwidth}
		\centering
		\includegraphics[width=\textwidth]{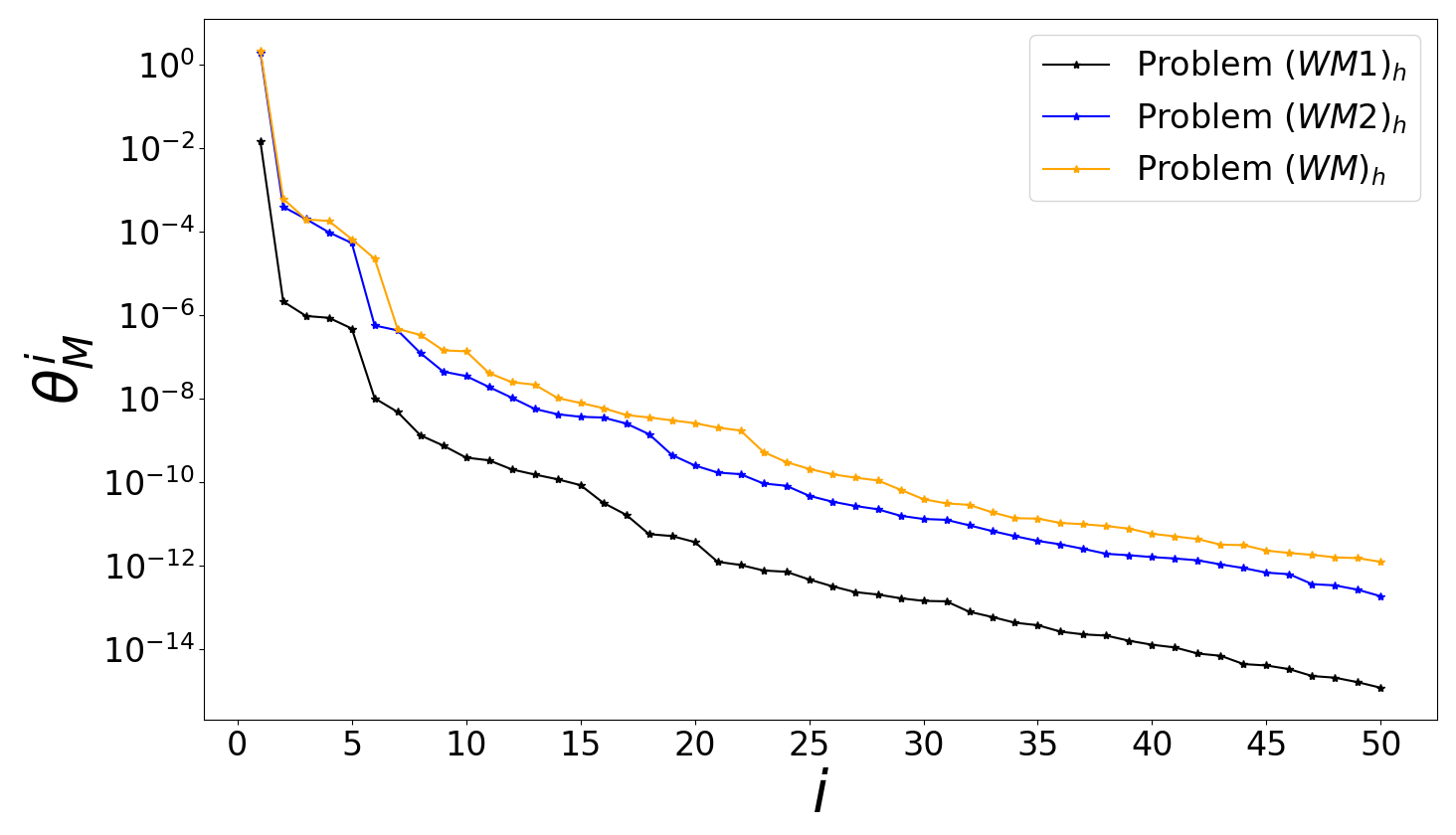}
		\caption[]%
		{{\small \revvv{Numerical} experiment (ii).}}
		\label{mechanical_model_4_materials_3_geometric_eigenvalues}
	\end{subfigure}
	\vskip\baselineskip
	\begin{subfigure}[b]{0.475\textwidth}
		\centering
		\includegraphics[width=\textwidth]{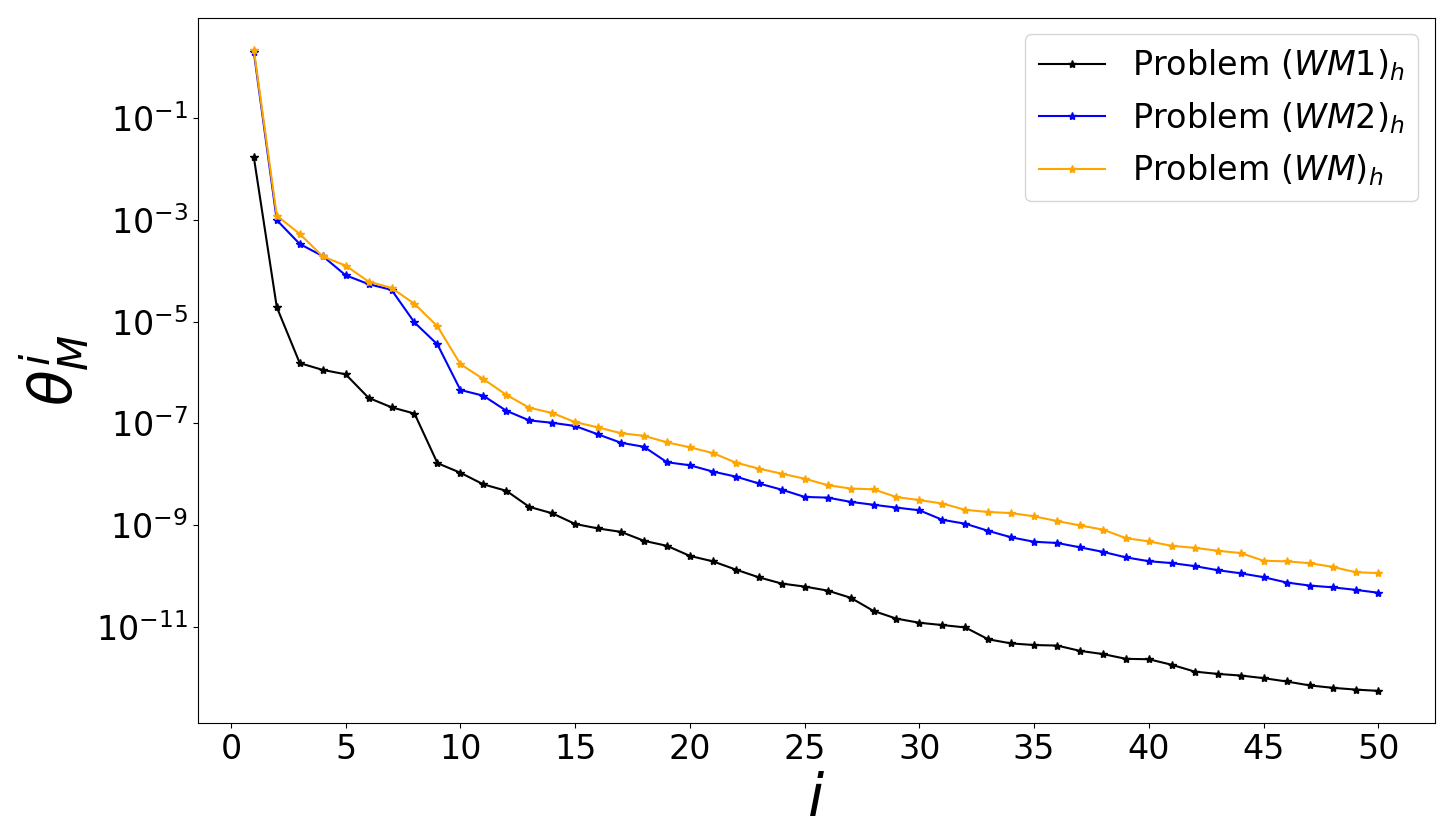}
		\caption[]%
		{{\small \revvv{Numerical} experiment (iii).}}
		\label{mechanical_model_4_materials_6_geometric_eigenvalues}
	\end{subfigure}
	\hfill
	\begin{subfigure}[b]{0.475\textwidth}
		\centering
		\includegraphics[width=\textwidth]{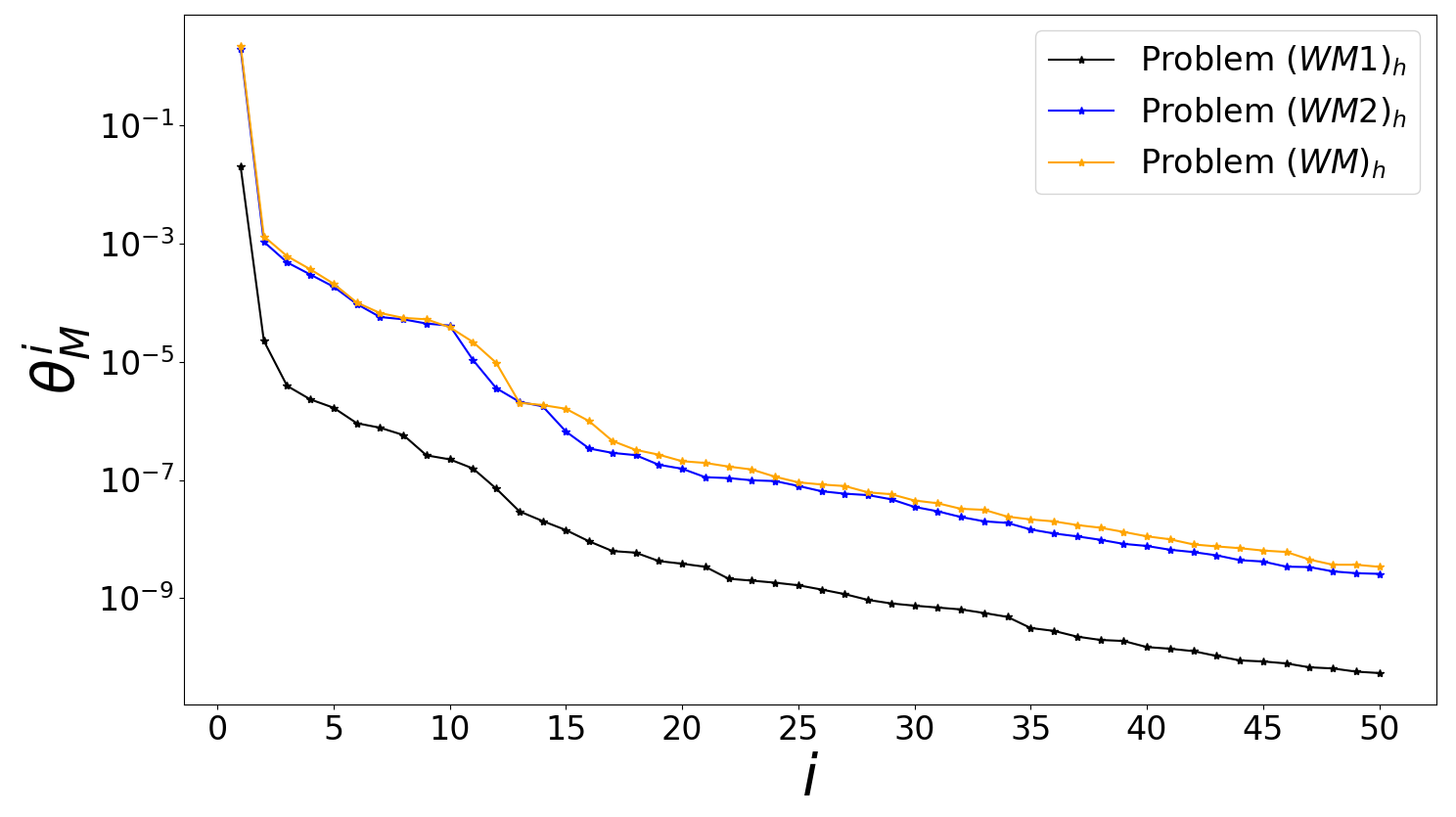}
		\caption[]%
		{{\small \revvv{Numerical} experiment (iv).}}
		\label{mechanical_model_all_parameters_eigenvalue}
	\end{subfigure}
	\caption[]%
	{Mechanical model: \revv{plot of the} eigenvalues \revv{$\lbrace \theta_M^i \rbrace_{i=1}^{50}$ sorted in descending order} for all the \revvv{numerical} experiments considered.}
	\label{eigenvalue_mechanical}
\end{figure}

\revvv{Table \ref{mechanical_training_validation_testing_table} shows the total number of samples provided by the full order model, the number of samples used for training and the one used for testing of ANN}. For all the \revvv{numerical} experiments, the POD space was computed by considering $1000$ snapshots. %POD space was constructed from these snapshots and applying criteria as per equation \eqref{eigenvalue_pod_criteria}.
The eigenvalue plot is shown in Figure \ref{eigenvalue_mechanical}. Like the thermal model, we observe that the \revvv{numerical} experiment (iv), characterized by the larger number of parameters, shows the lowest decay. On the other hand, as expected, among the different mechanical models we consider, the model $(WM)_h$ exhibits the slowest eigenvalues decay including it both thermal and mechanical effects.
%he eigenvalues decay is shown in Fig. \ref{eigenvalue_thermal_active}. We see that the decay related to the experiment (iv) is the slowest. %Therefore, a higher number of basis functions needs to be considered in order to obtain an accurate reconstruction of the field.
%This is due to the fact that in the experiment (iv) we consider a larger number of parameters, so the system exhibits a greater complexity and the modal content is more wide. %is more
%\textcolor{red}{The number of parameter tuples for error analysis and speedup analysis considered were $50$ in both the cases} \michele{I think that here we have the same issue detected for the thermal problem. We should refer eigenvalues decay plot and error/speed-up analysis to the same number of snapshots collected: 50 as for the termal model?}.
%TODO da scrivere nel paragrafo ANN We consider, $2$ hidden layer ANN network. The number of training parameters for the ANN, that is $n_{tr}$, and the number of neurons in the hidden layer, i.e. $H$, are determined by trial and error.

\begin{figure}[H]
        \centering
        \begin{subfigure}[b]{0.475\textwidth}
            \centering
            \includegraphics[width=\textwidth]{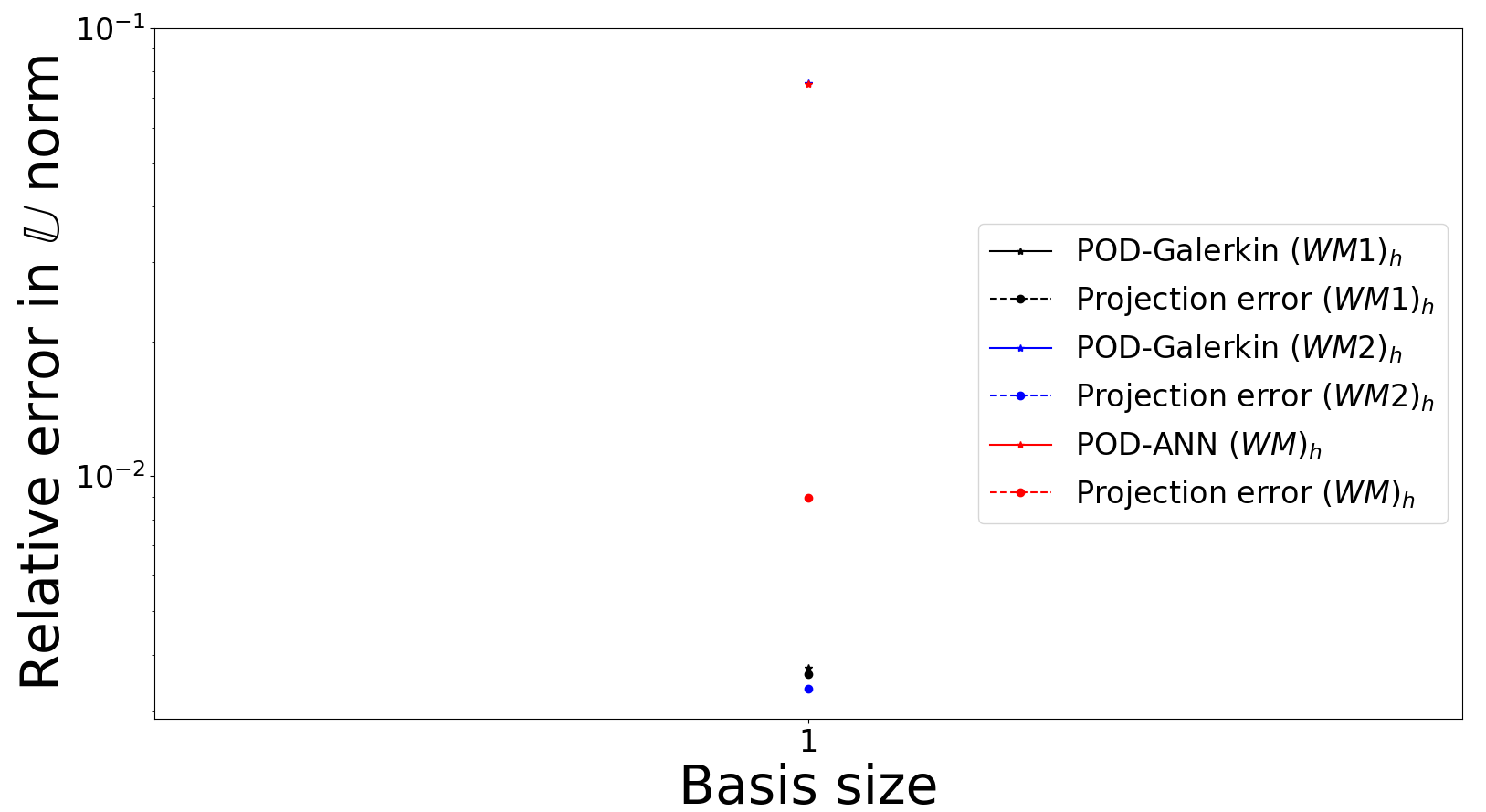}
            \caption[]%
            {{\small \revvv{Numerical} experiment (i). For POD-ANN $n_{tr} = 500$ and $H=60$.}}
            \label{Mechanical_model_4_materials_error}
        \end{subfigure}
        \hfill
        \begin{subfigure}[b]{0.475\textwidth}
            \centering
            \includegraphics[width=\textwidth]{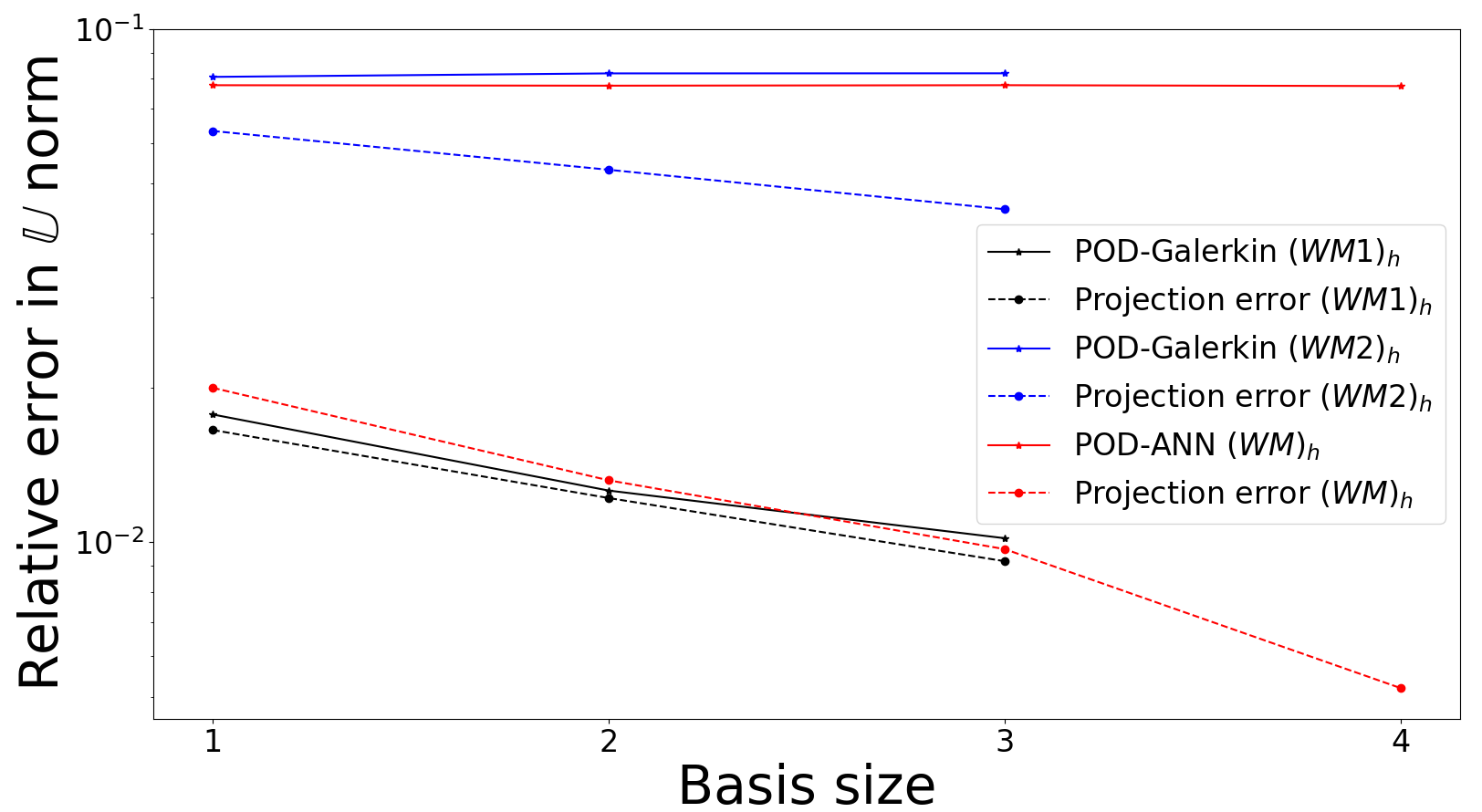}
            \caption[]%
            {{\small \revvv{Numerical} experiment (ii). For POD-ANN $n_{tr} = 500$ and $H=80$.}}
            \label{mechanical_model_4_materials_3_geometric_error}
        \end{subfigure}
        \vskip\baselineskip
        \begin{subfigure}[b]{0.475\textwidth}
            \centering
            \includegraphics[width=\textwidth]{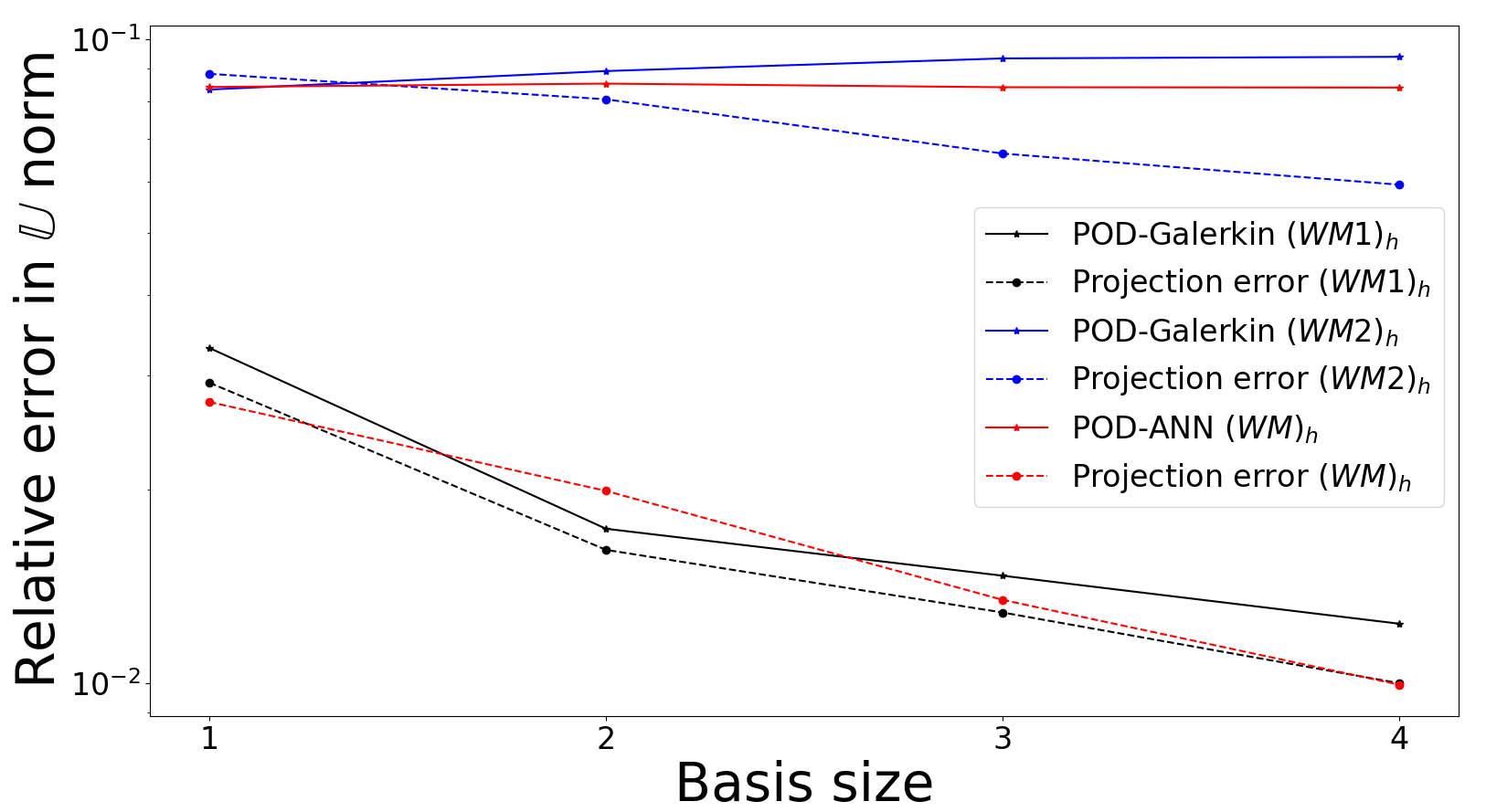}
            \caption[]%
            {{\small \revvv{Numerical} experiment (iii). For POD-ANN $n_{tr} = 1000$ and $H=170$.}}
            \label{mechanical_model_4_materials_6_geometric_error}
        \end{subfigure}
        \hfill
        \begin{subfigure}[b]{0.475\textwidth}
            \centering
            \includegraphics[width=\textwidth]{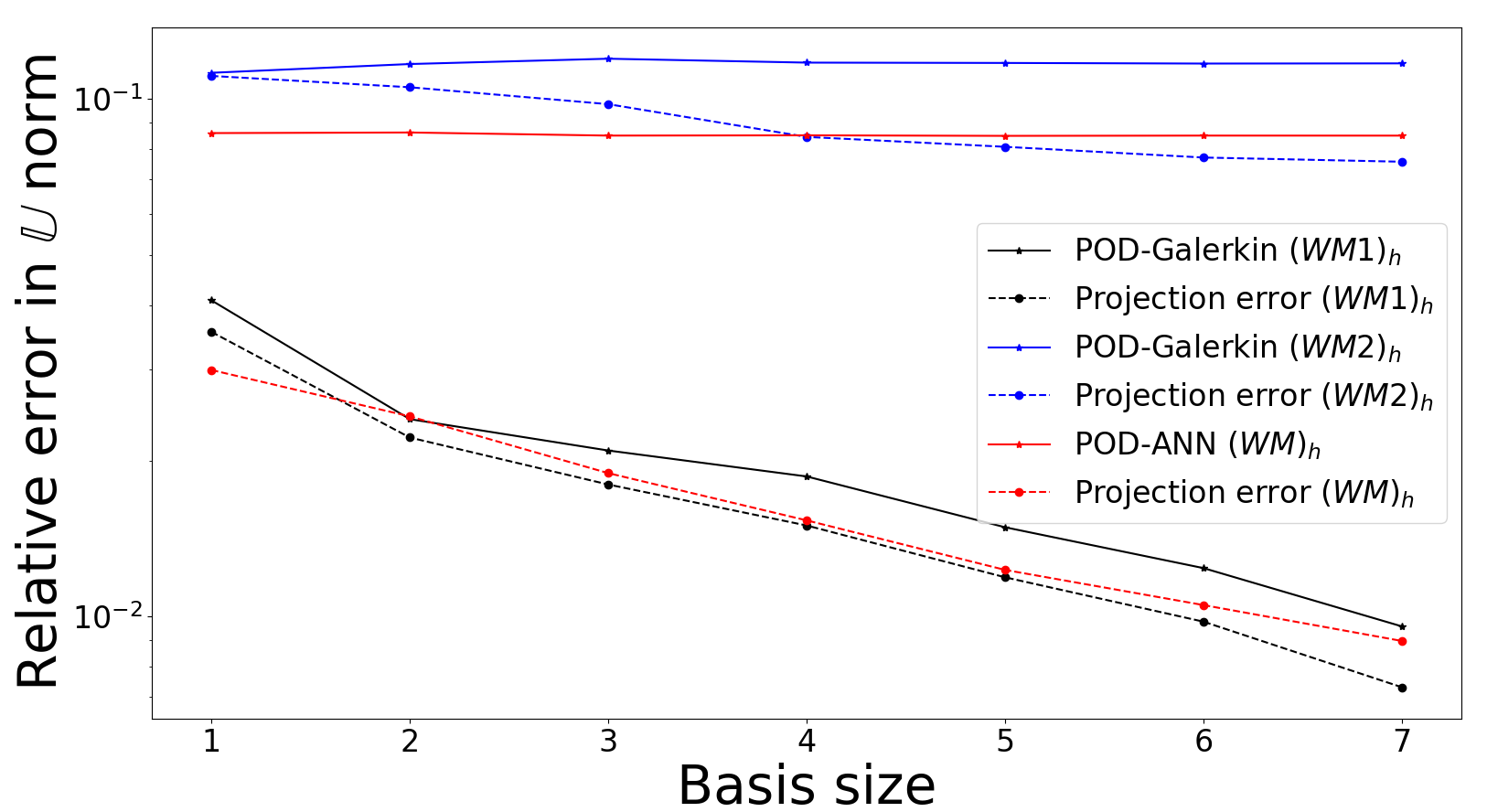}
            \caption[]%
            {{\small \revvv{Numerical} experiment (iv). For POD-ANN $n_{tr} = 2500$ and $H=130$.}}
            \label{mechanical_model_all_parameters_error}
        \end{subfigure}
        \caption[]%
        {Mechanical model: error analysis for POD-G and POD-ANN for all the \revvv{numerical} experiments considered. %\michele{How ever did you choose $n_{tr} = 100$ and $H = 170$? I would choose the values that have allow to obtain the best results for the thermal model (I mean, the thermal model should give guide lines for the mechanical one that is more complex). Still, did you perform a "trial and error" investigation also for this model?.} \nirav{I used Trial and Error approach. Thermal model is not a reference for thermo-mechanical model as the input parameters as well as outputs are different. Also, please note that since size of mechanical system of equations is larger than size of thermal model systme of equations, high $n_{tr}$ will be very costly.} \michele{Sorry Nirav but these figures are very very poor in terms of quality...We mandatory should improve them.}
}
        \label{mechanical_error_analysis_plots}
    \end{figure}

Figure \ref{mechanical_error_analysis_plots} shows the relative error (\ref{rel_error}) both for POD-ANN and POD-G. The projection error (\ref{projection_error}) is also depicted. \revvv{As observed for the thermal model, POD-G is able to provide more accurate results with respect to POD-ANN.}

\begin{figure}[H]
        \centering
        \begin{subfigure}[b]{0.32\textwidth}
            \centering
            \includegraphics[width=\textwidth]{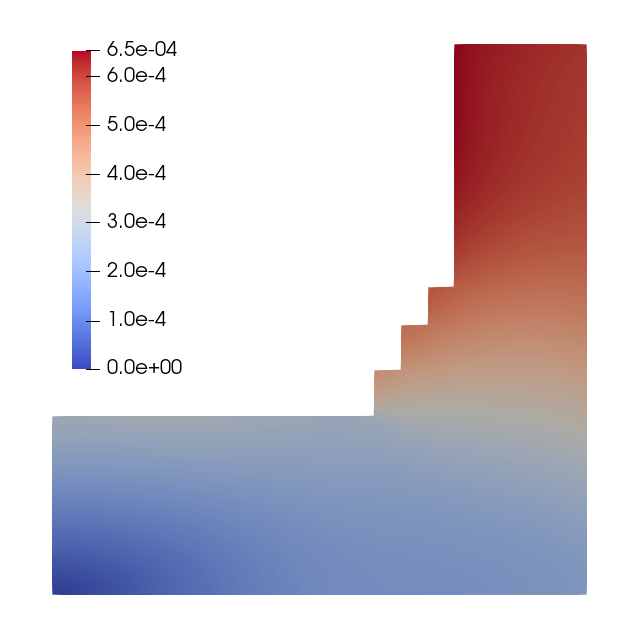}
            \caption[]%
            {{\small FOM solution related to the problem $(WM1)_h$}}
            \label{fem_ref_mechanical}
        \end{subfigure}
       \hfill
       \begin{subfigure}[b]{0.32\textwidth}
            \centering
            \includegraphics[width=\textwidth]{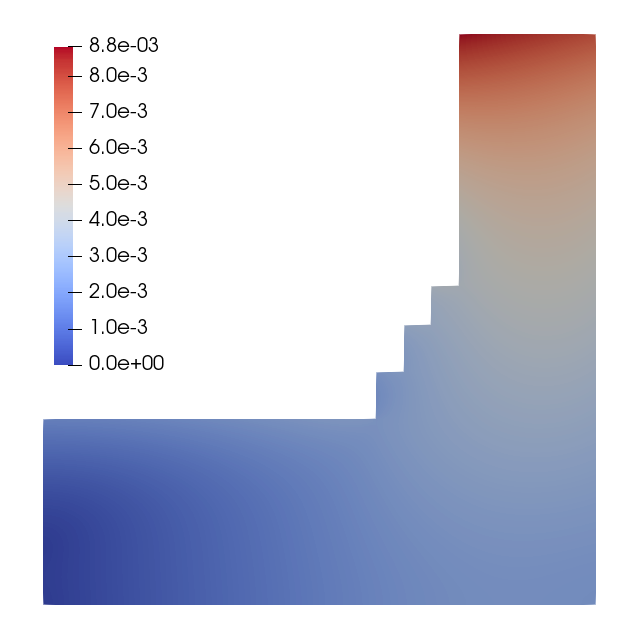}
            \caption[]%
            {{\small FOM solution related to the problem $(WM2)_h$}}
            \label{fem_ref_thermo-mechanical}
        \end{subfigure}
        \hfill
        %\begin{subfigure}[b]{0.32\textwidth}
        %    \centering
        %    \includegraphics[width=\textwidth]{images/podg_ref_thermo-mechanical.png}
        %    \caption[]%
        %    {{\small FOM solution related to the problem $(WM)_h$\newline}}
        %    \label{podg_ref_thermo-mechanical}
        %\end{subfigure}
        % \hfill
        \begin{subfigure}[b]{0.32\textwidth}
            \centering
            \includegraphics[width=\textwidth]{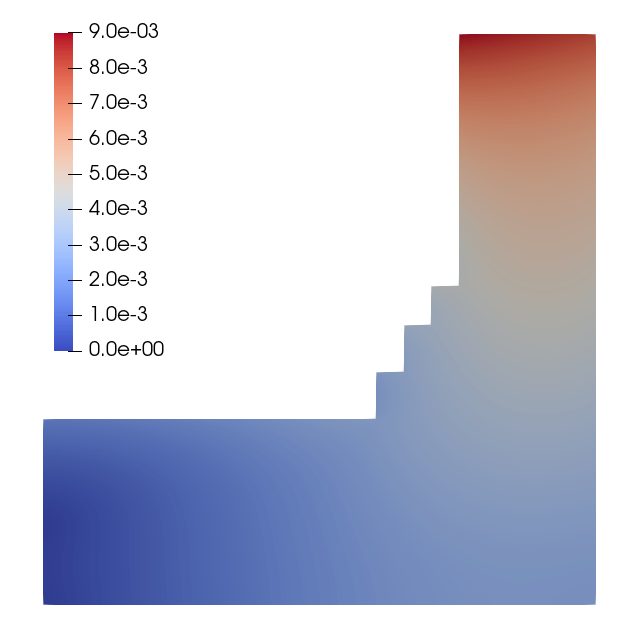}
            \caption[]%
            {{\small FOM solution related to the problem $(WM)_h$}}
            \label{fem_ref_thermo_mechanical}
        \end{subfigure}
        \hfill
        \begin{subfigure}[b]{0.32\textwidth}
            \centering
            \includegraphics[width=\textwidth]{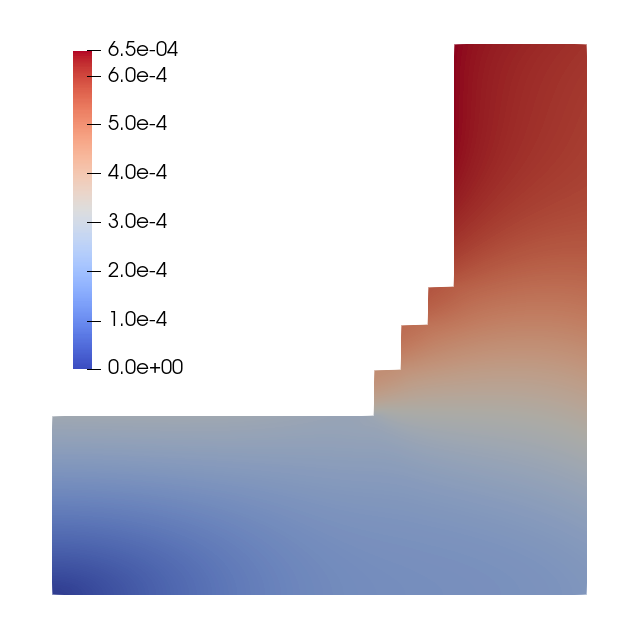}
            \caption[]%
            {{\small POD-G solution related to the problem $(WM1)_h$}}
            \label{podg_ref_mechanical}
        \end{subfigure}
        \hfill
        %\begin{subfigure}[b]{0.32\textwidth}
        %    \centering
        %    \includegraphics[width=\textwidth]{images/podg_error_ref_mechanical.png}
        %    \caption[]%
        %    {{\small Spatial distribution of error : Mechanical model}}
        %    \label{podg_error_ref_mechanical}
        %\end{subfigure}
        %\vskip\baselineskip
        %\begin{subfigure}[b]{0.32\textwidth}
        %    \centering
        %    \includegraphics[width=\textwidth]{images/podg_error_ref_thermo-mechanical.png}
        %    \caption[]%
        %    {{\small Spatial distribution of error : Coupling model}}
        %    \label{podg_error_ref_thermo-mechanical}
        %\end{subfigure}
        %\vskip\baselineskip
        \begin{subfigure}[b]{0.32\textwidth}
            \centering
            \includegraphics[width=\textwidth]{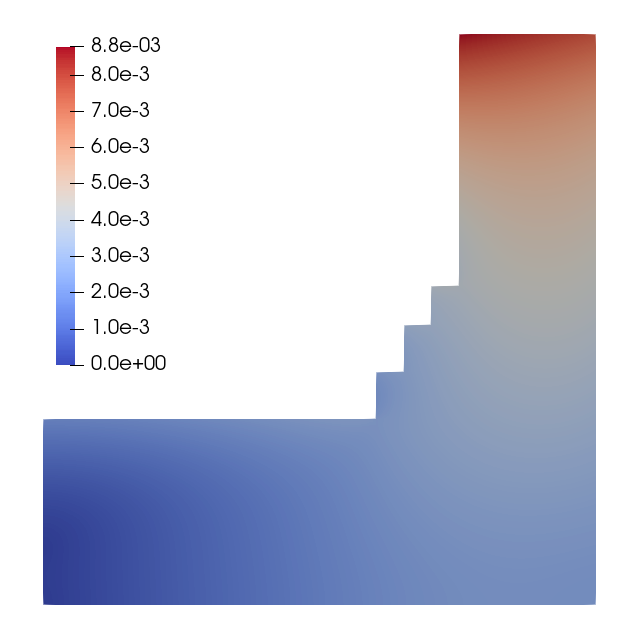}
            \caption[]%
            {{\small POD-G solution related to the problem $(WM2)_h$}}
            \label{podg_ref_thermo-mechanical}
        \end{subfigure}
         \hfill
        \begin{subfigure}[b]{0.32\textwidth}
            \centering
            \includegraphics[width=\textwidth]{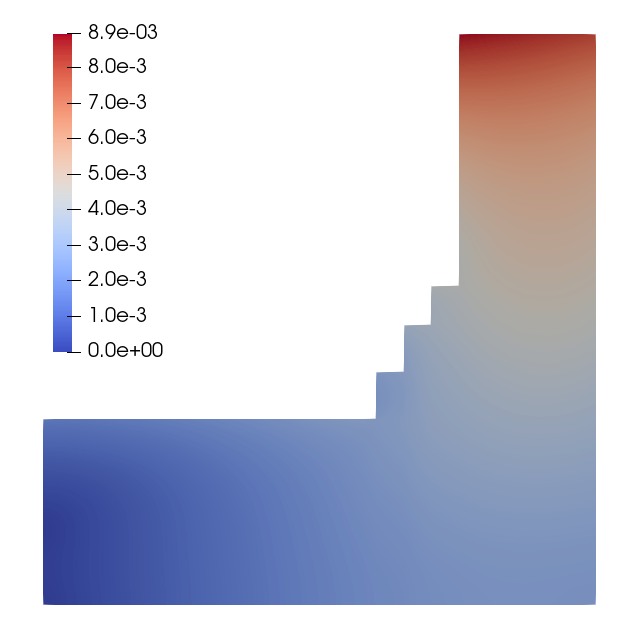}
            \caption[]%
            {{\small POD-ANN solution related to the problem $(WM)_h$}}
            \label{podann_ref_thermo-mechanical}
        \end{subfigure}
        \hfill
        %\begin{subfigure}[b]{0.32\textwidth}
        %    \centering
        %    \includegraphics[width=\textwidth]{images/podann_error_ref_thermo_mechanical.png}
        %    \caption[]%
        %    {{\small Spatial distribution of error : thermo-mechanical model}}
        %    \label{podann_error_ref_thermo-mechanical}
        %\end{subfigure}
        \caption[]%
        {Mechanical model: comparison between the displacement (in m) computed by FOM and by the POD-G and POD-ANN methods related to the \revvv{numerical} experiment (iv) for $\Xi = \lbrace 2.365, 0.6, 0.6, 0.5, 3.2, 14.10, 8.50, 9.2, 9.9, 10.6, 10, 2.08e9, 1.39e9, 1e-6 \rbrace$. We consider $7$ POD modes. For POD-ANN, we set $n_{tr} = 2500$ and $H=130$.}
        \label{solutions_ref_mech}
    \end{figure}

Figure \ref{solutions_ref_mech} shows the qualitative comparison between the computed FOM and MOR related to the \revvv{numerical} experiment (iv) for the parameters tuple %\michele{how many POD basis? For POD-ANN $n_{tr}=?$ and $H=?$?}. \nirav{$6$ POD basis with $n_{tr} = 2500$ and $H=130$. I added these values in caption of plots.} % \michele{inserire un commento in piu}. %The reduced basis space contained $4$ basis functions \michele{inserire un commento in piu}.  %the full order model solution, the full order solution projected on reduced basis space, POD-Galerkin solution and POD-ANN solution (figures \ref{solutions_ref_mech}, \ref{solutions_para_mech}) for the problem presented in the section \ref{benchmark_actual} for the parameter values :
\begin{gather*}
\Xi = \lbrace 2.365, 0.6, 0.6, 0.5, 3.2, 14.10, 8.50, 9.2, 9.9, 10.6, 10, 2.08e9, 1.39e9, 1e-6 \rbrace \ . %\\
%\Xi_4 = \lbrace 2.365, 0.6, 0.6, 0.45, 3.2, 14.10, 8.30, 9.2, 9.9, 10.6, 10, 2.08e9, 1.39e9, 1e-6 \rbrace \ .
\end{gather*}
We use 7 POD basis. The POD-ANN solution was computed with $n_{tr} = 2500$ and $H = 130$. We could observe that both MOR approaches are able to provide a good reconstruction of the displacement field. \revv{In order to justify our choice to consider separately $(WM1)_h$ and $(WM2)_h$ in the Galerkin projection framework, we also highlight that, as shown in Figure \ref{solutions_ref_mech}, the scale difference between their displacements (derived from mechanical loads for the first and from the thermal ones for the second) is of one order of magnitude. %This justifies our choice to consider separately $(WM1)_h$ and $(WM2)_h$ in the Galerkin projection framework \cite{greedy_paper,mor_book,variational_multiscale}. 
Thus, the use of the model $(WM)_h$ could lead to a less accurate reconstruction of the displacement field.}

Finally, we briefly discuss the efficiency of our MOR approach. \rev{We report in Table \ref{mechanical_offline_time} some
estimations related to the offline time for all the numerical experiments carried out. We observe that, unlike the thermal model, the time taken by POD is comparable for all the numerical experiments. This is not surprising because in this case we consider the same number of snapshots for the computation of the reduced space for all the numerical experiments. Like the thermal model, the ANN training is faster for the numerical experiment (i), probably because of the minor complexity with respect to the other numerical experiments. Unlike the thermal model, the total offline cost of the POD-ANN method is comparable with the one related the POD-G method. This is because we train two models for the POD-G method that take a similar amount of time as training one model for the POD-ANN method.} We report in Table \ref{mechanical_rb_time_plots} the online time related to the POD-G and POD-ANN methods for all the \revvv{numerical} experiments carried out. %\michele{Nirav, how long does a FOM simulation take?} \nirav{$(WM1)_h = 0.5 seconds, (WM2)_h = (WM)_h = 1 second$}.
\revvv{Like the thermal model, the online time of POD-G method increases significantly in presence of geometric parameters by moving from $8e-4$ s (\revvv{numerical} experiment (i)) to $2.6$/$6.9e-2$ s (\revvv{numerical} experiments (ii)-(iv)) for the model $(WM1)_h$ and from $4.5e-2$ s (\revvv{numerical} experiment (i)) to $1.9$/$2.6e-1$ s (\revvv{numerical} experiments (ii)-(iv)) for the model $(WM2)_h$. We could observe that the online time taken by the model $(WM2)_h$ is significantly greater than that taken by the model $(WM1)_h$. This is expected because for the model $(WM2)_h$ a reduced basis approximation of temperature needs to be computed due to the thermo-mechanical coupling. On the other hand, the POD-ANN, that does not need reduced basis approximation of temperature thanks to its non intrusive nature, is able to provide a higher computational efficiency. Moreover, like the thermal model, the POD-ANN online time remains relatively constant for all the \revvv{numerical} experiments under investigations, around $5e-4$, by showing a low sensitivity at varying of the kind and number of parameters considered.} %So the computational efficiency of POD-ANN is much higher (of almost two order of magnitude) than POD-G when geometrical parametrization is considered.

%online phase of POD-Galerkin thermo-mechanical model is slower than online phase of POD-Galerkin mechanical model. This is on the expected lines, as for the thermo-mechanical model a reduced basis approximation of temperature needs to be computed. Also, by eliminating the use of reduced basis approximation of temperature, POD-ANN is able to achieve higher speedup.

%As can be seen, the online time of POD-G method decreases significantly in presence of geometric parameters by moving from $7e-4$ s (Experiment 1) to $1.3e-2$ s (Experiments 2-4). On the other hand, the time taken by POD-ANN online stage remains relatively constant for all the experiments under investigations, around $5e-4$. So the computational efficiency of POD-ANN is much higher (of almost two order of magnitude) than POD-G when geometrical parametrization is considered.

%The error analysis with different $n_{tr}$ and $H$ is shown and comparison with POD-Galerkin method is shown in the figure \ref{mechanical_error_analysis_plots}. From the error analysis,
%\textcolor{red}{it can be seen that for POD-Galerkin and POD-ANN coupling is the main cause error. Also, by considering mechanical and thermo-mechanical models for POD-Galerkin, higher accuracy is obtained as compared to treat single mechanical model. On the other hand, by considering only mechanical model, POD-ANN eliminates the use of reduced basis temperature approximation for the reduced basis approximation of the displacement}.

\begin{table}[ht]
\centering
\begin{tabular}{|c|c|c|c|}
\hline
& $n_{tr}$ & Training & Testing \\ %Testing cases \\
\hline
\multirow{1}{16em}{Numerical experiment (i)} & 500 & 350 & 150 \\%&  50 \\
\hline
\multirow{1}{16em}{Numerical experiment (ii)} & 500 & 350 & 150 \\%& 100 \\
\hline
\multirow{1}{16em}{Numerical experiment (iii)} & 1000 & 700 & 300 \\%& 100 \\
\hline
\multirow{1}{16em}{Numerical experiment (iv)} & 2500 & 1750 & 750 \\%& 100 \\
\hline
\end{tabular}
\caption{\revvv{Mechanical model: number of total samples $n_{tr}$ by FOM and number of samples used for training and testing of ANN.}}
\label{mechanical_training_validation_testing_table}
\end{table}

\begin{table}
\begin{subtable}[h]{\textwidth}
\centering
\begin{tabular}{|c|c|c|c|c|c|c|}
\hline
 & \multicolumn{2}{c|}{$t_{off}^{POD-G}$} & \multicolumn{2}{c|}{$t_{POD}$} & \multicolumn{2}{c|}{$t_{FOM}$}\\
\hline
& $(WM1)_h$ & $(WM2)_h$ & $(WM1)_h$ & $(WM2)_h$ & $(WM1)_h$ & $(WM2)_h$ \\
\hline
Numerical experiment (i) & $\approx 4.2e2$ & $\approx 6.2e2$ & 1.6e1 & 1.8e1 & \multirow{4}{2em}{4e-1} & \multirow{4}{2em}{6e-1}\\
\cline{1-5}
Numerical experiment (ii) & $\approx 4.2e2$ & $\approx 6.2e2$ & 1.6e1 & 1.7e1 & & \\
\cline{1-5}
Numerical experiment (iii) & $\approx 4.2e2$ & $\approx 6.2e2$ & 1.8e1 & 1.7e1 & & \\
\cline{1-5}
Numerical experiment (iv) & $\approx 4.2e2$ & $\approx 6.2e2$ & 1.7e1 & 1.7e1 & & \\
\hline
\end{tabular}
\caption{Problem $(WM1)_h$ and Problem $(WM2)_h$ }
%\label{}
\end{subtable}
\hspace{0.5cm}
\begin{subtable}[h]{\textwidth}
\centering
\begin{tabular}{|c|c|c|c|c|c|}
\hline
 & $t_{off}^{POD-ANN}$ & \multicolumn{1}{c|}{$t_{POD}$} & $t_{tr}$ & \multicolumn{1}{c|}{$t_{FOM}$} & \multicolumn{1}{c|}{$t_{proj}$}\\
\hline
& $(WM)_h$ & $(WM)_h$ & $(WM)_h$ & $(WM)_h$ & $(WM)_h$\\
\hline
Numerical experiment (i) & $\approx 1.1e3$ & 1.8e1 & 3.4 & \multirow{4}{2em}{7e-1} & 8.1e-4 \\
\cline{1-4} \cline{6-6}
Numerical experiment (ii) & $\approx 1.1e3$ & 1.8e1 & 3.7 & & 9.2e-4\\
\cline{1-4} \cline{6-6}
Numerical experiment (iii) & $\approx 1.5e3$ & 1.6e1 & 2.9e1 & & 1.0e-3 \\
\cline{1-4} \cline{6-6}
Numerical experiment (iv) & $\approx 2.5e3$ & 1.8e1 & 7.6e1 & & 9.1e-4\\
\hline
\end{tabular}
\caption{Problem $(WM)_h$}
%\label{}
\end{subtable}
\caption{\rev{Mechanical model: time (in s) taken by (i) the entire offline stage ($t_{off}$), (ii) the computation of the POD modes ($t_{POD}$), (iii) the training of ANN ($t_{tr}$), (iv) the computation of a FOM solution ($t_{FOM}$) and (v) the projection of a FOM solution on the POD space ($t_{proj}$). Concerning POD-ANN, we use $n_{tr} = 500, H = 60$ for the \revvv{numerical} experiment (i), $n_{tr} = 500, H=80$ for the \revvv{numerical} experiment ii), $n_{tr} = 1000, H=170$ for the \revvv{numerical} experiment (iii) and  $n_{tr} = 2500, H = 130$ for the \revvv{numerical} experiment (iv).}}
%\caption{\rev{Mechanical model: Offline time taken (in s) for (i) eigenvalue decomposition of the snapshot matrix ($t_{POD}$) (ii) training of ANN after computation of training set ($t_{tr}$) (iii) computation of one snapshot by solving FEM formulation ($t_{FEM}$).}}
\label{mechanical_offline_time}
\end{table}

% \begin{table}
% \centering
% \begin{tabular}{|c|c|c|c|}
% \hline
% Problem & $(WM1)_h$ & $(WM2)_h$ & $(WM)_h$\\
% \hline
% Offline time & 4e-1 & 6e-1 & 7e-1\\
% \hline
% \end{tabular}
% \caption{\revv{Mechanical model: offline time (in s).}}
% \label{mechanical_offline_time}
% \end{table}

\begin{table}[H]
\centering
\resizebox{\textwidth}{!}{\begin{tabular}{|c|c|c|c|c|}
\hline
 & Basis size &  POD-G $(WM1)_h$ & POD-G $(WM2)_h$ & POD-ANN $(WM)_h$\\
%\hline
%& & POD-Galerkin (Mechanical) & POD-Galerkin (Coupling) & POD-ANN \\
\hline
\multirow{1}{16em}{\revvv{Numerical} experiment (i)} & 1 & 8e-4 & 4.5e-2 & 6.7e-4\\
\hline
\multirow{1}{16em}{\revvv{Numerical} experiment (ii)} %& 1 & 2.68e-2 & 1.93e-1 & 6.56e-4\\
%& 2 & 2.63e-2 & 1.89e-1 & 5.79e-4\\
& 3 & 2.6e-2 & 1.9e-1 & 5.3e-4\\
\hline
\multirow{1}{16em}{\revvv{Numerical} experiment (iii)} %& 1 & 5.26e-2 & 2.11e-1 & 5.35e-4\\
%& 2 & 5.29e-2 & 2.10e-1 & 5.24e-4\\
%& 3 & 5.64e-2 & 2.09e-1 & 5.24e-4\\
& 4 & 5.4e-2 & 2.1e-1 & 5.2e-4\\
\hline
\multirow{1}{16em}{\revvv{Numerical} experiment (iv)} %& 1 & 7.11e-2 & 2.35e-1 & 4.99e-4\\
%& 2 & 6.67e-2 & 2.56e-1 & 4.88e-4\\
%& 3 & 6.35e-2 & 2.55e-1 & 4.88e-4\\
%& 4 & 6.50e-2 & 2.56e-1 & 4.88e-4\\
%& 5 & 6.48e-2 & 2.56e-1 & 4.88e-4\\
%& 6 & 6.36e-2 & 2.57e-1 & 4.88e-4\\
& 7 & 6.9e-2 & 2.6e-1 & 4.9e-4\\
\hline
\end{tabular}}
\caption{Mechanical model: online time (in s) for all the \revvv{numerical} experiments under investigation. Concerning POD-ANN, we use $n_{tr} = 500, H = 60$ for the \revvv{numerical} experiment (i), $n_{tr} = 500, H = 80$ for the \revvv{numerical} experiment (ii), $n_{tr} = 1000, H = 170$ for the \revvv{numerical} experiment (iii) and  $n_{tr} = 2500, H = 130$ for the \revvv{numerical} experiment (iv).}
\label{mechanical_rb_time_plots}
\end{table}

\section{Some concluding remarks}\label{remark_chapter}
%This work presents two different MOR techniques, POD-G and POD-ANN,
In this work we propose a computational pipeline to obtain fast and reliable numerical simulations for one-way coupled steady state linear thermo-mechanical problems in a finite element environment. The test case is referred to a relevant industrial problem related to the investigation of the thermo-mechanical phenomena occurring in blast furnace heart walls. After introducing the main theoretical features of FOM, we detect customized benchmarks for the validation of its numerical implementation. Then we present our MOR framework: we apply POD for the computation of reduced basis space whilst for the evaluation of the modal coefficients we use two different methodologies, the one based on a classic Galerkin projection (POD-G) and the other one based on artificial neural networks (POD-ANN). \revvv{We found that POD-G is generally more accurate than POD-ANN although POD-ANN exhibits a very higher efficiency, especially when geometric parameters are considered}. % as relevant to industrial application to blast furnace hearth. The finite element formulation was applied to benchmark cases and real industrial application. The paper also discussed the considerations related to the application and the performance of the model order reduction method as applied to governing coupled parametric partial differential equation. In particular, we explained the POD-Galerkin and POD-Artificial Neural Network approaches and studied the computational performance related aspects.
\revvv{The higher efficiency of POD-ANN in the case of mechanical model can also be attributed to the fact that the computation of reduced basis approximation of temperature field is not required to compute the reduced basis approximation of the displacement field.}

We believe that insights given in this work could help to develop advanced numerical tools in order to deal with complex industrial problems. As a follow-up of this work, we are going to enhance training capacity of the deep learning methods, \rev{i.e. to reduce the offline cost, as well as to} move towards more complex thermo-mechanical problems involving heterogeneity, orthotropy and non-linearity. Some \rev{preliminary} efforts, in \rev{the latter} direction, have been already carried out \rev{both} at full order \cite{shah_ecmi_paper} \rev{and reduced order level \cite{nirav_doctoral_thesis}.} \rev{Concerning the POD-Galerkin method, it is expected that the computational efficiency furthermore decreases due to absence of affine expansion for assembling system of equations. Regarding the accuracy of the method, the interpolation of operators could introduce an additional source of error. On the other hand, the POD-ANN method can be properly set to take care of non-linearities by increasing the number of hidden layers and/or their depth. However, the accuracy
and the computational efficiency are not likely to change to the significant extent \cite{pod_nn_paper,federico_paper}.}

\section*{Acknowledgements}

%This work has been supported by the European Union{\textquotesingle}s Horizon 2020 research and innovation programme under the Marie Sk\l odowska-Curie Grant Agreement No. 765374. The authors and the institutes collaborating on this work do not intend and do not encourage the dual use of this work.
We are grateful to Dr. Federico Pichi (SISSA mathLab) for insights and crucial support in the numerical implementation of artificial neural network.

We would like to acknowledge the financial support of the European Union under the Marie Sklodowska-Curie Grant Agreement No. 765374. We also acknowledge the partial support by the European Union Funding for Research and Innovation - Horizon 2020 Program - in the framework of European Research Council Executive Agency: Consolidator Grant H2020 ERC CoG 2015 AROMA-CFD project 681447 ``Advanced Reduced Order Methods with Applications in Computational Fluid Dynamics'' and INDAM-GNCS project ``Advanced intrusive and non-intrusive model order reduction techniques and applications'', 2019. This work was also partially supported by FEDER and Xunta de Galicia [grant numbers ED431C 2017/60, ED431C 2021/15], and the Agencia Estatal de Investigación [PID2019-105615RB-I00/AEI/10.13039/501100011033].

This work has focused exclusively on civil applications. It is not to be used for any illegal, deceptive, misleading or unethical purpose or in any military applications. This includes any application where the use of this work may result in death, personal injury or severe physical or environmental damage.

%% The Appendices part is started with the command \appendix;
%% appendix sections are then done as normal sections
%% \appendix

%% \section{}
%% \label{}

%% If you have bibdatabase file and want bibtex to generate the
%% bibitems, please use
%%
%%  \bibliographystyle{elsarticle-harv} 
%%  \bibliography{<your bibdatabase>}

%% else use the following coding to input the bibitems directly in the
%% TeX file.

%\begin{thebibliography}{00}

%% \bibitem[Author(year)]{label}
%% Text of bibliographic item

%\bibitem[ ()]{}

\appendix

\section{Validation of the full order model}\label{appendix}
In this section we verify the numerical implementation of the FOM introduced in Secs. \ref{weak_form_chapter} and \ref{fem_chapter}. %Benchmark tests proposed are not taken from the literature but are built on our own. We consider analytical expressions for temperature and displacement, calculate corresponding data including boundary conditions and source terms to be coupled with the governing equations in order to detect the full order problem for which the selected analytical relationship are solutions and then we numerically solve the problem to carry out a comparison.
%In Sec. \ref{benchmark_test_section}, we report details about the construction of the computational grid used for numerical simulations. Then proper benchmark tests for $(WT)_h$, $(WM)_h$, $(WM1)_h$ and $(WM2)_h$ problems are investigated in \revv{\ref{benchmark_test}}. 
All the FOM computations have been performed by using the python finite element library FEniCS \cite{FEniCS}. %\michele{Nirav, could you insert the reference of the sec.? Thanks}.%The results of the numerical simulation corresponding to the thermo-mechanical deformation suffered by the wall of the blast furnace under real conditions are also presented, although, for the time being, without considering the different materials that make up its geometry.
%Given that a first objective is to apply a method of order reduction for the numerical resolution of the real problem, we focus first, in Sec. \michele{reference, thanks}, on the difficulty linked to the geometric shape of the vertical section.

%\subsection{Computational domain and mesh}\label{benchmark_test_section}
%The considered division of the computational domain into triangular subdomains verifies the assumption of affine parameter dependence, as will be explained later in Section \ref{geometric_parameter_subsection}.
%whilst we obtained an average quality over all elements of$0.95$. It is important to ensure that minimum mesh quality is sufficiently far from zero.

%\section{Benchmark tests}\label{benchmark_test}
%Benchmark tests proposed are not taken from the literature but are built on our own.
We use a mesh of $\omega$ containing $121137$ triangular elements and $61147$ vertices. %\nirav{The mesh size is measured as distance between vertices of an element of the mesh}. 
The minimum mesh size is $0.011$ m and the maximum one is $0.045$ m. Its minimum quality is $q_e = 0.25$ (eq. \ref{mesh_quality_formula}).

The pipeline that we follow for the design of reliable benchmark tests to be used for the FOM validation consists of three steps:

\begin{itemize}
\item We set analytical expressions for temperature and displacement.
\item We calculate corresponding model data, including boundary conditions and source terms, in order to identify the FOM for which the analytical relationships are solutions.
\item Finally, we numerically solve the problem and compare the computational solutions with the analytical ones.
\end{itemize}

We consider the physical properties reported in Table \ref{material_table} for all numerical simulations shown in this section. % unless otherwise stated. %Here, Yield stress \michele{Do we use this datum in our treatment? I seem not} refers to the stress at point in the stress-strain diagram where sudden strain takes place without any significant increase in the stress. %: \michele{maybe, it's better to organize the following scheme as a table} \nirav{I could not understand this. Could you clarify?}

\begin{table}[H]
\centering
\begin{tabular}{|l|l|}
\hline
\textbf{Property} & \textbf{Value} \\
\hline
Thermal conductivity $k$ & $10 \frac{W}{mK}$ \\[1em]
\hline
Convection coefficient $h_{c,-}$ & $2000 \frac{W}{m^2K}$\\[1em]
\hline
Convection coefficient $h_{c,f}$ & $200 \frac{W}{m^2K}$\\[1em]
\hline
Convection coefficient $h_{c,out}$ & $2000 \frac{W}{m^2K}$\\[1em]
\hline
Young's modulus $E$ & $5e9 Pa$\\
\hline
Poisson's ratio $\nu$ & $0.2$\\
%\hline
%Yield stress $\sigma_Y$ \michele{What is Yield stress?}\nirav{Done} & 50 MPa \\
\hline
Thermal expansion coefficient $\alpha$ & $10^{-6}/{K}$\\
\hline
Reference temperature $T_0$ & $298 K$\\
\hline
%Density of molten metal \michele{(I seem that this property does not appear in the description of the model)} $\rho_m$ & $7460 \frac{Kg}{m^3}$\\[1em]
%\hline
Gravitational acceleration $g$ & $9.81 \frac{m}{s^2}$\\[1em]
\hline
\end{tabular}
\caption{Physical properties values used for the FOM benchmark tests.}
\label{material_table}
\end{table}

\subsection{Thermal model}\label{thermal_benchmark_section}

We consider the following analytical expression for the temperature, % as analytical temperature, the known temperature,
\begin{equation}\label{benchmark_analytical_temperature}
T_a(r,y) = C' r^2 y, \quad  \text{with} \ C'=1 K/m^3 \ .
\end{equation}

Then: % the corresponding data for the thermal model $(T1)$ (eqs. \eqref{summarised_strong_energy_equations}-\eqref{summarised_strong_energy_equations_boundary}) are obtained as follows: %the corresponding data for the thermal model $(T1)$, defined by \eqref{summarised_strong_energy_equations}-\eqref{summarised_strong_energy_equations_boundary}, are calculated. Then, the computed temperature, obtained by solving the $(WT1)_h$ problem, is compared with the known analytical temperature $T_a$.

%Indeed,

\begin{itemize}
\item The corresponding source term $Q$ is obtained by using eq. \eqref{summarised_strong_energy_equations},
\begin{equation}\label{thermal_benchmark_source_term}
Q(r,y) = -k \frac{\partial^2 T_a}{\partial r^2} -k \frac{\partial^2 T_a}{\partial y^2} - \frac{k}{r} \frac{\partial T_a}{\partial r} = -4C'ky \ ,
\end{equation}

\item The heat flux $q_+$, as well as the temperatures $T_f, T_{out}$ and $T_-$, are derived from eq. \eqref{summarised_strong_energy_equations_boundary},
%\begin{equation}\label{thermal_benchmark_q_plus}
%q_+(r,y) = - k \frac{\partial T_a}{\partial y} = - C'kr^2 \ ,
%\end{equation}
\begin{subequations}
\begin{flalign}
\text{on $\gamma_+$ } & : q_+(r,y) = - k \frac{\partial T_a}{\partial y} = - C'kr^2 \ , \label{thermal_benchmark_q_plus} \\
\text{on} \ \gamma_{sf} \ & : \ T_f = T_a + \frac{k}{h_{c,f}} \left( \frac{\partial T_a}{\partial r} n_r + \frac{\partial T_a}{\partial y} n_y \right) \nonumber \\ & \ \ \ \ \ \ \ = C'r^2y + C' \frac{k}{h_{c,f}} (2ryn_r + r^2n_y) \ , \label{thermal_benchmark_t_f} \\
\text{on} \ \gamma_{out} \ & : \ T_{out} = T_a + \frac{k}{h_{c,out}} \frac{\partial T_a}{\partial r} = C'r^2y + C' \frac{2ryk}{h_{c,out}} \ , \label{thermal_benchmark_t_out} \\
\text{on} \ \gamma_- \ & : \ T_- = T_a - \frac{k}{h_{c,-}} \frac{\partial T_a}{\partial y} = C'r^2y - C' \frac{r^2 k}{h_{c,-}} \  , \label{thermal_benchmark_t_minus}
\end{flalign}
\end{subequations}
and it is verified that
\begin{subequations}
\begin{flalign}
\text{on} \ \gamma_s \ & : \ \frac{\partial T_a}{\partial r} = 0. \label{thermal_symmetry_boundary}
\end{flalign}
\end{subequations}

%\item the temperatures $T_f, T_{out}$ and $T_-$ are calculated from eq. \eqref{summarised_strong_energy_equations_boundary},
%\begin{subequations}
%\begin{flalign}
%\text{on} \ \gamma_{sf} \ & : \ T_f = T_a + \frac{k}{h_{c,f}} \left( \frac{\partial T_a}{\partial r} n_r + \frac{\partial T_a}{\partial y} n_y \right) \nonumber \\ & \ \ \ \ \ \ \ = C'r^2y + C' \frac{k}{h_{c,f}} (2ryn_r + r^2n_y) \ , \label{thermal_benchmark_t_f} \\
%\text{on} \ \gamma_{out} \ & : \ T_{out} = T_a + \frac{k}{h_{c,out}} \frac{\partial T_a}{\partial r} = C'r^2y + C' \frac{2ryk}{h_{c,out}} \ , \label{thermal_benchmark_t_out} \\
%\text{on} \ \gamma_- \ & : \ T_- = T_a - \frac{k}{h_{c,-}} \frac{\partial T_a}{\partial y} = C'r^2y - C' \frac{r^2 k}{h_{c,-}} \  , \label{thermal_benchmark_t_minus}
%\end{flalign}
%\end{subequations}
%and it is verified that
%\begin{subequations}
%\begin{flalign}
%\text{on} \ \gamma_s \ & : \ \frac{\partial T_a}{\partial r} = 0.
%\end{flalign}
%\end{subequations}

\end{itemize}

\begin{figure}[H]
\begin{minipage}{.3\textwidth}
\includegraphics[width=\textwidth]{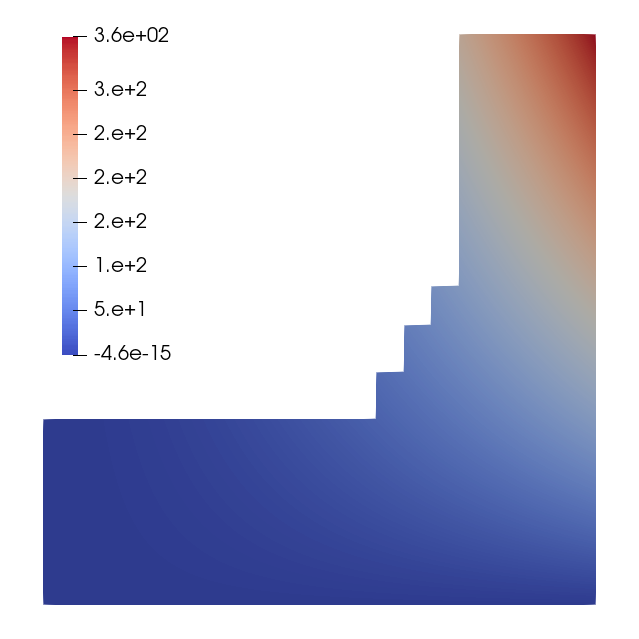}
%\caption*{$T_h (K)$}
\end{minipage}
\begin{minipage}{.3\textwidth}
\includegraphics[width=\textwidth]{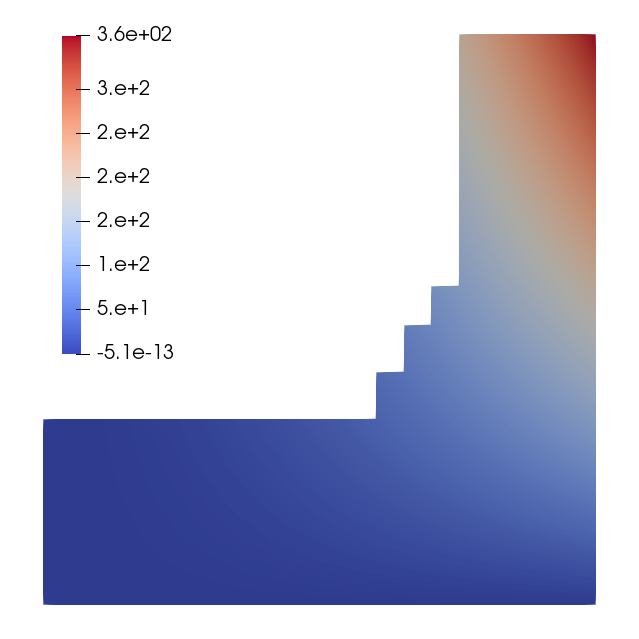}
%\caption*{$T_a (K)$}
\end{minipage}
\begin{minipage}{.3\textwidth}
\includegraphics[width=\textwidth]{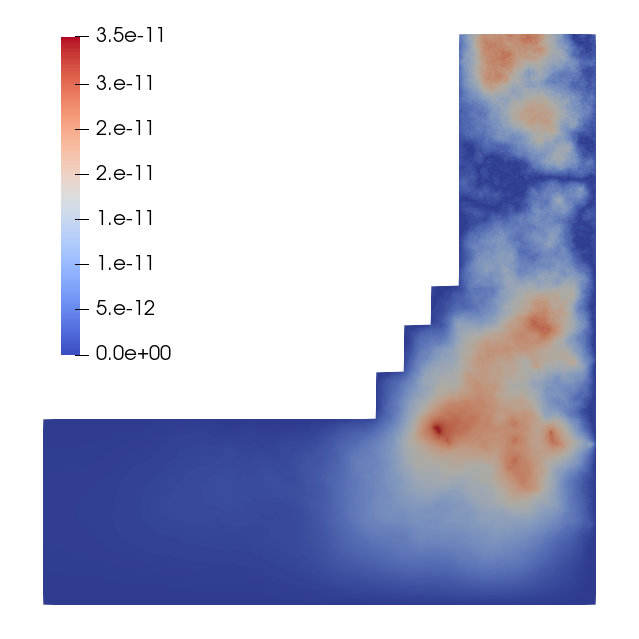}
%\caption*{$|T_h-T_a|$}
\end{minipage}
\caption{Benchmark for the thermal model ($WT$)$_h$: analytical temperature $T_a$ (left), numerical temperature $T_h$ (center), and corresponding absolute error $|T_a - T_h|$ (right) in K.}
\label{temperature_benchmark_comparison}
\end{figure}

We solve the $(WT)_h$ problem for the data $Q,q_+, T_f,T_{out}, T_-$ given by equations \eqref{thermal_benchmark_source_term}-\eqref{thermal_symmetry_boundary}. We choose a discretized space of polynomial of degree $3$. Analytical and numerical solutions are reported in Figure \ref{temperature_benchmark_comparison} (left and center). As we can see, a very good agreement is obtained. For a more quantitative comparison, we also display the absolute error in Figure \ref{temperature_benchmark_comparison} (right), and compute the relative error,
\begin{equation*}
\dfrac{||T_a - T_h||_{H^1_r(\omega)}}{||T_a||_{H^1_r(\omega)}} = 7e-13.
\end{equation*}

\begin{comment}
\begin{table}
\centering
\begin{tabular}{|c|c|}
\hline
Polynomial degree & Relative error\\
\hline
\hline
1 & 7.850e-04 \\
\hline
2 & 4.056e-07 \\
\hline
3 & 7.001e-13 \\
\hline
\end{tabular}
\hspace{1cm}
\begin{tabular}{|c|c|}
\hline
Step size & Relative error\\
\hline
\hline
0.25 & 0.0082 \\
\hline
0.2 & 0.0067 \\
\hline
0.15 & 0.0048 \\
\hline
0.10 & 0.0033 \\
\hline
\end{tabular}
\caption{Benchmark for the thermal problem:  convergence analysis. %\michele{These figures should be improved. The quality should be very higher.} \nirav{Revised} \michele{It is necessary to furthermore improve. Labels difficult to read, in the right figure the axis are missed.} \nirav{Removed figure and added table. It makes data clearly readable.}
}
\label{thermal_model_convergence}
\end{table}
\end{comment}

\begin{comment}
\begin{figure}[H]
\begin{subfigure}{.45\textwidth}
\includegraphics[width=\textwidth]{images/Convergence_temperature_benchmark_comparison.png}
\caption*{Relative error $p-$convergence}
\end{subfigure}
\hspace{1cm}
\begin{subfigure}{.45\textwidth}
\includegraphics[width=\textwidth]{images/energy_equation_h_convergence.png}
\caption*{Relative error $h-$convergence}
\end{subfigure}
\caption{Benchmark for the thermal problem:  convergence analysis. \michele{These figures should be improved. The quality should be very higher.} \nirav{Revised} \michele{It is necessary to furthermore improve. Labels difficult to read, in the right figure the axis are missed.}}
\label{thermal_model_convergence}
\end{figure}
\end{comment}

\subsection{Mechanical model}\label{mechanical_benchmark_section}

Firstly, we consider that the body is at reference temperature $T = T_0$, i.e. thermal stresses $(2 \mu + 3 \lambda) \alpha (T - T_0)$ are not present. Therefore, we refer to the problem ($WM1$)$_h$.% We postpone the discussion related to the thermal stresses to the next benchmark case.

We consider a known displacement function
\begin{equation}\label{benchmark_analytical_displacement}
\overrightarrow{u}_a = C (r y^2, r^2 y), \quad  \text{with} \ C = 1e-4 /m^2.
\end{equation}
%for mechanical problem $(M1)$ defined by equations \eqref{summarised_strong_momentum_equations}-\eqref{summarised_strong_momentum_equations_boundary}.
%In order that $\overrightarrow{u}_a$ is a solution of the mechanical model $(M1)$ (eqs. \eqref{summarised_strong_momentum_equations}-\eqref{summarised_strong_momentum_equations_boundary}),
Then:

\begin{itemize}

\item The components of the stress tensor $\bm{\sigma}$ are given by eq. \eqref{stress_cylindrical_coordinate} as, %\eqref{strain_tensor2} and \eqref{A_redefine}: % \michele{insert the eq. reference}: %= $\bm{\sigma}$($\overrightarrow{u}_a$) are given by: %associated with $\overrightarrow{u}_a$ according to \eqref{stress_strain_relation_axisymmetric} is:
\begin{subequations}\label{benchmark_mechanical_stress_tensor}
\begin{flalign}
\sigma_{rr} & = \frac{E}{(1-2\nu)(1+\nu)} (C y^2 + \nu C r^2) \ , \\
\sigma_{yy} & = \frac{E}{(1-2\nu)(1+\nu)} (2 \nu C y^2 + (1 - \nu) C r^2) \ , \\
\sigma_{\theta \theta} & = \frac{E}{(1-2\nu)(1+\nu)} (C y^2 + \nu C r^2) \ , \\
\sigma_{ry} & = \frac{2ECry}{(1+\nu)} \ .
\end{flalign}
\end{subequations}

\item The source term $\overrightarrow{f}_0 = [f_{0,r} \ f_{0,y}]$ is obtained from eqs. \eqref{summarised_strong_momentum_equations} and \eqref{benchmark_mechanical_stress_tensor} as,
\begin{subequations}
\begin{flalign}
f_{0,r} & = - \left( \frac{\partial \sigma_{rr}}{\partial r} + \frac{\partial \sigma_{ry}}{\partial y}  + \frac{\sigma_{rr} - \sigma_{\theta \theta}}{r} \right) \nonumber \\ & = - \left( \frac{2 E \nu C r}{(1 - 2 \nu)(1 + \nu)} + \frac{2 E C r}{(1+\nu)} \right) \ , \label{mechanical_benchmark_source_term_1} \\
f_{0,y} & = - \left( \frac{\partial \sigma_{ry}}{\partial r} + \frac{\partial \sigma_{yy}}{\partial y} + \frac{\sigma_{ry}}{r} \right) \nonumber \\ & = - \left( \frac{4 ECy}{(1+\nu)} + \frac{4 E \nu C y}{(1-2\nu)(1+\nu)} \right) \ . \label{mechanical_benchmark_source_term_2}
\end{flalign}
\end{subequations}

\item The boundary tractions are derived from eqs. \eqref{summarised_strong_momentum_equations_boundary} and \eqref{benchmark_mechanical_stress_tensor} as,
\begin{subequations}
\begin{flalign}
\textrm{on} \ \gamma_+ \ & : \ g_{+,r} = \frac{2ECry}{(1+\nu)} \ , \\ & \ \ \ g_{+,y} = \frac{E}{(1-2\nu)(1+\nu)} \left( 2 \nu C y^2 + (1-\nu) C r^2 \right) \ , \label{mechanical_benchmark_boundary_tractions_1} \\
\textrm{on} \ \gamma_- \ & : \ g_{-,r} = - \frac{2ECry}{(1+\nu)} \ , \label{mechanical_benchmark_boundary_tractions_2} \\
\textrm{on} \ \gamma_{sf} \ & : \ g_{sf,r} = \frac{E}{(1-2\nu)(1+\nu)} \left( Cy^2 + \nu C r^2 \right) n_r + \frac{2ECry}{(1+\nu)} n_y \ , \ \\ & \ \ \ g_{sf,y} = \frac{2ECry}{(1+\nu)} n_r + \frac{E}{(1-2\nu)(1+\nu)} \left( 2 \nu C y^2 + (1-\nu) C r^2 \right) n_y \ , \label{mechanical_benchmark_boundary_tractions_3} \\
\textrm{on} \ \gamma_{out} \ & : \ g_{out,r} = \frac{E}{(1-2\nu)(1+\nu)} \left( C y^2 + \nu C r^2 \right) \ , \\ & \ \ \ g_{out,y} = \frac{2ECry}{(1+\nu)} \ , \ \label{mechanical_benchmark_boundary_tractions_4}
\end{flalign}
\end{subequations}
and it is verified that
\begin{subequations}
\begin{flalign}
\textrm{on} \ \gamma_- \cup  \gamma_s \ & : \overrightarrow{u}_a \cdot \overrightarrow{n} = 0 \ . \ \label{mechanical_symmetry_boundary}
\end{flalign}
\end{subequations}

%$\overrightarrow{u}_a \cdot \overrightarrow{n} = 0$ on $\gamma_- \cup  \gamma_s$.

\end{itemize}

\begin{figure}[H]
\begin{minipage}{.3\textwidth}
\includegraphics[width=\textwidth]{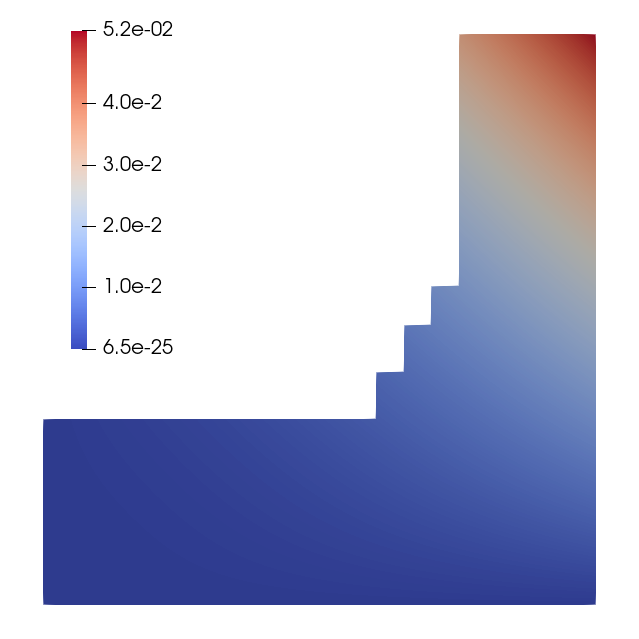}
%\caption*{$u_a$ (m)}
\end{minipage}
\begin{minipage}{.3\textwidth}
\includegraphics[width=\textwidth]{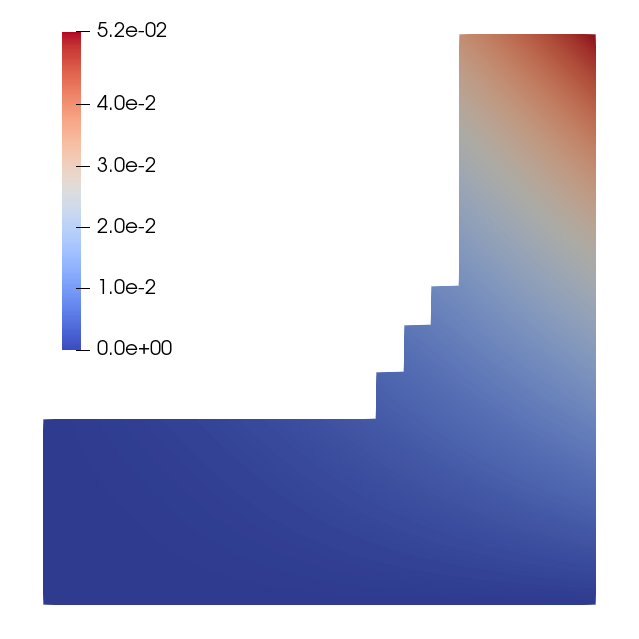}
%\caption*{$u_h$ (m)}
\end{minipage}
\begin{minipage}{.3\textwidth}
\includegraphics[width=\textwidth]{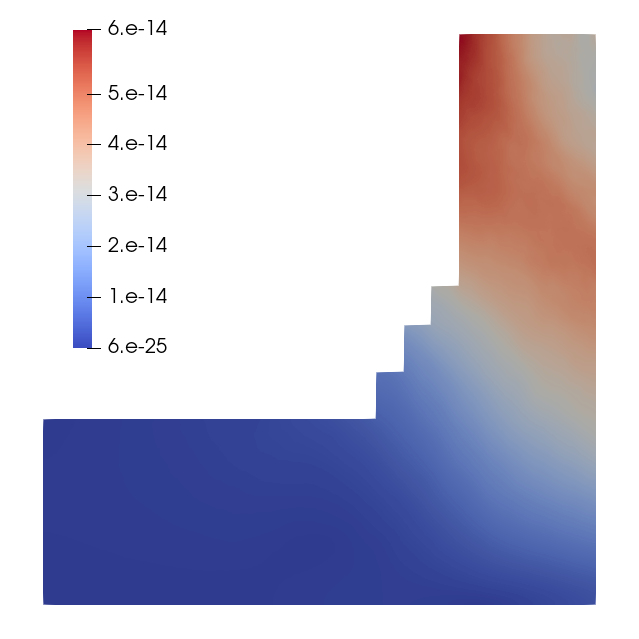}
%\caption*{$|u_h - u_a|$ (m)}
\end{minipage}
\caption{Benchmark for the mechanical model ($WM1$)$_h$: analytical displacement magnitude $|\protect\overrightarrow{u}_a|$ (left), numerical displacement magnitude $|\protect\overrightarrow{u}_h|$ (center), and absolute error magnitude $|\protect\overrightarrow{u}_a - \protect\overrightarrow{u}_h|$ (right) in m.} %\michele{Could you use exactly the same scale for the left and right images? Thanks.}\nirav{Done}}
\label{displacement_benchmark_comparison}
\end{figure}

We solve the $(WM1)_h$ problem for the data given by equations \eqref{mechanical_benchmark_source_term_1} - \eqref{mechanical_symmetry_boundary} by using a discretized space of polynomial of degree 3. The magnitude of the analytical and numerical displacement, as well as the associated absolute error, are represented in Figure \ref{displacement_benchmark_comparison}. Moreover, we compute the relative error

\begin{equation*}
\dfrac{||\overrightarrow{u}_a - \overrightarrow{u}_h||_{\mathbb{U}}}{||\overrightarrow{u}_a ||_{\mathbb{U}}} = 1.81e-12.
\end{equation*}

Like the thermal model, even in this case we could observe that the two solutions show a very good agreement. %Moreover, we obtain %a relative error $||\overrightarrow{u}_a - \overrightarrow{u}_h||_{\mathbb{U}} = 1.81e-12$ % \michele{This value is related to the Von Mises stress?}.
%Finally, we assess $p-$convergence and $h-$convergence (using a discretized space of polynomial of degree $1$) of the relative error in the $H^1_r$-norm in Figure \ref{thermal_model_convergence}.
%and assess $p-$convergence and $h$-convergence (using a discretized space of polynomial of degree $1$) with respect to $||\overrightarrow{u}_a - \overrightarrow{u}_h||_{\mathbb{U}}$ in Table %\ref{mechanical_model_convergence}.

\begin{figure}[H]
\begin{minipage}{.3\textwidth}
\includegraphics[width=\textwidth]{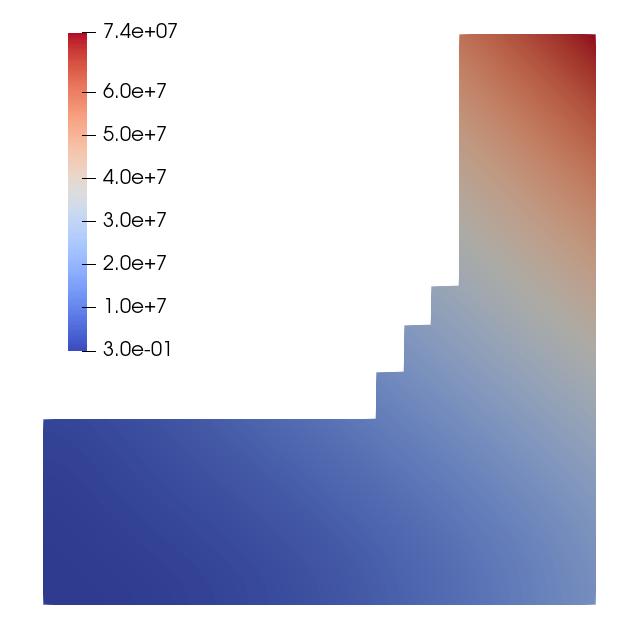}
%\caption*{${\sigma_{vm}}_h$ (N/m$^2$)}
\end{minipage}
\begin{minipage}{.3\textwidth}
\includegraphics[width=\textwidth]{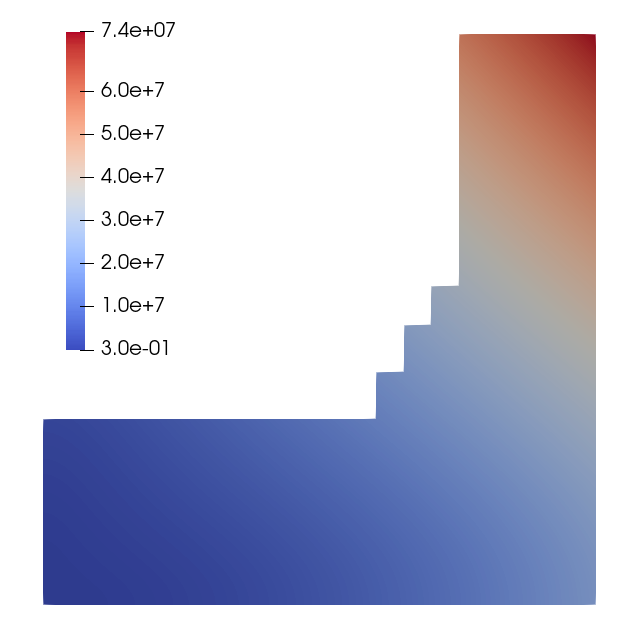}
%\caption*{${\sigma_{vm}}_a$ (N/m$^2$)}
\end{minipage}
\begin{minipage}{.3\textwidth}
\includegraphics[width=\textwidth]{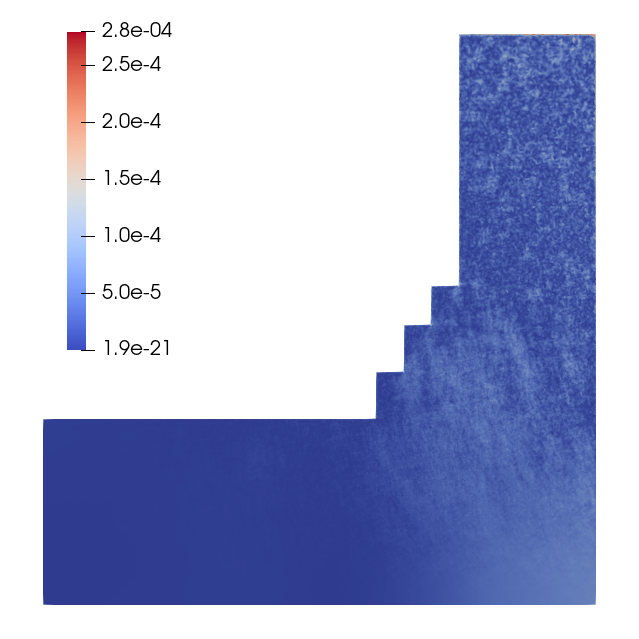}
%\caption*{$|{\sigma_{vm}}_h - {\sigma_{vm}}_a|$ (N/m$^2$)}
\end{minipage}
\caption{Benchmark for the mechanical model ($WM1$)$_h$: analytical Von Mises stress magnitude ${\sigma_{vm}}_a$ (left), numerical Von Mises stress magnitude ${\sigma_{vm}}_h$ (center), and absolute error $|{\sigma_{vm}}_h - {\sigma_{vm}}_a|$ (right) in N/m$^2$.} %\michele{Could you use exactly the same scale for the left and right images? Thanks}\nirav{Done}}
\label{displacement_benchmark_comparison_von_mises}
\end{figure}

For further comparison, we also computed the Von Mises stress:
\begin{equation}\label{von_mises_stress_eq_defn}
\sigma_{vm} = \sqrt{\frac{3}{2} \bm{\sigma}_d : \bm{\sigma}_d} \ ,
\end{equation}
%\end{itemize}
where $\bm{\sigma}_d$ is the deviatoric part of the stress tensor
\begin{equation}\label{deviatoric_stress_eq_defn}
\bm{\sigma}_d = \bm{\sigma} - \frac{1}{3} Tr(\bm{\sigma}) \bm{I} \ . %= \bm{\sigma} - \sigma_{hs} \bm{I} \, ,
\end{equation}
%and $\sigma_{hs}$ is the hydrostatic stress.
%\begin{itemize}
%\item Hydrostatic stress $\sigma_m$ :
%\begin{equation}\label{hydrostatic_stress_eq_defn}
%\sigma_m = \frac{1}{3}Tr(\bm{\sigma}) \ .
%\end{equation}
%\item Spherical part of the stress tensor $\bm{\sigma}_s$:
%\begin{equation}\label{spherical_stress_eq_defn}
%\bm{\sigma}_s = \frac{1}{3}Tr(\bm{\sigma}) \bm{I} = \sigma_m \bm{I} \ .
%\end{equation}
%\item Deviatoric part of the stress tensor $\bm{\The boundary tractions are derived fromsigma}_d$:
%\begin{equation}\label{deviatoric_stress_eq_defn}
%\bm{\sigma}_d = \bm{\sigma} - \frac{1}{3} Tr(\bm{\sigma}) \bm{I}  = \bm{\sigma} - \sigma_m \bm{I} \ .
%\end{equation}
%\item Von Mises effective stress $\sigma_{vm}$:
%\begin{equation}\label{von_mises_stress_eq_defn}
%\sigma_{vm} = \sqrt{\frac{3}{2} \bm{\sigma}_d : \bm{\sigma}_d} \ .
%\end{equation}
%\end{itemize}
%\end{comment}
%the Von Mises stress (Equation \eqref{von_mises_stress_eq_defn}) calculated from analytical displacement and the Von Mises stress calculated from computed displacement is given in Figure \ref{displacement_benchmark_comparison_von_mises}.
We display the magnitude of the analytical and numerical Von Mises stress, ${\sigma_{vm}}_a$ and ${\sigma_{vm}}_h$ respectively, and the corresponding absolute error in Figure \ref{displacement_benchmark_comparison_von_mises}. We see that a very good matching is obtained. %As we can see, the comparison %Notice that since the order of magnitude of Von Mises stress is high, the effect of round-off error in error plot becomes more pronounced \michele{with respect to what?}.

%Thus, $\overrightarrow{u}_a$ is the solution of the mechanical model $(M1)$ for the data given by equations \eqref{mechanical_benchmark_source_term_1} - \eqref{mechanical_benchmark_boundary_tractions_4} with $T=T­_0$. The comparison between the known analytical displacement with the computed displacement is given in Figure \ref{displacement_benchmark_comparison}. The comparison between the Von Mises stress (Equation \eqref{von_mises_stress_eq_defn}) calculated from analytical displacement and the Von Mises stress calculated from computed displacement is given in Figure \ref{displacement_benchmark_comparison_von_mises}. It can be further noted in Figure \ref{displacement_benchmark_comparison_von_mises} that since the order of magnitude of Von Mises stress is high, the effect of round-off error in error plot becomes more pronounced. For a more quantitative comparison, we computed the relative error in $\mathbb{U}$-norm (Equation \eqref{norms_mechanical_model}), that was $4.32e-12$. We also assess $p-$convergence behavior for the relative error in $\mathbb{U}-$norm and $h-$convergence behavior for the relative error in $\mathbb{U}-$norm using a discretized space of polynomial of degree $1$. The results are visualized in Figure \ref{mechanical_model_convergence}.

\begin{comment}
\begin{table}
\centering
\begin{tabular}{|c|c|}
\hline
Polynomial degree & Relative error\\
\hline
\hline
1 & 7.813e-04 \\
\hline
2 & 4.461e-07 \\
\hline
3 & 1.810e-12 \\
\hline
\end{tabular}
\hspace{1cm}
\begin{tabular}{|c|c|}
\hline
Step size & Relative error\\
\hline
\hline
0.25 & 0.0083 \\
\hline
0.2 & 0.0067 \\
\hline
0.15 & 0.0048 \\
\hline
0.10 & 0.0033 \\
\hline
\end{tabular}
\caption{Benchmark for the mechanical model: convergence analysis.}
\label{mechanical_model_convergence}
\end{table}
\end{comment}

\begin{comment}
\begin{figure}[H]
\begin{subfigure}{.45\textwidth}
\includegraphics[width=\textwidth]{images/Convergence_displacement_benchmark_comparison.png}
\caption*{Relative error $p-$convergence}
\end{subfigure}
\hspace{1cm}
\begin{subfigure}{.45\textwidth}
\includegraphics[width=\textwidth]{images/momentum_equation_h_convergence.png}
\caption*{Relative error $h-$convergence}
\end{subfigure}
\caption{Benchmark for the mechanical model: convergence analysis. \michele{These figures should be improved. The quality should be very higher.} \nirav{Revised} \michele{It is not sufficient.}}
\label{mechanical_model_convergence}
\end{figure}
\end{comment}

%\subsubsection{thermo-mechanical problem}

\begin{figure}[H]
\begin{minipage}{.3\textwidth}
\includegraphics[width=\textwidth]{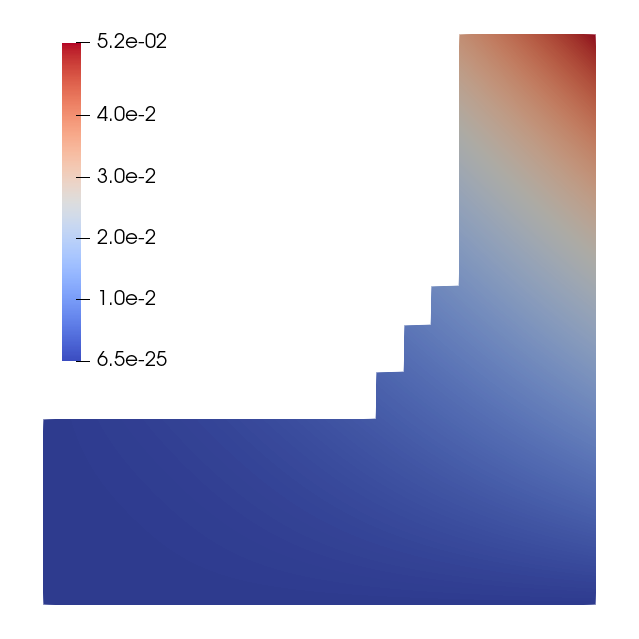}
%\caption*{$u_a$ (m)}
\end{minipage}
\begin{minipage}{.3\textwidth}
\includegraphics[width=\textwidth]{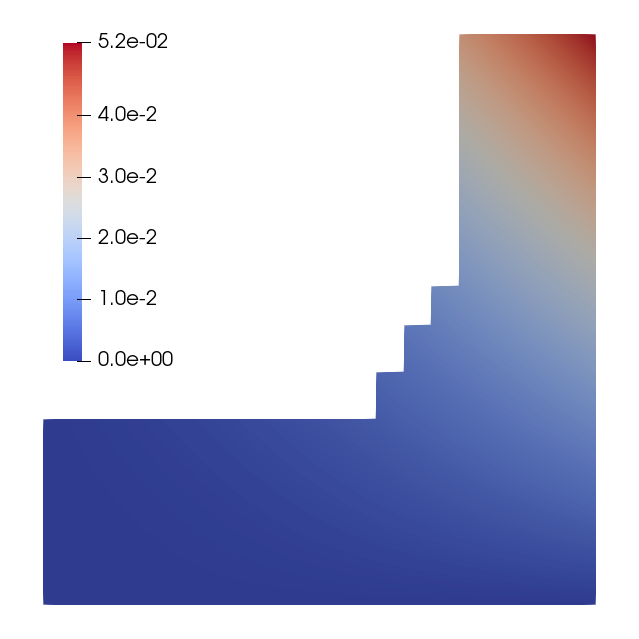}
%\caption*{$u_h$ (m)}
\end{minipage}
\begin{minipage}{.3\textwidth}
\includegraphics[width=\textwidth]{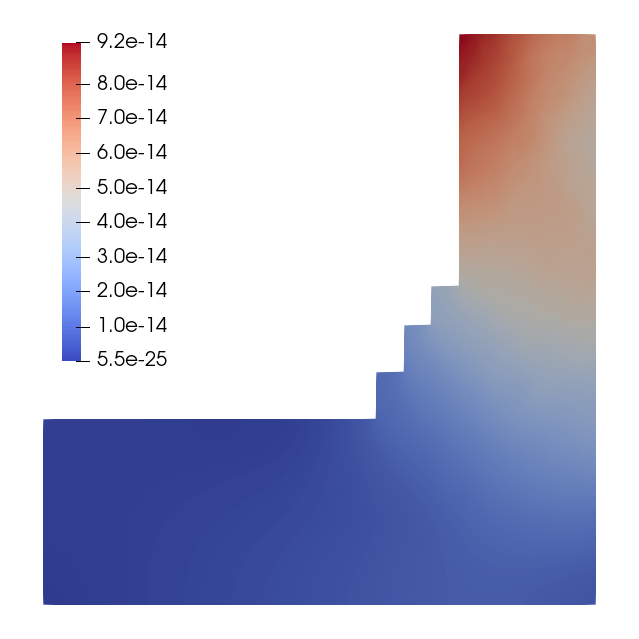}
%\caption*{$|u_h - u_a|$ (m)}
\end{minipage}
\caption{Benchmark for the mechanical problem ($WM$)$_h$: analytical displacement magnitude $|\protect\overrightarrow{u}_a|$ (left), numerical displacement magnitude $|\protect\overrightarrow{u}_h|$ (center), and absolute error $|\protect\overrightarrow{u}_a - \protect\overrightarrow{u}_h|$(right) in m.}
\label{coupling_benchmark_comparison_displacement}
\end{figure}

\begin{figure}[H]
\begin{minipage}{.3\textwidth}
\includegraphics[width=\textwidth]{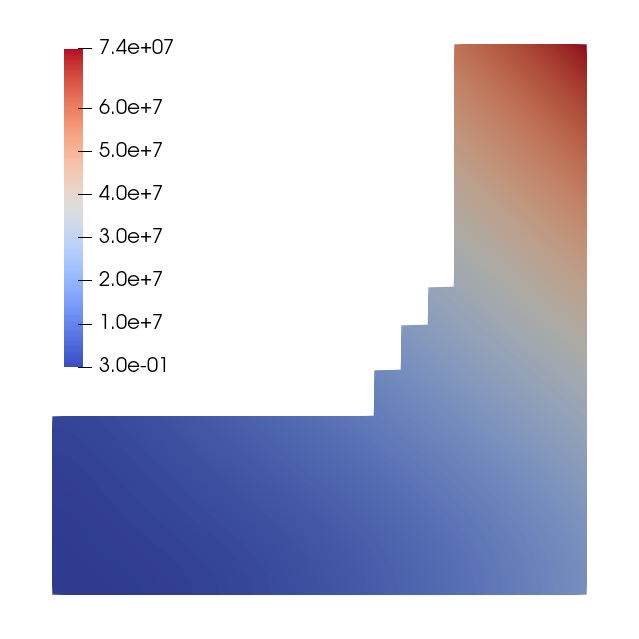}
%\caption*{${\sigma_{vm}}_h$ (N/m$^2$)}
\end{minipage}
\begin{minipage}{.3\textwidth}
\includegraphics[width=\textwidth]{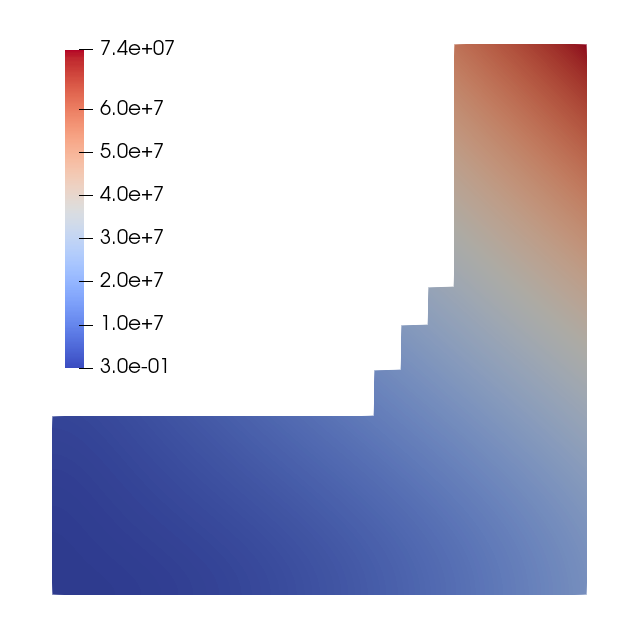}
%\caption*{${\sigma_{vm}}_a$ (N/m$^2$)}
\end{minipage}
\begin{minipage}{.3\textwidth}
\includegraphics[width=\textwidth]{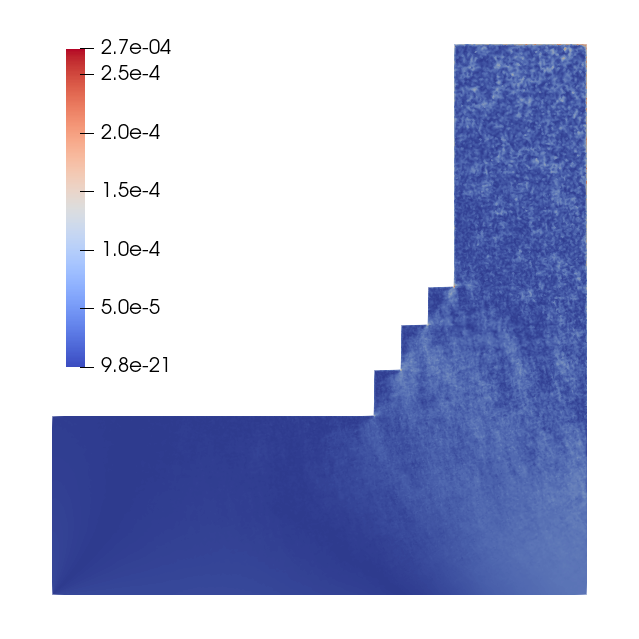}
%\caption*{$|{\sigma_{vm}}_h - {\sigma_{vm}}_a|$ (N/m$^2$)}
\end{minipage}
\caption{Benchmark for the mechanical problem ($WM$)$_h$: analytical Von Mises stress ${\sigma_{vm}}_a$ (left), numerical Von Mises stress ${\sigma_{vm}}_h$ (center) and absolute error $|{\sigma_{vm}}_h - {\sigma_{vm}}_a|$ (right) in N/m$^2$.}
\label{coupling_benchmark_comparison_von_mises}
\end{figure}

Now we address the coupling between the thermal and mechanical effects, so we refer to the problem ($WM$)$_h$. We assume for the temperature the analytical field used for the problem ($WT$)$_h$, $T_a$ (see eq. \eqref{benchmark_analytical_temperature}), and for the displacement the analytical field used for the problem ($WM1$)$_h$, $\overrightarrow{u}_a$ (see eq. \eqref{benchmark_analytical_displacement}). %We set $T_0 = 298 K$. %Firstly, we compute the hydrostatic thermo-mechanical stress (Equation \eqref{hydrostatic_stress_eq_defn}) when the domain is subjected to combined mechanical and thermal effects.
Then:

\begin{itemize}

\item We obtain the thermal stresses % \nirav{($(2 \mu + 3 \lambda) \alpha (T - T_0)\bm{I}$)
from eqs. \eqref{stress_tensor} and \eqref{lame_parameters}: % \michele{maybe it's better to insert an eq. ref.}: %\michele{I would write the following formula in a more clear way} \nirav{Done}
\begin{equation*}
(2 \mu + 3 \lambda) \alpha (T - T_0) = (2 \mu + 3 \lambda) \alpha (C'r^2y-T_0) = \frac{E}{(1-2\nu)} \alpha (C'r^2y-T_0) \ .
\end{equation*}

\item The components of the stress tensor $\bm{\sigma}$ are given by eqs. \eqref{stress_strain_relation_axisymmetric}: %\eqref{strain_tensor2} and \eqref{A_redefine}: % \michele{insert the eq.
\begin{subequations}\label{benchmark_coupling_stress_tensor}
\begin{flalign}
\sigma_{rr} & = \frac{E}{(1-2\nu)(1+\nu)} (C y^2 + \nu C r^2) - \frac{E}{(1-2\nu)} \alpha (C'r^2y-T_0) \ , \\
\sigma_{yy} & = \frac{E}{(1-2\nu)(1+\nu)} (2 \nu C y^2 + (1 - \nu) C r^2) - \frac{E}{(1-2\nu)} \alpha (C'r^2y-T_0) \ , \\
\sigma_{\theta \theta} & = \frac{E}{(1-2\nu)(1+\nu)} (C y^2 + \nu C r^2) - \frac{E}{(1-2\nu)} \alpha (C'r^2y - T_0) \ , \\
\sigma_{ry} & = \frac{2ECry}{(1+\nu)} \ .
\end{flalign}
\end{subequations}

\item The source term $\overrightarrow{f}_0 = [f_{0,r} \ f_{0,y}]$ is obtained from eqs. \eqref{summarised_strong_momentum_equations} and \eqref{benchmark_mechanical_stress_tensor} by:
\begin{subequations}
\begin{flalign}
f_{0,r} & = - \left( \frac{\partial \sigma_{rr}}{\partial r} + \frac{\partial \sigma_{ry}}{\partial y}  + \frac{\sigma_{rr} - \sigma_{\theta \theta}}{r} \right) \nonumber \\ & = - \left( \frac{2 E \nu C r}{(1 - 2 \nu)(1 + \nu)} + \frac{2 E C r}{(1+\nu)} - \frac{2C'ryE \alpha }{(1-2\nu)} \right) \ , \label{coupling_benchmark_source_term_1} \\
f_{0,y} & = - \left( \frac{\partial \sigma_{ry}}{\partial r} + \frac{\partial \sigma_{yy}}{\partial y} + \frac{\sigma_{ry}}{r} \right) \nonumber \\ &  = - \left( \frac{4 ECy}{(1+\nu)} + \frac{4 E \nu C y}{(1-2\nu)(1+\nu)} - \frac{C' r^2 E \alpha}{(1-2\nu)} \right) \ . \label{coupling_benchmark_source_term_2}
\end{flalign}
\end{subequations}

\item The boundary tractions are derived from eqs. \eqref{summarised_strong_momentum_equations_boundary} and \eqref{benchmark_mechanical_stress_tensor} as:
\begin{subequations}
\begin{flalign}
\textrm{on} \ \gamma_+ \ & : g_{+,r} = \frac{2ECry}{(1+\nu)} \ , \\
& \ \ g_{+,y} = \frac{E}{(1-2\nu)(1+\nu)} \left( 2\nu Cy^2 + (1-\nu) Cr^2 \right) \nonumber \\ & \ \ \ \ \ \ \ - \frac{E \alpha}{(1-2\nu)} (C' r^2 y - T_0) \ , \label{coupling_benchmark_boundary_tractions_1} \\
\textrm{on} \ \gamma_- \ & :  g_{-,r} = - \frac{2ECry}{(1+\nu)} \ , \label{coupling_benchmark_boundary_tractions_2} \\
\textrm{on} \ \gamma_{sf} \ & : \ g_{sf,r} = \frac{E}{(1-2\nu)(1+\nu)} (Cy^2 + \nu C r^2) n_r \nonumber \\ & \ \ \ \ \ \ \ \ \ - \frac{E \alpha}{(1-2\nu)}(C' r^2 y - T_0) n_r + \frac{2ECry}{(1+\nu)} n_y  \ , \ \\ & \ \ \ g_{sf,y} = \frac{E}{(1-2\nu)(1+\nu)} \left( 2\nu Cy^2 + (1-\nu) Cr^2 \right) n_y \nonumber \\ & \ \ \ \ \ \ \ \ \ \ - \frac{E \alpha}{(1-2\nu)} (C' r^2 y - T_0) n_y + \frac{2ECry}{(1+\nu)} n_r \ , \label{coupling_benchmark_boundary_tractions_3} \\
\textrm{on} \ \gamma_{out} \ & : \ g_{out,r} = \frac{E}{(1-2\nu)(1+\nu)} (Cy^2 + \nu C r^2) - \frac{E \alpha}{(1-2\nu)}(C' r^2 y - T_0) \ , \\ & \ \ \ g_{out,y} = \frac{2ECry}{(1+\nu)} \ . \label{coupling_benchmark_boundary_tractions_4}
\end{flalign}
\end{subequations}

\end{itemize}
%We display analytical ${\sigma_{vm}}_a$ and numerical ${\sigma_{vm}}_h$ Von Mises stress, and the corresponding absolute error in Figure \ref{displacement_benchmark_comparisreduced basison_von_mises}. %Notice that since the order of magnitude of Von Mises stress is high, the effect of round-off error in error plot becomes more pronounced \michele{with respect to what?}.
%The relative error in $\mathbb{U}$-norm (Equation \eqref{norms_mechanical_model})  was $4.32e-12$.
We display the magnitude of the analytical displacement and Von Mises stress comparing them with the corresponding numerical values in Figures \ref{coupling_benchmark_comparison_displacement} and \ref{coupling_benchmark_comparison_von_mises}, respectively, and compute the relative error

\begin{equation*}
\dfrac{||\overrightarrow{u}_a - \overrightarrow{u}_h||_{\mathbb{U}}}{||\overrightarrow{u}_a||_{\mathbb{U}}} = 2.2e-12.
\end{equation*}

We could see that, as for mechanical model ($WM1$)$_h$, the agreement between the two solutions is very good. %Then, we compute the relative error $||\overrightarrow{u}_a - \overrightarrow{u}_h||_{\mathbb{U}} = 2.2e-12$. %and report related $p-$convergence and $h-$convergence (using an discretized space of polynomial of degree 1) in Table \ref{coupling_model_convergence}.

\begin{comment}
\begin{table}
\centering
\begin{tabular}{|c|c|}
\hline
Polynomial degree & Relative error\\
\hline
\hline
1 & 7.814e-04 \\
\hline
2 & 4.461e-07 \\
\hline
3 & 2.237e-12 \\
\hline
\end{tabular}
\hspace{1cm}
\begin{tabular}{|c|c|}
\hline
Step size & Relative error\\
\hline
\hline
0.25 & 0.0091 \\
\hline
0.2 & 0.0074 \\
\hline
0.15 & 0.0054 \\
\hline
0.10 & 0.0037 \\
\hline
\end{tabular}
\caption{Benchmark for the thermo-mechanical problem: convergence analysis.}
\label{coupling_model_convergence}
\end{table}
\end{comment}

\begin{figure}[H]
\begin{minipage}{.3\textwidth}
\includegraphics[width=\textwidth]{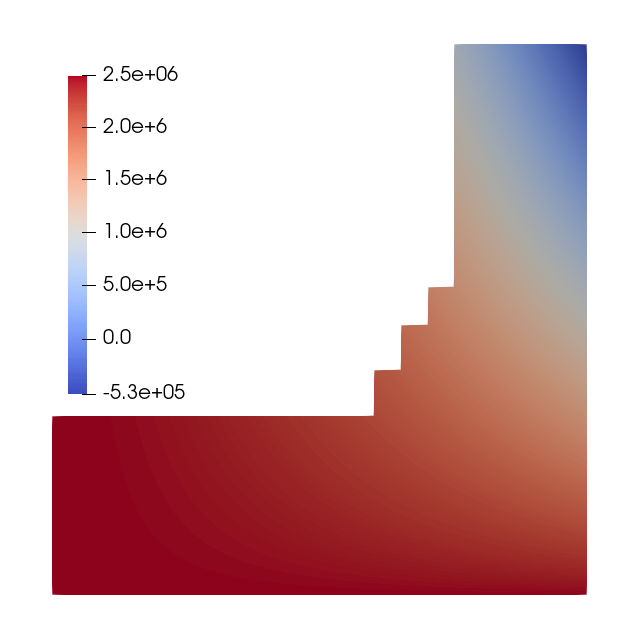}
%\caption*{$\frac{1}{3}tr\left(\bm{\sigma}(\overrightarrow{u})[T]-\bm{\sigma}(\overrightarrow{u})[T_0]\right)$ (N/m$^2$)}
\end{minipage}
\begin{minipage}{.3\textwidth}
\includegraphics[width=\textwidth]{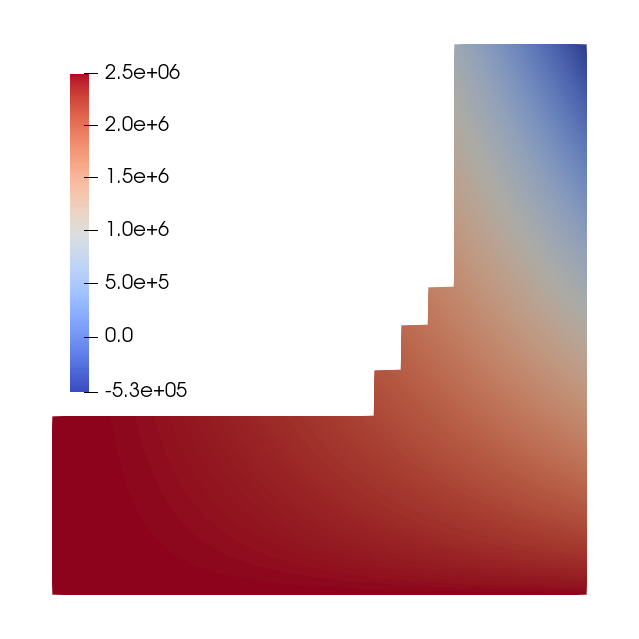}
%\caption*{$(2 \mu + 3 \lambda) \alpha (T - T_0)$ (N/m$^2$)}
\end{minipage}
\begin{minipage}{.3\textwidth}
\includegraphics[width=\textwidth]{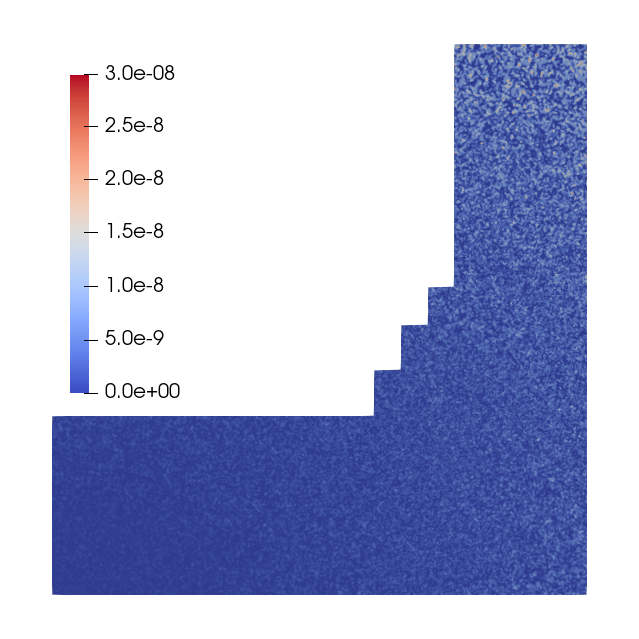}
%\caption*{Absolute error (N/m$^2$)}
\end{minipage}
\caption{Benchmark for the mechanical model ($WM2$)$_h$: difference between the hydrostatic stress computed with the $(WM)_h$ model and the $(WM1)_h$ model $\frac{1}{3}tr\left({\sigma}(\protect\overrightarrow{u_h})[T_h]-{\sigma}(\protect\overrightarrow{u_h})[T_0]\right)$ (left), thermal stress $(2 \mu + 3 \lambda) \alpha (T_h - T_0)$ (center) and corresponding absolute error $\left|\frac{1}{3}tr\left({\sigma}(\protect\overrightarrow{u_h})[T_h]-{\sigma}(\protect\overrightarrow{u_h})[T_0]\right) - (2 \mu + 3 \lambda) \alpha (T_h - T_0)\right|$ (right) in N/m$^2$.}
\label{coupling_benchmark_comparison_stress}
\end{figure}

Finally, we observe that the difference between the hydrostatic stress %$\sigma_{hs}$ (defined in eq. \eqref{deviatoric_stress_eq_defn}) %\michele{(why do you consider only the hydrostatic stress?)} \nirav{(because thermal stress is a hydrostatic stress)} (see eq. \ref{deviatoric_stress_eq_defn})
computed with the model ($WM$)$_h$ and the one computed with the model ($WM1$)$_h$, i.e. the hydrostatic stress related to the model ($WM2$)$_h$, should be equal to the thermal stress:
\begin{equation}\label{eq:hydr2}
\frac{1}{3}tr\left(\bm{\sigma}(\overrightarrow{u_h})[T_h]-\bm{\sigma}(\overrightarrow{u_h})[T_0]\right) = (2 \mu + 3 \lambda) \alpha (T_h - T_0) \bm{I}.
\end{equation}
The right and hand sides of eq. \eqref{eq:hydr2} are shown in Figure \ref{coupling_benchmark_comparison_stress}. We obtain a very good agreement. %It can be further noted in figures \ref{coupling_benchmark_comparison_von_mises} and \ref{coupling_benchmark_comparison_stress} that since the order of magnitude of Von Mises stress and Hydrostatic stress is high the effect of round-off error in error plot becomes more pronounced.
%Finally, we computed the relative error in $\mathbb{U}$-norm (Equation \eqref{norms_mechanical_model}), that was $2.237e-12$ \michele{related to???} \nirav{Related to analytical displcement}. We also assess $p-$convergence behavior and $h-$convergence behavior of the computed displacement field, using an discretized space of polynomial of degree 1, with the analytical one in Figure \ref{coupling_model_convergence}.

\begin{comment}
\begin{figure}[H]
\begin{subfigure}{.45\textwidth}
\includegraphics[width=\textwidth]{images/Convergence_coupling_displacement_benchmark_comparison.png}
\caption*{Relative error $p-$convergence}
\end{subfigure}
\hspace{1cm}
\begin{subfigure}{.45\textwidth}
\includegraphics[width=\textwidth]{images/coupling_h_convergence.png}
\caption*{Relative error $h-$convergence}
\end{subfigure}
\caption{thermo-mechanical model: convergence analysis. \michele{These figures should be improved. The quality should be very higher.} \nirav{Revised} \michele{it is not sufficient}}
\label{coupling_model_convergence}
\end{figure}
\end{comment}%reduced basis

\bibliographystyle{elsarticle-harv}
\bibliography{bibliography}

%\end{thebibliography}
\end{document}